\documentstyle{article}

\newcommand{\N}{{\rm I\kern-.5ex N}}
\newcommand{\Z}{{\sf \vrule height 1.55ex depth-1.2ex width.03em\kern-.11em Z 
\kern-.9ex Z\kern-.11em\vrule height 0.3ex depth0ex width.03em}}
\newcommand{\Q}{{\rm\kern.2ex\vrule height1.55ex depth-.05ex width.03em\kern-
.7ex Q}}
\newcommand{\R}{{\rm I\kern-.5ex R}}
\newcommand{\Rvar}{{\rm I\kern-.5ex R}}
\newcommand{\C}{{\rm\kern.3ex\vrule height1.55ex depth-.05ex width.03em\kern-
.7ex C}}

\newcommand{\spat}{\hspace{4ex}}

\newcommand{\cst}{C$^*$}

\newcommand{\otp}{\dot{\otimes}}
\newcommand{\cL}{{\cal L}}
\newcommand{\cF}{{\cal F}}
\newcommand{\cG}{{\cal G}}
\newcommand{\cH}{{\cal H}}
\newcommand{\cN}{{\cal N}}
\newcommand{\cM}{{\cal M}}
\newcommand{\tNfi}{\tilde{\cN}_\vfi}

\newcommand{\od}{\odot}
\newcommand{\ot}{\otimes}

\newcommand{\la}{\Lambda}

\newcommand{\om}{\omega}
\newcommand{\io}{\iota}
\newcommand{\vfi}{\varphi}
\newcommand{\vep}{\varepsilon}

\newcommand{\sde}{\delta}
\newcommand{\de}{\Delta}

\newcommand{\th}{\theta}

\newcommand{\si}{\sigma}
\newcommand{\Mfi}{{\cal M}_{\vfi}}
\newcommand{\Nfi}{{\cal N}_{\vfi}}

\newcommand{\lafi}{\la_\vfi}

\newcommand{\pifi}{\pi_\vfi}

\newcommand{\text}[1]{\mbox{#1}}

\newcommand{\qed}{\ \hfill \rule{2mm}{2mm}}
\newenvironment{demo}{\medskip\noindent\bf Proof :\ \  \rm}{\qed\bigskip\par }

\newtheorem{definition}{Definition}[section]
\newtheorem{proposition}[definition]{Proposition}
\newtheorem{lemma}[definition]{Lemma}
\newtheorem{corollary}[definition]{Corollary}
\newtheorem{remark}[definition]{Remark}
\newtheorem{theorem}[definition]{Theorem}

\newtheorem{notation}[definition]{Notation}
\newtheorem{result}[definition]{Result}

\begin{document}
\begin{center}
\Large\bf
Regular $\protect\boldmath C^*$-valued weights.

\end{center}

\begin{center}
\rm J. Kustermans  \footnote{Research Assistant of the
National Fund for Scientific Research (Belgium)}

Departement Wiskunde

Katholieke Universiteit Leuven

Celestijnenlaan 200B

B-3001 Heverlee

Belgium
\bigskip

January 1997
\end{center}

\subsection*{Abstract}
We introduce the notion of a \cst-valued weight between two \cst-algebras  as 
a generalization of an ordinary weight on a \cst-algebra and as a \cst-
version of operator valued weights on von Neumann algebras.

Also, some form of lower semi-continuity will be discussed together with
an extension to the Multiplier algebra.

A strong but useful condition for \cst-valued weights, the so-called 
regularity, is introduced. At the same time, we propose a construction
procedure for such regular \cst-valued weights.

This construction procedure will be used to define the tensor product of 
regular \cst-valued weights.

\section*{Introduction.}

In the articles \cite{Haa1} and \cite{Haa2}, Uffe Haagerup introduced the 
notion of operator valued weights on von Neumann algebras as generalizations 
of weights which can take values in another von Neumann algebra.

\medskip

An application of this theory can be found in the theory of the Kac-algebras 
$(M,\de,\vfi)$ (see \cite{E-S}). There, we have a semi-finite faithful normal 
weight $\vfi$ on the von Neumann algebra $M$ which is left invariant
in the sense that $(\io \overline{\ot} \vfi)(\de(a)) =  \vfi(a) \, 1$ for 
every $a \in M^+$. But one has to give a meaning to the expression $(\io 
\overline{\ot} \vfi)(\de(a))$. This can be done, using the theory of operator 
valued weights. In this case, $\io \overline{\ot} \vfi$ is an operator valued 
weight from $M \overline{\ot} M$ into $M$.

\medskip

Now we turn to the case of quantum groups in the framework of \cst-algebras 
in stead of in the framework of von Neumann algebras.
In this case, we have a \cst-algebra $A$, a comultiplication $\de$ from $A$ 
into the Multiplier algebra $M(A \ot A)$ and a densely defined lower semi-
continuous weight $\vfi$ on $A$ which is left invariant in the sense
that $(\io \ot \vfi)(\de(a)) = \vfi(a) \, 1$ for every $a \in \Mfi$.

Also in this case we have to give a meaning to the expression
$(\io \ot \vfi)(\de(a))$ for $a \in \Mfi$. This can be done by regarding
$\io \ot \vfi$ as a \cst-valued weight from $A \ot A$ into $A$ (which can be 
extended to certain elements of $M(A \ot A)$).

We want to perform the construction of $\io \ot \vfi$ in a purely \cst-
algebraic setting, not in a von Neumann algebra setting. For instance, we
do not want to extend $\vfi$ to a normal weight on a bigger von Neumann 
algebra and work with this extension. Partly, because there are no guarantees 
that $\vfi$ can be extended in the case of \cst-algebraic quantum groups.

If we look at the left invariance of $\vfi$, we see that $(\io \ot 
\vfi)(\de(a))$ has to belong to $M(A)$ for every $a \in \Mfi$.
But also in the case of quantum groups there is some interest to let
\cst-valued weights take values in the set of affiliated elements. This can 
be seen as follows.

A \cst-algebraic quantum group will (possibly) have the analogue of a modular 
function for classical groups. This will be a strictly positive element 
$\sde$  affiliated with $A$ such that
$(\vfi \ot \io)(\de(a)) = \sde \, \vfi(a)$ for every $a \in \Mfi^+$.
We have to give a meaning to the expression $(\vfi \ot \io)(\de(a))$
but we see already that it will be unbounded in many cases.

\medskip

In the first section, we will give a possible definition for \cst-valued
weights. Loosely speaking, they will be completely positive linear mappings 
between \cst-algebras which are unbounded. We will restrict
the domain of our \cst-valued weights and let them take values
in the Multiplier algebra of another \cst-algebra. Later, we will discuss
some extensions and even let them take values in the set of affiliated 
elements.

We will also introduce a KSGNS-construction of a \cst-valued weight,
modelled on the KSGNS-construc- tion for completely positive mappings,
(for instance, see \cite{Lan}).

\medskip

In the second section we introduce a special family of completely positive 
mappings relative to a \cst-valued weight (and even a little bit more 
general). These completely positive mappings allow us to introduce a notion 
of lower semi-continuity  for \cst-valued weight in the third section. We 
also introduce the notion of regular \cst-valued weights in this third 
section.

\medskip

In the fourth section, the  extension of a lower semi-continuous weight to 
the Multiplier algebra is discussed.

\medskip

The fifth section serves as a first step in a construction procedure for \cst-
valued weights.

\medskip

We propose a construction procedure for regular \cst-valued weights
in the sixth section. Along the way, we prove an important result
which will be used in the next section to prove some nice results
about regular \cst-valued weights. There is also a short discussion
about a further extension of regular \cst-valued weights which lets it take 
values in the set of affiliated elements.

\medskip

In the last section, we introduce the tensorproduct of two regular \cst-
valued weights using the construction procedure of section 6.

\bigskip

At the end of this introduction, we will fix some notations and conventions.

\medskip

The main technical tools for this paper come from the theory of Hilbert \cst-
modules over \cst-algebras. A nice overview of this theory can be found in 
\cite{Lan}.

\medskip

\medskip

All our Hilbert \cst-modules are right modules and have innerproducts
which are linear in the first variable.

Consider Hilber \cst-modules $E,F$ over a \cst-algebra $A$. We will use
the following notations:
\begin{itemize}
\item $L(E,F)$ is the set of linear mappings from $E$ into $F$.
\item ${\cal B}(E,F)$ is the set of bounded linear mappings from $E$ into $F$.
\item $\cL(E,F)$ is the set of adjointable mappings from $E$ into $F$.
\end{itemize}

\medskip\medskip

Consider two \cst-algebras $A$ and $B$ and $\rho$ a completely positive 
mapping from $A$ into $M(B)$. We call $\rho$ strict if and only if is 
strictly continuous on bounded sets.

Let $(e_k)_{k \in K}$ be an approximate unit for $A$. Theorem 6.5 of 
\cite{Lan} implies that $\rho$ is strict if and only if
$(\rho(e_k)\,)_{k \in K}$ is strictly convergent.

If $\rho$ is strict, than it has a unique extension $\overline{\rho}$
which is a completely positive linear mapping, strictly continuous on
bounded sets. For every $a \in M(A)$, we define $\rho(a) = 
\overline{\rho}(a)$.

\section{The definition of a \cst-valued weight and its KSGNS-construc- 
tion.} \label{art1}

In this first section, we will introduce the definition of a \cst-valued 
weight between \cst-algebras. This definition is a generalization of the 
definition of usual weights on \cst-algebras but we will assume
complete positivity in stead of just positivity (which is to be expected). 
After that, we will introduce a KSGNS-construction for such a \cst-valued 
weight similar to the KSGNS-construction for completely positive mappings.

\medskip

We will start of with the definition.

\begin{definition}
Consider two C$^*$-algebras $A$ and $B$ and a hereditary cone $P$ in
$A^+$. \newline
Put ${\cal N} = \{ a \in A \mid a^* a \in P \}$
and ${\cal M} = \text{span } P = {\cal N}^*{\cal N}$.
Suppose that $\vfi$ is a linear mapping from ${\cal M}$ into $M(B)$ such that
$$\sum_{i,j=1}^n b_j^* \, \vfi(a_j^* \, a_i)\, b_i \geq 0 $$
for every $n \in \N$ and all
$a_1,\ldots\!,a_n \in  {\cal N}$, $b_1,\ldots\!,b_n \in B$.
Then we call $\vfi$ a C$^*$-valued weight from $A$ into $M(B)$.
\end{definition}

We will introduce the following notations.

\begin{notation}
Consider two \cst-algebras $A$ and $B$ and a \cst-valued weight $\vfi$ from 
$A$ into $M(B)$. We will use the following notations:
\begin{itemize}
\item The domain of $\vfi$ will be  denoted by $\Mfi$.
\item We define $\Nfi = \{ a \in A \mid a^* a \in \Mfi^+ \}$.
\end{itemize}
We have that $\Mfi$ is a sub-$^*$-algebra of $A$, $\Mfi^+$ is a hereditary 
cone in $A^+$ and $\Mfi = \text{span } \Mfi^+$. \newline
Furthermore, $\Nfi$ is a left ideal in $M(A)$ and $\Mfi = \Nfi^* \Nfi$.
\end{notation}

\begin{remark}
Consider two \cst-algebras $A,B$ and a \cst-valued weight from $A$ into 
$M(B)$. As usual, we say that $\vfi$ is densely defined $\Leftrightarrow$ \ 
$\Mfi$ is dense in $A$ \ $\Leftrightarrow$ \ $\Nfi$ is dense in $A$ 
$\Leftrightarrow$ \ $\Mfi^+$ is dense in $A^+$.
\end{remark}

As the generalization of the GNS-construction for weights, we have the
KSGNS-construction for \cst-valued weights:

\begin{definition}
Consider two C$^*$-algebras $A$ and $B$ and a C$^*$-valued weight
$\vfi$ from $A$ into $M(B)$.

A KSGNS-construction for $\vfi$ is by definition
a triplet $(E,\la,\pi)$ where
\begin{itemize}
\item $E$ is a Hilbert C$^*$-module over $B$.
\item $\la$ is a linear mapping from $\Nfi$ into ${\cal L}(B,E)$ such that
\begin{enumerate}
\item The set $\langle \, \la(a) b \mid a \in \Nfi , b \in B \, \rangle$ is 
dense in $E$
\item We have for every $a_1,a_2 \in \Nfi$ and $b_1,b_2 \in B$ that
$\langle \la(a_1) b_1 , \la(a_2) b_2 \rangle = b_2^* \,  \vfi(a_2^* \,  a_1) 
\, b_1$.
\end{enumerate}
\item $\pi$ is a $^*$-homomorphism from $A$ into ${\cal L}(E)$ such that
$\pi(x) \la(a) = \la(x a)$ for every $x \in A$ and $a \in \Nfi$.
\end{itemize}
\end{definition}

It is clear that a KSGNS-construction is unique up to unitary equivalence.

\begin{result} \label{res1.1}
Consider two C$^*$-algebras $A$ and $B$ and a C$^*$-valued weight
$\vfi$ from $A$ into $M(B)$ and let $(E,\la,\pi)$ be a KSGNS-construction
for $\vfi$. Then
\begin{enumerate}
\item We have that $\vfi(a_2^* \, a_1) = \la(a_2)^* \la(a_1)$ for every
$a_1,a_2 \in \Nfi$.
\item We have that $\|\la(a)\|^2 = \|\vfi(a^* a)\|$ for every $a \in \Nfi$.
\item If $\pi$ is non-degenerate, then
$\pi(x) \la(a) = \la(x a)$ for every $x \in M(A)$ and $a \in \Nfi$.
\end{enumerate}
\end{result}

\begin{demo}
The first property follows immediately from the definition of the KSGNS-
construction. The second follows from the first.
For the proof of the third, suppose that $\pi$ is non-degenerate.

Take $x \in M(A)$ and $a \in \Nfi$.
We have for every $e \in A$ that $e x$ belongs to $A$, thus
$$\pi(e) (\pi(x) \la(a)) = \pi(e x) \la(a) = \la(e x a) = \pi(e) \la(x a)  .$$
From the non-degeneracy of $\pi$, we infer that $\pi(x) \la(a) = \la(x a)$.
\end{demo}

\bigskip

The proof of the existence of a KSGNS-construction is a generalization of the 
KSGNS-construction for completely positive maps.
For a large part, we will mimic the proof of theorem 5.6 of \cite{Lan} but we 
will also add some features for this specific case (lemma \ref{lem1.2} and 
definition \ref{def1.1}).

\bigskip

For the most part of this section, we will fix two C$^*$-algebras $A$ and $B$ 
and a C$^*$-valued weight $\vfi$ from $A$ into $M(B)$. In the next part, we 
will gradually construct a KSGNS-construction for $\vfi$.

\medskip

First, we define the complex vector space $F = \Nfi \od B$. We turn $F$ into 
a semi innerproduct module over $B$ such that
\begin{itemize}
\item We have for every $a \in \Nfi$, $b,c \in B$ that $(a \ot b) \, c =
a \ot (b c)$.
\item We have for every $a_1,a_2 \in \Nfi$, $b_1,b_2 \in B$ that
$\langle a_1 \ot b_1 , a_2 \ot b_2 \rangle = b_2^* \, \vfi(a_2^* \, a_1) \, 
b_1$.
\end{itemize}
We put $N = \{ x \in F \mid \langle x , x \rangle = 0 \}$, then $N$ is a 
submodule of $F$. By the discussion in chapter 1 of \cite{Lan}, we know that
$\frac{F}{N}$ can be naturally turned into a innerproduct module over $B$. We 
define $A \ot_\vfi B$ to be the completion of $\frac{F}{N}$.
So $A \ot_\vfi B$ is a Hilbert \cst-module over $B$.

For every $a \in \Nfi$ and $b \in B$, we define $a \otp b$ to be the
equivalence class in $\frac{F}{N}$ associated with $a \ot b$.

Then we have the following properties:
\begin{itemize}
\item The mapping $\Nfi \times B \rightarrow A \ot_\vfi B : (a,b) \mapsto a 
\otp b$ is bilinear.
\item For every $a \in \Nfi$ and $b,c \in B$, we have that
$(a \otp b) \, c = a \otp (b c)$.
\item For every $a_1,a_2 \in \Nfi$, $b_1,b_2 \in B$, we have that
$\langle a_1 \otp b_1 , a_2 \otp b_2 \rangle = b_2^* \, \vfi(a_2^* \,  a_1) 
\, b_1$.
\item The set $\langle a \otp b \mid a \in \Nfi , b \in B \rangle$ is dense 
in $A \ot_\vfi B$.
\end{itemize}
These properties determine the Hilbert \cst-module $A \ot_\vfi B$ completely.

\medskip

It is easy to see that $\| a \otp b \| \leq \| \vfi(a^* a) \|^\frac{1}{2} \, 
\|b\|$ for every $a \in \Nfi$ and $b \in B$.
Therefore, we have for every $a \in \Nfi$ that the mapping
$B \rightarrow A \ot_\vfi B : b \mapsto a \otp b$ is continuous.

\medskip

First, we prove the following lemma:

\begin{lemma} \label{lem1.2}
Let $a$ be an element in $\Nfi$ and define the mapping $t$ from $B$ into $A 
\ot_\vfi B$ such that $t(b) = a \otp b$ for every $b \in B$.
Then $t$ belongs to $\cL(B,A \ot_\vfi B)$ and $t^*(c \otp d) = \vfi(a^* c) d$ 
for every $c \in \Nfi$, $d \in B$.
\end{lemma}
\begin{demo}
We have for every $c \in \Nfi$ and $b,d \in B$ that
$$\langle t(b) , c \otp d \rangle = \langle a \otp b , c \otp d \rangle
= d^*  \vfi(c^* a) b  . $$
Therefore proposition \ref{prop4.1} implies that $t$ is an element of 
$\cL(B,A \ot_\vfi B)$. It is also clear from the above equation that
$t^*(c \otp d) = \vfi(a^* c) d$ for every $c \in \Nfi$ and $d \in B$.
\end{demo}

This lemma justifies the following definition.

\begin{definition} \label{def1.1}
We define the linear mapping $\lafi$ from $\Nfi$ into $\cL(B, A \ot_\vfi B)$ 
such that $\lafi(a) b = a \otp b$ for every $a \in \Nfi$ and $b \in B$. We 
also have that $\lafi(a)^* (c \otp d) = \vfi(a^* c) d$ for every
$a,c \in \Nfi$ and $d \in B$.
\end{definition}

Using this definition, it is straightforward to check the following 
properties:
\begin{itemize}
\item The set $\langle \, \lafi(a) b \mid a \in \Nfi, b \in B \rangle$ is 
dense in $A \ot_\vfi B$.
\item We have for every $a_1,a_2 \in \Nfi$ and $b_1,b_2 \in B$ that
$\langle \lafi(a_1) b_1 , \lafi(a_2) b_2 \rangle = b_2^* \, \vfi(a_2^* \, 
a_1) \, b_1$.
\end{itemize}

\medskip

\begin{lemma}
Consider $a_1,\ldots\!,a_n \in \Nfi$, $b_1,\ldots\!,b_n \in B$ and
$x \in M(A)$. Then
$$\sum_{i,j=1}^n b_j^* \, \vfi((x a_j)^* (x a_i)) \, b_i
\leq \|x\|^2 \,\, \sum_{i,j=1}^n b_j^* \, \vfi(a_j^* \, a_i) \, b_i .$$
\end{lemma}
\begin{demo}
There exist an element $y \in M(A)$ such that $\|x\|^2 \, 1 - x^* x = y^* y$. 
Because $\Nfi$ is a left ideal in $M(A)$, we have that $y a_1,\ldots\!,y a_n$ 
belong to $\Nfi$.

This implies that
\begin{eqnarray*}
& & \|x\|^2 \,\, \sum_{i,j=1}^n b_j^* \, \vfi(a_j^* \,  a_i) \, b_i
- \sum_{i,j=1}^n b_j^* \, \vfi((x a_j)^* (x a_i)) \, b_i
= \sum_{i,j=1}^n b_j^* \, \vfi(a_j^* (\|x\|^2 \, 1 - x^* x) a_i) \, b_i \\
& & \ \ \ \ =  \sum_{i,j=1}^n b_j^* \, \vfi(a_j^* y^* y a_i) \,  b_i
 =  \sum_{i,j=1}^n b_j^* \, \vfi((y a_j)^* (y a_i)) \, b_i \geq 0
\end{eqnarray*}
\end{demo}

\begin{lemma}
Consider $x \in M(A)$. Then there exists a unique $T \in \cL(A \ot_\vfi B)$ 
such that $T \lafi(a) = \lafi(x a)$ for every $a \in \Nfi$.
Moreover, $T^* \lafi(a) = \lafi(x^* a)$ for every $a \in \Nfi$.
\end{lemma}
\begin{demo}
Choose $y \in M(A)$.

Take $a_1,\ldots\!,a_n \in \Nfi$ and $b_1,\ldots\!,b_n \in B$. We can
restate the previous lemma in the following form:
$$ \langle \sum_{i=1}^n \lafi(y a_i) b_i , \sum_{i=1}^n \lafi(y a_i) b_i 
\rangle \leq \|y\|^2 \,\, \langle \sum_{i=1}^n \lafi( a_i) b_i , \sum_{i=1}^n 
\lafi( a_i) b_i \rangle  \ \ \ \ \ , $$
which implies that
$$ \| \sum_{i=1}^n \lafi(y a_i) b_i \| \leq \|y\| \, \| \sum_{i=1}^n  
\lafi(a_i) b_i\| .$$
As usual, this last equality implies the existence of a unique continuous
linear map $T_y$ from $A \ot_\vfi B$ into $A \ot_\vfi B$
such that $T_y (\lafi(a) b) = \lafi(y a) b$ for every $a \in \Nfi$ and
$b \in B$.

We have for every $a_1,a_2 \in \Nfi$ and $b_1,b_2 \in B$ that
\begin{eqnarray*}
\langle T_x( \lafi(a_1)b_1 ) , \lafi(a_2) b_2 \rangle & = &
\langle \lafi(x a_1) b_1 , \lafi(a_2) b_2 \rangle
=  b_2^* \, \vfi(a_2^* (x a_1)) \, b_1
=  b_2^* \, \vfi((x^* a_2)^* a_1) \, b_1 \\
& = & \langle \lafi(a_1) b_1 , \lafi(x^* a_2) b_1 \rangle
= \langle \lafi(a_1) b_1 , T_{x^*}( \lafi(a_2) b_2 ) \rangle
\end{eqnarray*}
This implies that $\langle T_x(v) , w \rangle = \langle v , T_{x^*}(w) 
\rangle$ for every $v,w \in A \ot_\vfi B$. Therefore, we get that $T_x$ 
belongs to $\cL(A \ot_\vfi B)$ and $T_x^* = T_{x^*}$.
\end{demo}

This lemma justifies the following definition.

\begin{definition}
We define the mapping $\pifi$ from $A$ into $\cL(A \ot_\vfi B)$
such that $\pifi(x) \lafi(a) = \lafi(x a) $ for every
$x \in A$, $a \in \Nfi$. Then $\pifi$ is a $^*$-homomorphism.
\end{definition}

This discussion allows us to formulate the following proposition.

\begin{proposition}
We have that $(A \ot_\vfi B, \lafi, \pifi)$ is a KSGNS-construction for
$\vfi$.
This triplet is called the canonical KSGNS-construction for $\vfi$.
\end{proposition}

\bigskip

In a last part, we look at the case where $\vfi$ takes values in $B$ and 
investigate the connection with the mapping $\lafi$. First we prove the 
following lemma.

\begin{lemma} \label{lem1.1}
Consider a \cst-algebra $C$ and a Hilbert \cst-module $F$ over $C$. Let $t$ 
be an element in $\cL(C,F)$. Then $t^* t$ belongs to $C$
$\Leftrightarrow$
There exists an element $x \in F$ such that $t c = x c $ for every $c \in C$ 
$\Leftrightarrow$ $t$ belongs to ${\cal K}(C,F)$ .
\end{lemma}
\begin{demo}
\begin{itemize}
\item Suppose that $t^* t$ belongs to $C$. Define $S \in \cL(C \oplus F , C 
\oplus F)$  with $S^* = S$ such that
$$ S = \left( \begin{array}{cc} 0 & t^* \\ t  & 0 \end{array} \right) .$$
We have that
$$ S^2 = \left( \begin{array}{cc} t^* t & 0 \\ 0 & t t^* \end{array} \right) 
.$$
By proposition 1.4.5 of \cite{Ped}, there exists an element $U \in \cL(C 
\oplus F, C \oplus F)$ such that $S = U (S^2)^\frac{1}{4}$.
This implies that
$$  \left( \begin{array}{cc} 0 & t^* \\ t  & 0 \end{array} \right)
 = U \, \left( \begin{array}{cc} (t^* t)^\frac{1}{4} & 0 \\ 0 &
(t t^*)^\frac{1}{4} \end{array} \right).$$
From this, we get the existence of $u \in \cL(C,F)$ such that $t = u (t^* 
t)^\frac{1}{4}$.
Because $t^* t$ belongs to $C$, we have that $(t^* t)^\frac{1}{4}$ belongs to 
$C$. Put $x = u( (t^* t)^\frac{1}{4} ) \in F$. Then $t(c) = x c$ for every $c 
\in C$.
\item Suppose there exists an element $x \in F$ such that $t(c) = x c$ for 
every $c \in C$. In that case, we have for every $c \in C$ that
$$ c^* (t^* t) c  = \langle c , (t^* t) c \rangle = \langle t c , t c \rangle 
= \langle x c , x c \rangle = c^* \langle x , x \rangle c $$
which implies that $t^* t = \langle x , x \rangle \in C$.
\end{itemize}
\end{demo}

This implies immediately the following proposition.

\begin{proposition}
Consider two \cst-algebras $A,B$ and a \cst-valued weight $\vfi$ from $A$ 
into $M(B)$ with KSGNS-construction $(E,\la,\pi)$. Let $a$ be an element of 
$\Nfi$. \newline
Then $\vfi(a^* a)$ belongs to $B$  $\Leftrightarrow$ There exists an element 
$x \in E$ such that $\la(a) b = x b$ for every $b \in B$ $\Leftrightarrow$
$\la(a)$ belongs to ${\cal K}(B,E)$.
\end{proposition}

Let $a$ be an element in $\Nfi$ such that $\vfi(a^* a) \in B$. By the 
previous proposition, we get the existence of a unique element $x \in E$ such 
that $\la(a) b = x b$ for every $b \in B$. As we have seen in the previous 
lemma, we have in this case that $\vfi(a^* a) = \langle x , x \rangle$ which 
looks like the GNS-construction for weights.

\section{A special family of completely positive mappings.}
\label{art2}

In this section, we will consider two \cst-algebras $A$ and $B$,
a Hilbert \cst-module $E$ over $B$ and a dense left ideal $N$ in $A$.
Furthermore, let $\la$ be a linear mapping from $N$ into $\cL(B,E)$
and $\pi$ a $^*$-homomorpism from $A$ into $\cL(E)$ such that
\begin{itemize}
\item The set $\langle \, \la(a) b \mid a \in N, b \in B \rangle$
is dense in $E$.
\item We have that $\pi(x) \la(a) = \la(x a)$ for every $x \in A$ and $a \in 
N$.
\end{itemize}

We introduce a special family of completely positive mappings relative to 
this objects and investigate some properties of them.

\begin{definition}
We define ${\cal H}$ to be the set
\begin{eqnarray*}
& & \{ \, \,  \rho \text{ a strict completely positive linear mapping from } 
A \text{ into } M(B) \ \mid
\text{ There exists an operator } \\
& & \, \, \, \, T \in {\cal L}(E)^+    \text{ such that } b_2^* \,  
\rho(a_2^* \, a_1) \, b_1 = \langle T \la(a_1) b_1 , \la(a_2) b_2 \rangle  
\text{ for every } a_1,a_2 \in N , b_1,b_2 \in B \,\, \}
\end{eqnarray*}
\end{definition}

\begin{notation}
Consider $\rho \in {\cal H}$. Then there exist a unique element
$T \in \cL(E)^+$ such that \newline $\langle T \la(a_1) b_1 , \la(a_2) b_2 
\rangle
= b_2^* \, \rho(a_2^* \, a_1) \, b_1$ for every $a_1,a_2 \in N$ and $b_1,b_2 
\in B$. We define $T_\rho = T$.

This implies that $\rho(a_2^* \, a_1) = \la(a_2)^* T_\rho \, \la(a_1)$
for every $a_1,a_2 \in N$.
\end{notation}

\begin{result}
Consider $\rho \in {\cal H}$. Then $T_\rho$ belongs to $\pi(A)'$.
\end{result}
\begin{demo}
Choose $x \in A$. Using the definition of $T_\rho$, we have for every 
$a_1,a_2 \in N$ and $b_1,b_2 \in B$ that
\begin{eqnarray*}
\langle T_\rho \pi(x) \la(a_1) b_1 , \la(a_2) b_2 \rangle
& = & \langle T_\rho \la(x a_1) b_1 , \la(a_2) b_2 \rangle
= b_2^* \, \rho( a_2^* (x a_1) ) \, b_1  \\
& = & b_2^* \, \rho((x^* a_2)^* a_1) \, b_1
= \langle T_\rho \la(a_1) b_1 , \la(x^* a_2) b_2 \rangle \\
& = & \langle T_\rho \la(a_1) b_1 , \pi(x^*) \la(a_2) b_2 \rangle
= \langle \pi(x) T_\rho \la(a_1) b_1 , \la(a_2) b_2 \rangle
\end{eqnarray*}
This implies that $T_\rho \pi(x) = \pi(x) T_\rho$.
\end{demo}

\begin{result}
Consider $\rho \in {\cal H}$ and $a \in N$. Then
$\|T_\rho^\frac{1}{2} \la(a)\|^2 = \|\rho(a^* a)\| \leq \|\rho\| \, \|a\|^2$.
\end{result}

This follows immediately from the fact that
$$(T_\rho^\frac{1}{2} \la(a))^* (T_\rho^\frac{1}{2} \la(a)) =
\la(a)^* T_\rho \la(a) = \rho(a^* a) .$$

The following result has its analogue in the theory of weights (see 
proposition 2.4 of \cite{Comb}).

\begin{proposition} \label{prop2.1}
Assume that $\pi$ is non-degenerate. Consider $\rho \in {\cal H}$ and
$S \in \cL(E) \cap \pi(A)'$ such that $S^* S = T_\rho$.
Then there exists a unique element $v \in \cL(B,E)$ such that
$S \la(a) = \pi(a) v$ for every $a \in N$. \newline
Furthermore, the equality  $\|v\|^2 = \|\rho\|$ holds.
We have also that $\rho(x) = v^* \pi(x) v$ for every $x \in M(A)$.
\end{proposition}
\begin{demo}
The non-degeneracy of $\pi$ clearly implies the unicity of $v$. We will turn 
our attention to the existence.
Take an approximate unit $(e_k)_{k \in K}$
for $A$ in $N$ (which is possible because $N$ is a dense left ideal in $A$).

\begin{enumerate}
\item Choose $b \in B$.

We have for every $k,l \in K$ that
\begin{eqnarray*}
\|S \la(e_k) b - S \la(e_l) b\|^2
& = &  \|\langle S \la(e_k - e_l) b , S \la(e_k - e_l)b \rangle\|
= \|\langle S^* S \la(e_k - e_l) b , \la(e_k - e_l)b \rangle\| \\
& = & \|\langle T_\rho \la(e_k - e_l) b , \la(e_k - e_l)b \rangle\|
= \| b^* \rho((e_k-e_l)^2) b \|
\end{eqnarray*}
Because the net $(\,(e_k - e_l)^2 )_{(k,l) \in K \times K}$ is bounded
and converges strictly to $0$, we have also
that the net $(\,\rho((e_k - e_l)^2)\, )_{(k,l) \in K \times K}$
converges strictly to $0$ (remember that $\rho$ is assumed to be strict).

Therefore, the net $(b^* \rho((e_k-e_l)^2) b )_{(k,l) \in K \times K}$ 
converges to $0$, which implies that the net
$(S\la(e_k)b - S\la(e_l)b)_{(k,l) \in K \times K}$ converges to 0.
Consequently, the net $(S \la(e_k) b)_{k \in K}$ is Cauchy and hence
convergent in $B$.

From this all, we get the existence of a linear mapping $v$ from $B$ into
$E$ such that $(S \la(e_k))_{k \in K}$ converges strongly to $v$.

(We borrowed the idea for the preceding part from theorem 5.6 of \cite{Lan}.)

\item Choose $a \in N$.
It is clear that $(\pi(a) S \la(e_k))_{k \in K}$ converges strongly
to $\pi(a) v$ \ \ \ \ (a).

We have for every $k \in K$ that
\begin{eqnarray*}
\|\pi(a) S \la(e_k) - S \la(a)\|^2 & = & \| S \pi(a) \la(e_k) - S \la(a) \|^2 
= \| S \la(a e_k) -  S \la(a) \|^2 \\
& = & \| \la(a e_k - a)^* S^* S \la(a e_k - a) \|
= \| \la(a e_k - a)^* T_\rho \la(a e_k - a) \| \\
& = & \| \rho( (a e_k - a)^* (a e_k -a))\|  .
\end{eqnarray*}
This implies that the net $(\pi(a) S \la(e_k))_{k \in K}$ converges
to $S \la(a)$ \ \ \ \ (b).

Combining this result with (a), we get the equality $\pi(a) v =  S \la(a)$.

\item We have for every $k \in K$ that
$$\| S \la(e_k) \|^2 = \| \la(e_k)^* S^* S \la(e_k) \| = \| \la(e_k)^* T_\rho 
\la(e_k) \| = \| \rho((e_k)^2) \| .$$

Hence, $\| S \la(e_k) \|^2 \leq \|\rho\|$ for every
$k \in K$. Because $(S \la(e_k))_{k \in K}$ converges strongly to $v$,
this implies that $v$ is bounded and $\|v\|^2 \leq \|\rho\|$.

By 2, we have for every $k \in K$ that
$$\|v\|^2 \geq \| \pi(e_k) v\|^2 = \|S \la(e_k)\|^2 =  \| \rho((e_k)^2) \| \ 
.$$
 Because $\rho$ is positive, we know that
$(\| \rho((e_k)^2) \|)_{k \in K}$ converges to $\|\rho\|$.
This implies that $\|v\|^2 \geq \|\rho\|$.

Combining these two inequalities, we get that $\|v\|^2 = \|\rho\|$.

\item We still have to prove that $v$ is adjointable.

Therefore, choose $a \in N$ and $x \in E$.
We know from (b) that $(\, (S \la(e_k))^* \pi(a^*)\,)_{k \in K}$ converges to 
$(S \la(a))^*$. This implies that the net
$(\,(S \la(e_k))^* \pi(a^*) x\,)_{k \in K}$  is convergent.

Because $\pi$ is assumed to be non-degenerate, the set  $\langle \, \pi(a^*) 
x \mid a \in N , x \in E \rangle$ is dense in $E$. By 3, we also know that 
$(\,(S \la(e_k))^*)_{k \in K}$ is bounded.

These three facts imply easily the existence of a unique linear map
$w$ from $E$ into $B$ such that $(\,(S \la(e_k))^*)_{k \in K}$ converges
strongly to $w$.
By definition, we have that $(S \la(e_k))_{k \in K}$ converges strongly
to $v$.

These two results imply that
$\langle v b , x \rangle = \langle b , w x \rangle$ for $x \in E$, $b \in B$. 
Consequently, $v$ belongs to $\cL(B,E)$ and $v^* = w$.

\item We have for every $a_1,a_2 \in \Nfi$ that
\begin{eqnarray*}
\rho(a_2^* \, a_1) & = & \la(a_2)^* T_\rho \, \la(a_1) = \la(a_2)^* S^* S 
\la(a_1) = (S \la(a_2))^* (S \la(a_1)) \\
& = & (\pi(a_2) v)^* (\pi(a_1) v) = v^* \pi(a_2^* \, a_1) v  \ .
\end{eqnarray*}
This implies that $\rho(x) = v^* \pi(x) v$ for every $x \in \Mfi$. The usual 
continuity and strict continuity arguments imply that $\rho(x) = v^* \pi(x) 
v$ for every $x \in M(A)$.
\end{enumerate}
\end{demo}

\begin{notation} \label{not2.1}
Suppose that $\pi$ is non-degenerate. Consider $\rho \in {\cal H}$, we define 
$v_\rho$ to be the unique element in $\cL(B,E)$ such that $T_\rho^\frac{1}{2} 
\la(a) = \pi(a) v_\rho$ for every $a \in N$.
Furthermore, $\|v_\rho\|^2 = \|\rho\|$. \newline
We have also that $\rho(x) = v_\rho^* \pi(x) v_\rho$ for every $x \in M(A)$.
\end{notation}

\medskip

\begin{remark} \rm
We have the following obvious properties:
\begin{enumerate}
\item 0 belongs to ${\cal H}$ and $T_0 = 0$.
\item For every $\rho_1,\rho_2 \in {\cal H}$, we have that $\rho_1+\rho_2 $ 
belongs to ${\cal H}$ and $T_{\rho_1+\rho_2} = T_{\rho_1}+T_{\rho_2}$.
\item For every $\rho \in {\cal H}$ and every $\lambda \in \R^+$, we have
that $\lambda T$ belongs to ${\cal H}$ and $T_{\lambda \rho} = \lambda 
T_\rho$.
\end{enumerate}
\end{remark}

We will need the following ordering $\leq$ on ${\cal H}$.

\begin{definition}
Consider $\rho_1,\rho_2 \in {\cal H}$. We say that $\rho_1 \leq \rho_2$
\begin{description}
\item[$\Leftrightarrow$] $\rho_2 - \rho_1$ is completely positive.
\item[$\Leftrightarrow$] We have for every $n \in \N$, $a_1,\ldots\!,a_n \in 
A$ and $b_1,\ldots\!,b_n \in B$ that
$$\sum_{i,j=1}^n b_j^* \, \rho_1(a_j^* \, a_i) \, b_i
\leq \sum_{i,j=1}^n b_j^*  \, \rho_2(a_j^* \, a_i) \, b_i .$$
\item[$\Leftrightarrow$] $T_{\rho_1} \leq T_{\rho_2}$.
\end{description}
\end{definition}

It is straigthforward to check the equivalence of the different conditions in 
the previous definition.

\begin{remark} \rm
We also have the following property:

Let $\rho_1, \rho_2$ be elements in ${\cal H}$ with $\rho_1 \leq \rho_2$. 
Then $\rho_2 - \rho_1$ belongs to ${\cal H}$ and $T_{\rho_2-\rho_1}
= T_{\rho_1} - T_{\rho_2}$.
\end{remark}

\medskip

The following sets are also generalizations of objects, known in ordinary
weight theory (see e.g definition 2.1.7 of \cite{Verd}).

\begin{definition}
We define the following sets
$${\cal F} = \{ \rho \in {\cal H} \mid T_\rho \leq 1 \}$$
and
$${\cal G} = \{ \lambda \rho \mid \lambda \in [0,1[ \, , \rho \in {\cal F} \} 
\subseteq {\cal F}.$$
\end{definition}

\begin{remark} \rm Let $\rho$ be an element in ${\cal H}$. It is not 
difficult to check that $\rho$ belongs to ${\cal F}$ if and only if
$$\sum_{i,j=1}^n b_j^* \, \rho(a_j^* \, a_i) \, b_i \leq \langle \sum_{i=1}^n 
\la(a_i) b_i , \sum_{i=1}^n \la(a_i) b_i \rangle $$
for every $n \in \N$, $a_1,\ldots\!,a_n \in N$ and $b_1,\ldots\!,b_n \in B$.
\end{remark}

Like in the case of weights, the reason for introducing the set ${\cal G}$ is 
to obtain a set which is directed for the ordering $\leq$.
This is the content of the following proposition. The basic idea of its proof 
is the same as the proof for ordinary weights (see proposition 2.1.8 of 
\cite{Verd}), but some
extra work has to be done in this case.

\medskip

First, we prove the following lemma (which we took from proposition
2.1.8 of \cite{Verd}).

\begin{lemma}
Consider a unital \cst-algebra $C$. Let $T_1,T_2$ be elements in $C$ with
$0 \leq T_1,T_2 \leq 1$, let $\gamma$ be a number in $[0,1[$.
Then there exist an element $T \in C$ with $0 \leq T \leq 1$ and such
that $\gamma \, T_1 \leq T$,  $\gamma \, T_2 \leq T$
and $T \leq \frac{\gamma}{1-\gamma} (T_1+T_2)$.
\end{lemma}
\begin{demo}
For the moment, fix $i \in \{1,2\}$. We  define
$$S_i = \frac{\gamma \, T_i}{1-\gamma \, T_i} \in C.$$

It is then easy to check that
$$ \gamma \, T_i = \frac{S_i}{1+S_i}   \text{\ \ \ \ \ \ \ (a)} 
\hspace{1.5cm}  \text{ and } \hspace{1.5cm}
S_i \leq \frac{\gamma}{1-\gamma}\, T_i \text{\ \ \ \ \ \ \ (b)} \ .$$

Next, we define
$$T = \frac{S_1 + S_2}{1+S_1+S_2} \in C .$$
We get immediately that $0 \leq T \leq 1$.
By (b), we have that $T \leq S_1+S_2 \leq \frac{\gamma}{1-\gamma} (T_1 + 
T_2)$.

We know that the function $\R^+ \rightarrow \R^+ : t \mapsto \frac{t}{1+t}$ 
is operator monotone (see \cite{Ped}). This implies
for every $i=1,2$ that
$$ T \geq \frac{S_i}{S_i+1} = \gamma \, T_i \ ,$$
where we used (a) in the last equality.
\end{demo}

\begin{proposition}
The set ${\cal G}$ is directed for $\leq$ \ .
\end{proposition}
\begin{demo}
Choose $\rho_1, \rho_2 \in {\cal F}$ and $\lambda_1,\lambda_2 \in [0,1[$.
Then there exist a number $\gamma \in [0,1[$ such that $\lambda_1,\lambda_2 < 
\gamma$.

We also know that $T_{\rho_1},T_{\rho_2}$ belong to $\pi(A)' \cap \cL(E)$
and $0 \leq T_{\rho_1},T_{\rho_2} \leq 1$. The previous lemma implies
the existence of $T \in \pi(A)' \cap \cL(E)$ with $0 \leq T \leq 1$ and
such that $\gamma \, T_{\rho_1} \leq T$,  $\gamma \, T_{\rho_2} \leq T$
and $T \leq \frac{\gamma}{1-\gamma} (T_{\rho_1}+T_{\rho_2})$.

\medskip

Put $\lambda = \max(\frac{\lambda_1}{\gamma},\frac{\lambda_2}{\gamma}) \in 
[0,1[$. \ Then $\lambda T \geq \lambda \gamma \, T_{\rho_1} \geq \lambda_1 
T_{\rho_1}
= T_{\lambda_1 \rho_1}$ and analogously, $\lambda T \geq T_{\lambda_2 
\rho_2}$. \ \ \ \ \ (a)

\medskip

In the rest of the proof, we want to construct an element $\rho$ in ${\cal 
F}$ such that $T_{\rho} = T$. \smallskip

In two steps, we will prove the following inequality:

\medskip

Consider $a_1,\ldots\!,a_n \in N$ and $b_1,\ldots\!,b_n \in B$. Then
$$\|\sum_{i=1}^n \la(b_i)^* T \la(a_i)\| \leq 4 \frac{\gamma}{1-\gamma}
(\|\rho_1\| + \|\rho_2\|) \, \| \sum_{i=1}^n b_i^* a_i \| . $$

\begin{enumerate}
\item Choose $\om \in B_+^*$.

Define the function $\th$ from $A$ into $\C$ such that
$\th(x) = \frac{\gamma}{1-\gamma}(\om(\rho_1(x)) + \om(\rho_2(x))\,)$ for
every $x \in A$.

Then $\th$ belongs to $A_+^*$ and $\|\th\| \leq \frac{\gamma}{1 - \gamma} 
(\|\rho_1\| + \|\rho_2\|) \|\om\|$ \ \ \ \ \ (b).

Next, define the mapping $S$ from $N \times N$ into $\C$ such that
$S(x,y) = \om(\la(y)^* T \la(x))$ for every $x,y \in N$. Then $S$ is 
sesquilinear. Furthermore:
\begin{itemize}
\item We have for $x \in N$ that
$$ 0 \leq \la(x)^* T \la(x) \leq \frac{\gamma}{1-\gamma} \, 
\la(x)^*(T_{\rho_1}+T_{\rho_2})\la(x)
\leq \frac{\gamma}{1-\gamma} (\rho_1(x^* x) + \rho_2(x^* x)) , $$
which implies that $0 \leq S(x,x) \leq \th(x^* x)$.
\item We have for every $x,y \in N$ and $a \in A$ that
\begin{eqnarray*}
\la(y)^* T \la(a x) & = & \la(y)^* T \pi(a) \la(x)
= \la(y)^* \pi(a) T \la(x) \\
& = & (\pi(a^*) \la(y))^* T \la(x) = \la(a^* y)^* T \la(x) ,
\end{eqnarray*}
which implies that $S(a x,y) = S(x, a^* y)$.
\end{itemize}

This allows us to apply lemma \ref{lemA6}. Therefore, we get the existence of 
$\psi \in B_+^*$ with $\psi \leq \th$ and such that
$S(x,y) = \psi(y^* x)$ for every $x,y \in N$.

This implies that
\begin{eqnarray*}
|\om(\,\sum_{i=1}^n \la(b_i)^* T \la(a_i)\,)| & = &
|\sum_{i=1}^n S(a_i,b_i)| = |\psi(\sum_{i=1}^n b_i^* a_i)| \leq \|\psi\| \, 
\|\sum_{i=1}^n b_i^* a_i\| \\
& \leq & \|\th\| \, \|\sum_{i=1}^n b_i^* a_i\|
\leq \frac{\gamma}{1-\gamma}(\|\rho_1\|+\|\rho_2\|) \|\om\| \, \|\sum_{i=1}^n 
b_i^* a_i \|
\end{eqnarray*}
where we used inequality (b) in the last inequality.

\item Choose $\om \in B^*$ with $\|\om\| \leq 1$. Then there exist 
$\om_1,\ldots\!,\om_4 \in B_+^*$ such that $\|\om_1\|,\ldots\!,\|\om_4\| \leq 
1$ and $\om = \sum_{k=1}^4 i^k \, \om_k$.
Referring to 1), it is not difficult to see that
$$|\om(\sum_{i=1}^n \la(b_i)^* T \la(a_i))| \leq 4 \frac{\gamma}{1-\gamma} 
(\|\rho_1\|+\|\rho_2\|) \, \|\sum_{i=1}^n b_i^* a_i\| .$$

This implies that
$$\|\sum_{i=1}^n \la(b_i)^* T \la(a_i)\| \leq 4 \frac{\gamma}{1-\gamma}
(\|\rho_1\|+\|\rho_2\|) \, \|\sum_{i=1}^n b_i^* a_i\| .$$

\end{enumerate}

As usual, this inequality guarantees the existence of a unique continuous 
linear map $\rho$ from $A$ into $M(B)$ such that $\rho(a_2^* \, a_1) = 
\la(a_2)^* T \la(a_1)$ for $a_1,a_2 \in N$. It is clear that $\rho$ is 
completely positive.

Choose $x \in N$. We have for every $b \in B$ that
\begin{eqnarray*}
b^* \rho(x^* x) b & = & \langle T \la(x) b , \la(x) b \rangle
\leq \frac{\gamma}{1-\gamma} (\langle T_{\rho_1} \la(x) b , \la(x) b \rangle 
+ \langle T_{\rho_2} \la(x) b , \la(x) b \rangle) \\
& = & \frac{\gamma}{1-\gamma} (b^* \rho_1(x^* x) b + b^* \rho_2(x^* x) b) , 
\end{eqnarray*}
which implies that $\rho(x^* x) \leq \frac{\gamma}{1-\gamma}
(\rho_1(x^* x) + \rho_2(x^* x))$.
Hence, $\rho \leq \frac{\gamma}{1-\gamma} (\rho_1 + \rho_2)$.
Therefore, the strictness of $\rho_1,\rho_2$ implies that $\rho$ is strict.

\smallskip

It is now easy to see that $\rho$ belongs to ${\cal F}$ and
$T_\rho = T$.
Moreover, inequality (a) implies that $T_{\lambda \rho} = \lambda T_\rho = 
\lambda T \geq
T_{\rho_1},T_{\rho_2}$, which implies that $\lambda \rho \geq \lambda_1 
\rho_1, \lambda_2 \rho_2$.
\end{demo}

\bigskip

The following easy lemma will be used several times.

\begin{lemma} \label{lem2.1}
Consider $a \in M(A)^+$ and $b \in B$. Let $x$ be an element in $B$
such that $b^* \rho(a) b \leq x$ for every $\rho \in \cG$.
Then the net $(b^* \rho(a) b)_{\rho \in \cG}$ converges to $x$ 
$\Leftrightarrow$
For every $\vep > 0$ there exist an element $\eta \in \cF$ such that
$\| x - b^* \eta(a) b \| \leq \vep$.
\end{lemma}
\begin{demo}
One implication is trivial. We will prove the one from the right to the left.

We have immediately that $b^* \rho(a) b \leq x$ for every $\rho \in \cF$.

\medskip

Choose $\vep > 0$. By assumption, there exist an element $\eta \in \cF$
such that $\| x - b^* \eta(a) b \| \leq \frac{\vep}{2}$.

Furtermore, there exist a number $\lambda \in [0,1[$ such that
$|\lambda - 1| \leq \frac{\vep}{2(\|x\|+1)}$.

Put $\rho_0 = \lambda \, \eta$. So $\rho_0$ belongs to $\cG$.
We have also that
\begin{eqnarray*}
\| b^* \rho_0(a) b - x \| & \leq & \| b^* \rho_0(a) b -  b^* \eta(a) b \|
+ \| b^* \eta(a) b - x \| \\
& \leq &  |\lambda-1| \, \|b^* \eta(a) b\| + \frac{\vep}{2}
\leq  \frac{\vep}{2(\|x\|+1)} \, \|x\| + \frac{\vep}{2} \leq \vep .
\end{eqnarray*}

We have for every $\rho \in \cG$ with $\rho_0 \leq \rho$ that
$0 \leq b^* \rho_0(a) b \leq b^* \rho(a) b \leq x$, which implies
that $$\| b^* \rho(a) b - x \| \leq \|b^* \rho_0(a) b - x \| \leq \vep .$$

Consequently, $(b^* \rho(a) b)_{\rho \in \cG}$ converges to $x$.
\end{demo}

\section{Lower semi-continuity of \cst-valued weights.}
\label{art3}

In this section, we will introduce a possible definition for lower semi-
continuity of a \cst-valued weight. Also, some easy consequences will be 
derived for such lower semi-continuous \cst-valued weights. At the end of 
this section, we will introduce an ever stronger condition than lower semi-
continuity, the so-called regularity.

\begin{notation} \label{not3.1}
Consider two \cst-algebras $A$ and $B$ and a densely defined \cst-valued 
weight $\vfi$ from $A$ into $M(B)$ with KSGNS-construction $(E,\la,\pi)$.  
\newline The ingredients $A,B,E,\Nfi,\la,\pi$ satisfy the conditions of the 
beginning of  section 2.
We define ${\cal H}_\vfi$, ${\cal F}_\vfi$, ${\cal G}_\vfi$ to be
the objects ${\cal H}$, ${\cal F}$, ${\cal G}$ constructed in section 2, 
using this specific ingredients.
\end{notation}

It is not difficult to see that the definition of ${\cal H}_\vfi$ , ${\cal 
F}_\vfi$ and ${\cal G}_\vfi$ is independent of the choice of the KSGNS-
construction. The notations $T_\rho$ and $v_\rho$ however, will always be 
relative to a specific  KSGNS-construction.

\medskip

We will use the following definition of lower semi-continuity.

\begin{definition} \label{def3.1}
Consider two \cst-algebras $A$ and $B$ and a densely defined \cst-valued 
weight $\vfi$ from $A$ into $M(B)$. We call $\vfi$ lower semi-continuous if 
and only if
\begin{enumerate}
\item We have that $ \Mfi^+ = \{ x \in A \mid (\rho(x))_{\rho \in {\cal 
G}_\vfi} \text{ is strictly convergent in } M(B) \} $.
\item For every $x \in \Mfi^+$, the net $(\rho(x))_{\rho \in {\cal G}_\vfi}$ 
converges strictly to $\vfi(x)$.
\end{enumerate}
\end{definition}

\begin{remark} \rm
Consider two \cst-algebras $A$ and $B$ and a densely defined lower
semi-continuous \cst-valued weight $\vfi$ from $A$ into $M(B)$.
Because every element of $\Mfi$ is a linear combination of elements from 
$\Mfi^+$, the net $(\rho(x))_{\rho \in {\cal G}_\vfi}$ converges strictly to 
$\vfi(x)$ for every $x \in \Mfi$.
\end{remark}

\medskip

This definition is rather heavy. In fact, in ordinary weight theory this 
definition is a major result proved by Combes (see \cite{Comb}). Therefore, 
it might be an interesting question whether it is possible to give some kind 
of lower semi-continuity definition in the classical sense which is 
equivalent to this one.

On the other hand, this definition will probably be workable because most of 
the C$^*$-valued weights will be defined starting from some sort of family of 
completely positive mappings.

\bigskip

For the most part of this section, we consider two \cst-algebras $A$ and $B$ 
and a densely defined lower semi-continuous \cst-valued weight from $A$ into 
$M(B)$. We also fix a KSGNS-construction $(E,\la,pi)$ for $\vfi$.

\medskip

Referring to lemma \ref{lemA4}, the following proposition follows immediately.

\begin{proposition}
Let $x$ be an element in $A^+$. Then $x$ belongs to $\Mfi^+$ 
$\Leftrightarrow$ We have for every $b \in B$ that the net $(b^* \rho(x) 
b)_{\rho \in {\cal G}_\vfi}$ is convergent in $B$.
\end{proposition}

As usual, we have some kind of monotone convergence properties:

\begin{result} \label{res3.1}
Consider a net $(x_i)_{i \in I}$ in $\Mfi^+$ and an element $x$ in $\Mfi^+$ 
such that $(x_i)_{i \in I}$ converges strictly to $x$ and
$x_i \leq x$ for every $i \in I$.
Then $(\vfi(x_i))_{i \in I}$ converges strictly to $\vfi(x)$.
\end{result}
\begin{demo}
It is clear that $0 \leq \vfi(x_i) \leq \vfi(x)$ for every $i \in I$.

Choose $b \in B$.

Take $\vep > 0$. Then there exist $\rho \in {\cal G}_\vfi$ such that
$\| b^* \vfi(x) b - b^* \rho(x) b \| \leq \frac{\vep}{2}$.

Because $(\rho(x_i))_{i \in I}$ converges strictly to $\rho(x)$, there
exists an element $i_0 \in I$ such that
\newline $\|b^* \rho(x_i) b - b^* \rho(x) b\| \leq \frac{\vep}{2}$ for every 
$i \in I$ with $i \geq i_0$.

Choose $j \in J$ with $j \geq i_0$. Then
$$\| b^* \rho(x_j) b - b^* \vfi(x) b \|
\leq \| b^* \rho(x_j) b - b^* \rho(x) b \| +
\| b^* \rho(x) b - b^* \vfi(x) b \|
\leq \frac{\vep}{2} + \frac{\vep}{2} = \vep .$$
Because $0 \leq b^* \rho(x_j) b \leq b^* \vfi(x_j) b \leq b^* \vfi(x) b $, we 
get that
$$\| b^* \vfi(x_j) b - b^* \vfi(x) b \|
\leq \| b^* \rho(x_j) b - b^* \vfi(x) b \| \leq \vep.$$

Consequently, we see that $(b^* \vfi(x_i) b)_{i \in I}$ converges to $b^* 
\vfi(x) b$. Lemma \ref{lemA3} implies that $(\vfi(x_i))_{i \in I}$ converges 
strictly to $\vfi(x)$.
\end{demo}

\begin{result} \label{res3.2}
Consider a net $(x_i)_{i \in I}$ in $\Mfi^+$ and an element $x$ in $A^+$ such 
that $(x_i)_{i \in I}$ converges strictly to $x$ and
$x_i \leq x$ for every $i \in I$.
Then $x$ belongs to $\Mfi^+$ $\Leftrightarrow$ The net $(b^* \vfi(x_i) b)_{i 
\in I}$ is convergent for every $b \in B$.
\end{result}
\begin{demo}
One implication follows from the previous result, we will turn to the other 
one. Therefore, assume that the net $(b^* \vfi(x_i) b)_{i \in I}$ is 
convergent for every $b \in B$.

\medskip

Choose $c \in B$.

By assumption there exists an element $d \in B^+$ such that
$(c^* \vfi(x_i) c)_{i \in I}$ converges to $d$  \ \ \ \ \ (a).

First, we will prove that $c^* \rho(x) c \leq d $ for every $\rho \in
\cG_\vfi$ \ \ \ \ \ (b).

\begin{list}{}{\setlength{\leftmargin}{.4 cm}}

\item Take $\rho \in \cG_\vfi$. Choose $n \in \N$.
From (a), we get the existence of an element $i_0 \in I$
such that $\| d - c^* \vfi(x_i) c \| \leq \frac{1}{n}$ for every $i \in I$ 
with $i \geq i_0$.
Therefore, we have for every $i \in I$ with $i \geq i_0$
that $c^* \rho(x_i) c \leq c^* \vfi(x_i) c \leq d + \frac{1}{n} \, 1$

Because $(c^* \rho(x_i) c)_{\i \in I}$ converges to $c^* \rho(x) c$,
this implies that $c^* \rho(x) c \leq d + \frac{1}{n} \, 1$.

If we let $n$ tend to $\infty$, we get that $c^* \rho(x) c \leq d$.

\end{list}

Choose $\vep > 0$.
By (a), there exists an element $j \in I$ such that
$\| c^* \vfi(x_j) c - d \| \leq \frac{\vep}{2}$.

The lower semi-continuity of $\vfi$ implies the existence of $\rho_0 \in 
\cG_\vfi$ such that $\| c^* \vfi(x_j) c - c^* \rho_0(x_j) c \| \leq 
\frac{\vep}{2}$ for every $\rho \in \cG_\vfi$ with $\rho \geq \rho_0$.
So we get that $\|c^* \rho(x_j) c - d \| \leq \vep$ for every
$\rho \in \cG_\vfi$ with $\rho \geq \rho_0$ \ \ \ \ \ (c).

\smallskip

Take $\eta \in \cG_\vfi$ with $\eta \geq \rho_0$.
Using (b), we have that $0 \leq c^* \eta(x_j) c \leq c^* \eta(x) c \leq d$. 
Consequently, (c) implies that
$$\| c^* \eta(x) c - d \| \leq \| c^* \eta(x_j) c - d \| \leq \vep \ . $$

\medskip

Hence, we see that $(c^* \rho(x) c)_{\rho \in \cG_\vfi}$ converges to
$d$. The lower semi-continuity of $\vfi$ guarantees that $x$ belongs
to $\Mfi^+$.
\end{demo}

\medskip

Like for weights, these kind of convergence properties imply that the 
representation of a KSGNS-construction is non-degenerate (lemma 2.1 od 
\cite{Comb}).

\begin{result}
The mapping  $\pi$ is non-degenerate and $\pi(x) \la(a) = \la(x a)$ for every 
$x \in M(A)$ and $a \in \Nfi$.
\end{result}
\begin{demo}
Take an approximate unit $(e_i)_{i \in I}$ for $A$.

Choose $a \in \Nfi$ and $b \in B$. Then $(a^* e_i a)_{i \in I}$ is a net in 
$\Mfi^+$ such that $a^* e_i a \leq a^* a$ for every $i \in I$. Therefore, 
result \ref{res3.1} implies
that $(b^* \, \vfi(a^* e_i a) \, b)_{i \in I}$ converges to $\vfi(a^* a)$. \ 
\ \ (a)

Similarly, we get that $(b^* \, \vfi(a^* e_i^2 \, a) \, b)_{i \in I}$ 
converges to $b^* \vfi(a^* a) b$. \ \ \ (b)

We have for every $i \in I$ that
\begin{eqnarray*}
\| \pi(e_i) \la(a) b - \la(a) b \|^2 & = &
\|\langle \la(e_i a) b - \la(a) b , \la(e_i a) b - \la(a) b \rangle \| \\
& = & \| b^* \, \vfi(a^* e_i^2 \, a) \, b - b^* \, \vfi(a^* e_i a) \, b - b^* 
\, \vfi(a^* e_i a) \, b + b^* \, \vfi(a^* a) \, b \|
\end{eqnarray*}
Therefore, (a) and (b) imply that
$(\pi(e_i) \la(a) b)_{i \in I}$ converges to $\la(a) b$.

Because $(\pi(e_i))_{i \in I}$ is bounded, we can conclude from this
that $(\pi(e_i))_{i \in I}$ converges strongly to 1.
\end{demo}

\begin{result} \label{res3.3}
The net $(T_\rho)_{\rho \in {\cal G}_\vfi}$ converges strongly to 1.
\end{result}
\begin{demo}
Choose $a_1,a_2 \in \Nfi$ and $b_1,b_2 \in B$.

We have for every $\rho \in {\cal G}_\vfi$ that $\langle T_\rho \la(a_1) b_1 
, \la(a_2) b_2 \rangle = b_2^* \, \rho(a_2^* \, a_1) \, b_1$.

Therefore, the lower semi-continuity of $\vfi$ implies that $(\langle T_\rho 
\la(a_1) b_1 , \la(a_2) b_2 \rangle)_{\rho \in {\cal H}_\vfi}$ converges to 
\newline $b_2^* \, \vfi(a_2^* \, a_1) \, b_1$, which is equal to
$\langle \la(a_1) b_1 , \la(a_2) b_2 \rangle$.

Because $(T_\rho)_{\rho \in {\cal G}_\vfi}$ is bounded, this implies that
$(\langle T_\rho v , w \rangle)_{\rho \in {\cal G}_\vfi}$ converges to
$\langle v , w \rangle$ for all $v,w \in E$.

Using lemma \ref{lemA1}, we conclude that $(T_\rho)_{\rho \in {\cal G}_\vfi}$ 
converges strongly to 1.
\end{demo}

\begin{proposition} \label{prop3.1}
The mapping $\la$ is closed for the strict topology on $A$ and the strong 
topology on $L(B,E)$.
\end{proposition}
\begin{demo}
Choose a net $(x_i)_{i \in I}$ in $\Nfi$, $x \in A$ and $t \in L(B,E)$
such that $(x_i)_{i \in I}$ converges strictly to $x$ and
$(\la(x_i))_{i \in I}$ converges strongly to $t$.

\medskip

Take $\rho \in {\cal G}_\vfi$, $c \in \Nfi$ and $b,d \in B$.

We have for every $i \in I$ that
$\langle T_\rho (\la(x_i)b) , \la(c) d \rangle
= d^* \rho(c^* x_i) b$, which implies that
$(\langle T_\rho (\la(x_i)b) , \la(c) d \rangle)_{i \in I}$ converges
to $d^* \rho(c^* x) b$.
It is also clear that
$(\langle T_\rho (\la(x_i)b) , \la(c) d \rangle)_{i \in I}$ converges
to $\langle T_\rho \, t(b) , \la(c) d \rangle$.

Combining these two facts, we conclude that
$\langle T_\rho \, t(b) , \la(c) d \rangle = d^* \rho(c^* x) b$ \ \ \ \ (a).

\smallskip

Now we will use this last equality to prove that $x \in \Nfi$ and $\la(x) = 
t$.
\begin{itemize}

\item Take $b \in B$. For the moment, fix $\rho \in {\cal G}_\vfi$.

We have immediately that the net $(b^* \,\rho(x_i^* x) \, b)_{i \in I}$ 
converges to $b \, \rho(x^* x) \, b$.
Using (a), we have that  \newline $b \, \rho(x_i^* x) \, b$ $= \langle T_\rho 
\, t(b) , \la(x_i) b \rangle$ for every $i \in I$. This implies that the net 
$(b^* \rho(x_i^* x) b)_{i \in I}$ converges
to $\langle T_\rho t(b) , t(b) \rangle$.
Combining these two results, we conclude that
$\langle T_\rho t(b) , t(b) \rangle = b^* \rho(x^* x) b $.

This last equality implies that $(b^* \, \rho(x^* x) \, b)_{\rho \in {\cal 
G}_\vfi}$  converges to $\langle t(b), t(b) \rangle$. Because $\vfi$ is 
assumed to be lower semi-continuous, this implies that $x^* x$ belongs to 
$\Mfi^+$. So $x$ belongs to $\Nfi$.

\item Choose $c \in \Nfi$ and $b,d \in B$.

By (a), we have for every $\rho \in {\cal G}_\vfi$ that
$$ \langle T_\rho \, t(b) , \la(c) d \rangle =
d^* \, \rho(c^* x) \, b = \langle T_\rho \la(x) b , \la(c) d \rangle .$$

Because $(T_\rho)_{\rho \in {\cal G}_\vfi}$ converges strongly to 1, we can 
conclude that $\langle  t(b) , \la(c) d \rangle
= \langle \la(x) b , \la(c) d \rangle$.

Consequently, $t = \la(x)$.
\end{itemize}
\end{demo}

\bigskip \bigskip

At the end of this section, we introduce a stronger condition than lower semi-
continuity, the so-called regularity condition. A first version of this 
condition was introduced for weights by J. Verding (see definition 2.1.13 of 
\cite{Verd}).

\begin{definition}
Consider two \cst-algebras $A$ and $B$ and a \cst-valued weight $\vfi$ from 
$A$ into $M(B)$ with KSGNS-construction $(E,\la,\pi)$.
We say that $\vfi$ is regular if and only if
\begin{enumerate}
\item $\vfi$ is densely defined and lower semi-continuous.
\item There exist a net $(u_i)_{i \in I}$ in $M(A)$ satisfying the following 
properties:
\begin{description}
\item[{\rm a)}] We have for every $i \in I$ that:
\begin{itemize}
\item $\Nfi u_i \subseteq \Nfi$
\item There exist a unique operator $S_i \in \cL(E)$ such that $S_i \la(a) = 
\la(a u_i)$ for every $a \in \Nfi$.
\item $\|u_i\| \leq 1$ and $\|S_i\| \leq 1$
\item There exists a uniqe strict completely positive linear mapping
$\rho_i$ from $A$ into $M(B)$ such that $b_2^* \, \rho_i(a_2^* \, a_1) \, b_1
= \langle S_i \la(a_1) b_1 , S_i \la(a_2) b_2 \rangle$ for every $a_1,a_2 \in 
\Nfi$ and $b_1,b_2 \in B$.
\end{itemize}
\item[{\rm b)}] Moreover,
\begin{itemize}
\item $(u_i)_{i \in I}$ converges strictly to 1
\item $(S_i)_{i \in I}$ converges strongly to 1.
\end{itemize}
\end{description}
\end{enumerate}
Such a net $(u_i)_{i \in I}$ is called a truncating net for $\vfi$.
If the truncating net can be chosen in such a way that every element belongs 
to $A$, then we call $\vfi$ strongly regular.
\end{definition}

Again, it is not too difficult to check that this definition of regularity 
(and strong regularity) is independent of the choice of the KSGNS-
construction.

\medskip

This regularity condition is a very strong condition. Nevertheless, there
are some interesting \cst-valued weights, which satisfy this condition.

\begin{itemize}
\item Consider two \cst-algebras $A$ and $B$ and a strict completely positive 
mapping $\rho$ from $A$ into $M(B)$. Then $\rho$ is a regular
\cst-valued weight from $A$ into $B$ with a truncating net consisting of one 
element, the unit in $M(A)$.

In particular, the non-degenerate $^*$-homomorphisms from
$A$ into $M(B)$ are regular.
\item Also, the previous example implies that all positive linear functionals 
on a \cst-algebra are regular.
\item Consider a \cst-algebra $A$ and a KMS-weight $\vfi$ on $A$.
This means that $\vfi$ is a densely defined lower semi-continuous weight on 
$A$ such that there exist a norm-continuous one-paramater group $\si$ on $A$ 
such that
\begin{enumerate}
\item $\vfi$ is invariant under $\si$.
\item We have that $\vfi(a^* a) = \vfi(\si_\frac{i}{2}(a)       
\si_\frac{i}{2}(a)^*)$ for every $a \in D(\si_\frac{i}{2})$.
\end{enumerate}
Then $\vfi$ is strongly regular and has a truncating net consisting of 
elements which are analytic with respect to $\si$.

This fact is proven by Jan Verding. The proof of proposition 2.1.18 in 
\cite{Verd} contains a mistake but has been corrected by him.
There will be a modified proof of this fact in \cite{JK1}.

\medskip

It is not yet clear whether every densely defined lower semi-continuous 
weight is regular, but intuition tells us that this will probably not be 
true. It will be easily true if $A$ is commutative.
\end{itemize}

\section{Extension of a lower semi-continuous \cst-valued weight to the 
Multiplier algebra.}
\label{art5}

It is natural to look for extensions of \cst-valued weights to the Multiplier 
algebra. In this section, we will do this for lower semi-continuous \cst-
valued weights. The KSGNS-construction of this extension will also be 
investigated.

\medskip

We start with a lemma which will also be used in
a later section.

\begin{lemma} \label{lem5.1}
Consider two \cst-algebras $A$ and $B$ and an increasing net $(\rho_i)_{i \in 
I}$ of strict completely mappings from $A$ into $M(B)$.
Let $b$ be an element in $B$ and define
$$P = \{ x \in M(A)^+ \mid \text{The net  } (b^* \rho_i(x) b)_{i \in I} 
\text{ is convergent }  \} . $$
Then $P$ is a hereditary cone in $M(A)^+$.
\end{lemma}
\begin{demo}
It is easy to see that $P$ is a cone in $M(A)^+$.
We turn now to the hereditarity of $P$.

Choose $x \in P$ and $y \in M(A)^+$ such that $y \leq x$.

Take $\vep > 0$. Then there exists an element $i_0 \in I$
such that we have for every $i \in I$ with $i \geq i_0$ \newline that
$\| b^* \rho_i(x) b - b^* \rho_{i_0}(x) b \| \leq \vep$.

We have for every $i \in I$ with $i \geq i_0$ that
$$ 0 \leq b^* (\rho_i(y) - \rho_{i_0}(y)) b
\leq b^* (\rho_i(x) - \rho_{i_0}(x)) b  , $$
which implies that
$$ \| b^* \rho_i(y) b - b^* \rho_{i_0}(y) b \| \leq
\| b^* \rho_i(x) b - b^* \rho_{i_0}(x) b \| \leq \vep . $$
Therefore, we see that the net $(b^* \rho_i(y) b)_{i \in I}$
is Cauchy and hence convergent in $B$.

Hence, we get that $y \in P$.
\end{demo}

\bigskip

For the rest of this section, we consider two \cst-algebras $A$ and $B$ and a 
densely defined lower semi-continuous \cst-valued weight $\vfi$ from $A$ into 
$M(B)$. We will also fix a KSGNS-construction $(E,\la,\pi)$ for $\vfi$.

\medskip

\begin{remark} \rm
Define $$P = \{ x \in M(A)^+ \mid \text{The net  } (b^* \rho(x) b)_{\rho \in 
\cG_\vfi} \text{ is convergent for every } b \in B \} . $$
We know by the previous lemma that $P$ is a hereditary cone in $M(A)^+$
such that $\Mfi^+ \subseteq P$.

Lemma \ref{lemA4} implies for every $x \in P$
that the net $(\rho(x))_{\rho \in \cG_\vfi}$ is strictly convergent in $M(B)$.

So we get for every  $x \in \text{span } P$ that
the net $(\rho(x))_{\rho \in \cG_\vfi}$ is strictly convergent in
$M(B)$.

If we define a mapping $\psi$ from $\text{span}\,P$ into $M(B)$
such that the net $(\rho(x))_{\rho \in \cG_\vfi}$ converges strictly to 
$\psi(x)$ for every $x \in \text{span }$, we get a \cst-valued weight $\psi$ 
from $M(A)$ into $M(B)$.

This remarks justify the following definition.
\end{remark}

\begin{definition} \label{def5.2}
We define the \cst-valued weight $\overline{\vfi}$ from $M(A)$ into
$M(B)$ such that
\begin{enumerate}
\item We have that $$\cM_{\overline{\vfi}}^+
= \{ x \in M(A)^+ \mid \text{The net  } (b^* \rho(x) b)_{\rho \in \cG_\vfi} 
\text{ is convergent for every } b \in B \} . $$
\item The net $(\rho(x))_{\rho \in \cG_\vfi}$ converges
strictly to $\overline{\vfi}(x)$ for every $x \in \cM_{\overline{\vfi}}$.
\end{enumerate}
\end{definition}

Because $\vfi$ is lower semi-continuous, it is clear that $\overline{\vfi}$ 
extends $\vfi$.

\medskip

\begin{notation}
We will use the following notations.

\begin{enumerate}
\item We define $\overline{\cM}_\vfi = \cM_{\overline{\vfi}}$ and
$\overline{\cN}_\vfi = \cN_{\overline{\vfi}}$.
\item For every $x \in \overline{\cM}_\vfi$, we put $\vfi(x) = 
\overline{\vfi}(x)$.
\end{enumerate}

It is clear that $\Mfi^+ = \overline{\cM}_\vfi^+ \cap A$
and $\Nfi = \overline{\cN}_\vfi \cap A$.
\end{notation}

\bigskip

The proof of the following results is the same as the proofs of result 
\ref{res3.1} and result  \ref{res3.2}.

\begin{result} \label{res5.1}
Let $(x_i)_{i \in I}$ be a net in $\overline{\cM}_\vfi^+$ and $x$ an element 
in $\overline{\cM}_\vfi^+$ such that $(x_i)_{i \in I}$ converges strictly to 
$x$ and $x_i \leq x$ for every $i \in I$.
Then $(\vfi(x_i))_{i \in I}$ converges strictly to $\vfi(x)$.
\end{result}

\begin{result} \label{res5.2}
Let $(x_i)_{i \in I}$ be a net in $\overline{\cM}_\vfi^+$ and $x$ an element 
in $M(A)^+$ such that $(x_i)_{i \in I}$ converges strictly to $x$ and $x_i 
\leq x$ for every $i \in I$.
Then $x$ belongs to $\overline{\cM}_\vfi^+$ $\Leftrightarrow$ The net $(b^* 
\vfi(x_i) b)_{i \in I}$ is convergent for every $b \in B$.
\end{result}

\bigskip

We want to construct a KSGNS-construction for $\overline{\vfi}$. Therefore, 
we want to extend the map $\la$. This will be done in the next part.

\begin{lemma}
The mapping $\la$ is closable for the strict topology on $M(A)$ and the 
strong topology on $L(B,E)$.
\end{lemma}

This follows immediately from the fact that $\la$ is closed for the strict 
topology on $A$ and the strong topology on $L(B,E)$ (proposition 
\ref{prop3.1}).

\begin{definition} \label{def5.1}
We define $\overline{\la}$ to be the closure of $\la$ with respect to the 
strict topology on $M(A)$ and the strong topology on $L(B,E)$.
For every $a$ in the domain of $\overline{\la}$, we define $\la(a) =
\overline{\la}(a)$.
\end{definition}

\medskip

\begin{proposition}
We have that $\overline{\la}$ is a linear mapping from $\overline{\cN}_\vfi$ 
into $\cL(B,E)$ such that:
\begin{enumerate}
\item We have that $\pi(x) \la(a) = \la(x a)$
for every $x \in M(A)$ and $a \in \overline{\cN}_\vfi$.
\item We have that $\langle \la(a_1) b_1 , \la(a_2) b_2 \rangle = b_2^* \, 
\vfi(a_2^* \, a_1) \, b_1$ for every $a_1,a_2 \in \overline{\cN}_\vfi$ and 
$b_1,b_2 \in B$.
\end{enumerate}
Consequently, the triplet $(E,\overline{\la},\overline{\pi})$ is a KSGNS-
construction for $\overline{\vfi}$.
\end{proposition}
\begin{demo}
We will split the proof up in several parts.
\begin{enumerate}
\item Choose $x \in A$ and $a$ in the domain of $\overline{\la}$.
Then there exists a net $(a_k)_{k \in K}$ in $\Nfi$ such that $(a_k)_{k \in 
K}$ converges strictly to $a$ and $(\la(a_k))_{k \in K}$ converges
strongly to $\la(a)$.

It is clear that $(x a_k)_{k \in K}$ converges to $x a$.
We also have for every $k \in K$ that $x a_k$ belongs to $\Nfi$
and $\la(x a_k) = \pi(x) \la(a_k)$. So we get that
the net $(\la(x a_k))_{k \in K}$ converges strongly to $\pi(x) \la(a)$.

By the norm-strong closedness of $\la$, we see that $x a$ belongs to $\Nfi$ 
and $\la(x a) = \pi(x) \la(a)$.

\item Choose $a$ in the domain of $\overline{\la}$.

Take an approximate unit $(e_k)_{k \in K}$ of $A$. By 1, we know for every $k 
\in K$ that $e_k \, a$ belongs to  $\Nfi$ and $\la(e_k a) = \pi(e_k) \la(a)$. 
From this
we conclude that $a^* e_k^2 \, a$ belongs to $\Mfi^+$
and
$$ b^* \vfi(a^* e_k^2 \, a) b = \langle \la(e_k \, a) b , \la(e_k \, a) b 
\rangle = \langle \pi(e_k) \la(a) b , \pi(e_k) \la(a) b \rangle$$
for every $b \in B$ and $k \in K$.

Hence, we see that the net $(b^* \vfi(a^* e_k^2 \, a) b)_{k \in K}$
converges to $\langle \la(a) b , \la(a) b \rangle$ for every
$b \in B$.

Because we also have that $a^* e_k^2 \, a \leq a^* a$ for every $k \in K$ and 
because the net $(a^* e_k^2 \, a)_{k \in K}$ converges strictly to $a^* a$, 
results \ref{res5.1} and \ref{res5.2} imply
that $a^* a$ belongs to $\overline{\cM}_\vfi^+$
and $b^* \vfi(a^* a) b = \langle \la(a) b , \la(a) b \rangle$ for every $b 
\in B$. Of course, we get also that $a$ belongs to $\overline{\cN}_\vfi$

\item Choose $a \in \overline{\cN}_\vfi$.

Take an approximate unit $(e_k)_{k \in K}$ for $A$. Then $e_k \, a$ belongs 
to $\Nfi$ for every $k \in K$.

We see that $(a^* e_k \, a)_{k \in K}$ is a net in $\Mfi^+$ such that $a^* 
e_k \, a \leq a^* a$ for every $k \in k$. Because $(a^* e_k \, a)_{k \in K}$ 
converges strictly to $a^* a$, result \ref{res5.1} implies that the net 
$(\vfi(a^* e_k \, a))_{k \in K}$ converges strictly to $\vfi(a^* a)$.

Choose $b \in B$.

Take $k,l \in K$ with $l \geq k$. Then
$0 \leq e_l - e_k \leq 1$, which implies that
$0 \leq (e_l - e_k)^2 \leq e_l - e_k$.

Therefore,
\begin{eqnarray*}
\| \la(e_l \, a) b - \la(e_k \,a) b \|^2 & = &
 \| \langle \la((e_l - e_k) \, a) b , \la((e_l - e_k) \, a) b \rangle \|
= \| b^* \vfi(a^* (e_l - e_k)^2 \, a) b \|  \\
& \leq & \| b^* \vfi(a^* (e_l - e_k) \, a) b \|
= \| b^* \vfi(a^* e_l \, a) b - b^* \vfi(a^* e_k \, a) b \| .
\end{eqnarray*}

This implies easily that $(\la(e_k \, a)b)_{k \in K}$ is Cauchy and hence 
convergent in $B$.

From this all, we infer the existence of a linear operator $t$ from $B$ into 
$E$ such that $(\la(e_k \, a))_{k \in K}$ converges strongly to $t$. By the 
definition of $\overline{\la}$, we find that $a$ belongs to the domain of 
$\overline{\la}$.

\end{enumerate}

In the preceding part of the proof, we have proven the following facts.
\begin{description}
\item[\rm \ a)] $\overline{\cN}_\vfi$ equals the domain of $\overline{\la}$
\item[\rm \ b)] We have for every $a \in \overline{\cN}_\vfi$ and $b \in B$ 
that $b^* \vfi(a^* a) b = \langle \la(a) b , \la(a) b \rangle $.
\item[\rm \ c)] For every $x \in A$ and $a \in \overline{\cN}_\vfi$, we have 
that $\la(x a) = \pi(x) \la(a)$.
\end{description}
We still have to proof some minor details. Statement 2 of the proposition 
follows from b) by polarization. Statement 1 of the proposition follows from 
c) by the same method as in the proof of result \ref{res1.1}.3 .

Finally, choose $a \in \overline{\cN}_\vfi$.

We have for every $c \in \Nfi$, $d \in B$ and $b \in B$
that $\langle \la(a)(b) , \la(c) d \rangle
= d^* \vfi(c^* a) \, b$.  Therefore, lemma \ref{prop4.1}  of appendix 2 
implies that $\la(a)$ belongs
to $\cL(B,E)$.
\end{demo}

\begin{result}
We have the following equalities:
\begin{enumerate}
\item We have for every $a_1,a_2 \in \overline{\cN}_\vfi$ that
$\vfi(a_2^* a_1) = \la(a_2)^* \la(a_1)$.
\item For every $a \in \overline{\cN}_\vfi$, we have that $\|\la(a)\|^2 = 
\|\vfi(a^* a)\|$.
\end{enumerate}
\end{result}

\bigskip

\begin{remark} \rm
By definition, we have that $\Nfi$ is a strict-strong core for 
$\overline{\la}$. But we have even more:

Consider $a \in \overline{\cN}_\vfi$. Then there exists a net $(a_k)_{k \in 
K}$ in $\Nfi$ such that
\begin{itemize}
\item We have for every $k \in K$ that $\|a_k\| \leq \|a\|$ and
$\|\la(a_k)\| \leq \|\la(a)\|$.
\item The net $(a_k)_{k \in K}$ converges strictly to $a$ and the net
$(\la(a_k))_{k \in K}$ converges strongly$^*$ to $\la(a)$.
\end{itemize}
This follows immediately by multiplying $a$ to the left by an approximate 
unit of $A$.
\end{remark}

\bigskip

Looking at prosition \ref{prop2.1}, we have the following generalization to 
the multiplier algebra.

\begin{proposition} \label{prop5.1}
Consider $\rho \in \cF_\vfi$. Let $S$ be an element in $\cL(E) \cap \pi(A)'$ 
such that $S^* S = T_\rho$ and let $v$ be the unique element in $\cL(B,E)$ 
such that $S \la(a) = \pi(a) v$ for every $a \in \Nfi$.
Then $S \la(a) = \pi(a) v$ for every $a \in \overline{\cN}_\vfi$.
\end{proposition}
\begin{demo}
We have for every $e \in A$ that $e a$ belongs to $\Nfi$ and $\la(e a) = 
\pi(e) \la(a)$, which implies
that $$\pi(e) S \la(a) = S \pi(e) \la(a) = S \la(e a) = \pi(e a) v =
\pi(e) \pi(a) v.$$

The non-degeneracy of $\pi$ implies that $S \la(a) = \pi(a) v$.
\end{demo}

\begin{corollary} \label{cor5.1}
Consider $\rho \in \cF_\vfi$. Then $T_\rho^\frac{1}{2} \la(a) = \pi(a) 
v_\rho$ for every $a \in \overline{\cN}_\vfi$. We have moreover that
$\la(b)^* T_\rho \, \la(a) = \rho(b^* a)$ for every $a,b \in 
\overline{\cN}_\vfi$.
\end{corollary}

The last statement of this corollary follows easily from the last statement 
of notation \ref{not2.1}.

\bigskip

We also want to mention the following result.

\begin{lemma} \label{lem5.2}
Let $(a_j)_{j \in J}$ be a net in $\overline{\cN}_\vfi$ and an element $a \in 
\overline{\cN}_\vfi$
such that $(a_j)_{j \in J}$ converges strictly to $a$ and $(\la(a_j))_{j \in 
J}$ is bounded.
Then $(\la(a_j)^*)_{j \in J}$ converges strongly to $\la(a)^*$.
\end{lemma}
\begin{demo}
Choose $c \in \Nfi$, $d \in B$ and $\rho \in \cG_\vfi$.

We have for every $j \in J$ that $\la(a_j)^* T_\rho \, \la(c) d = \rho(a_j^* 
c) d$.
Therefore $(\la(a_j)^* T_\rho \, \la(c) d)_{j \in J}$ converges to
$\rho(a^* c)d$, which is equal to $\la(a)^* T_\rho \, \la(c) d$.

Because $(\la(a_j)^*)_{j \in J}$ is also bounded, this implies
that $(\la(a_j)^*)_{j \in J}$ converges strongly to $\la(a)^*$.
\end{demo}

\bigskip

Referring to lemma \ref{lem1.1}, we get the following proposition.

\begin{proposition}
Let $a$ be an element in $\overline{\cN}_\vfi$.
Then $\vfi(a^* a)$ belongs to $B$  $\Leftrightarrow$ There exists an element 
$x \in E$ such that $\la(a) b = x b$ for every $b \in B$ $\Leftrightarrow$
$\la(a)$ belongs to ${\cal K}(B,E)$.
\end{proposition}

As before, the following remark applies. Let $a$ be an element in 
$\overline{\cN}_\vfi$ sucht that $\vfi(a^* a)$ belongs to $B$. Then there 
exists $x \in E$ such that $\la(a) b = x b$ for every $b \in B$. In this 
case, we have that $\vfi(a^* a) = \langle x , x \rangle$.

\begin{proposition} \label{prop5.2}
Let $a$ be an element in $\overline{\cN}_\vfi$ such that $\vfi(a^*a)$ belongs 
to $B$. Then we have for every $\rho \in \cF_\vfi$ that $\rho(a^* a)$ belongs 
to $B$ and the net $(\rho(a^* a))_{\rho \in \cG_\vfi}$ converges to $\vfi(a^* 
a)$.
\end{proposition}
\begin{demo}
Take a KSGNS-construction $(E,\la,\pi)$ for $\vfi$.

By the previous proposition, we know that there exists an element $x \in
E$ such that $\la(a) b = x b$ for every $b \in B$.

Choose $\rho \in \cF_\vfi$. Then we have for every $b \in B$
that
$$ b^* \langle T_\rho x , x \rangle b = \langle T_\rho x b , x b \rangle
= \langle T_\rho \la(a) b , \la(a) b \rangle = b^* \rho(a^* a) b $$
which implies that  $\langle T_\rho x , x \rangle = \rho(a^* a)$.
Hence, we see immediately that $\rho(a^* a)$ belongs to $B$.

We have also that $(\rho(a^* a))_{\rho \in \cG_\vfi}$ converges to $\langle x 
, x \rangle$ which is equal to $\vfi(a^* a)$.
\end{demo}

\begin{corollary}
Let $a$ be an element in $\overline{\cN}_\vfi$. Then $\vfi(a^* a)$ belongs to 
$B$ $\Leftrightarrow$ We have for every $\rho \in \cG_\vfi$ that $\rho(a^* 
a)$ belongs to $B$ and $(\rho(a^* a))_{\rho \in \cG_\vfi}$ is convergent in 
$B$.
\end{corollary}

\section{A first step towards a construction procedure for
C$^*$-valued weights.}  \label{art6}

We will consider the following objects in this section.

Suppose that $A$,\ $B$ are two \cst-algebras and that $E$ is a Hilbert \cst-
module over $B$. Let $N_0$ be a dense subalgebra of $A$ and $\la_0$ a linear 
mapping from $N_0$ into $L(B,E)$ such that
$\langle \, \la_0(a) b \mid a \in N_0 , b \in N \rangle$ is dense in $E$.

Furthermore, we assume the existence of a net $(T_i)_{i \in I}$ in $\cL(E)^+$ 
such that
\begin{itemize}
\item We have for every $i \in I$ that $\|T_i\| \leq 1$.
\item $(T_i)_{i \in I}$ converges strongly to 1.
\item For every $i \in I$ there exists a unique strict completely positive
mapping $\rho_i$ from $A$ into $M(B)$ such that
$b_2^* \, \rho_i(a_2^* \,  a_1) \, b_1 = \langle T_i \la_0(a_1)b_1 , 
\la_0(a_2)b_2 \rangle$ for every $a_1,a_2 \in N_0$ and $b_1,b_2 \in B$.
\end{itemize}

It is clear from proposition \ref{prop4.1} that $\la_0(a)$ belongs to 
$\cL(B,E)$ for every $a \in N_0$.

\begin{lemma} \label{lem6.2}
We have that $\la_0$ is closable for the strict topology on $A$ and the 
strong topology on $L(B,E)$.
\end{lemma}
\begin{demo}
Choose a net $(a_j)_{j \in J}$ in $N_0$, $t \in L(B,E)$ such that
$(a_j)_{j \in J}$ converges strictly to 0 and $(\la_0(a_j))_{j \in J}$ 
converges strongly to $t$. Choose $b,d \in B$ and $c \in N_0$.

\medskip

Fix $i \in I$.
We have clearly that $(\langle T_i \la_0(a_j) b , \la_0(c) d \rangle)_{j \in 
J}$ converges to $\langle T_i \, t(b) , \la_0(c) d \rangle$.

Because $\langle T_i \la_0(a_j) b , \la_0(c) d \rangle = d^* \rho_i(c^* a_j) 
b$ for every $j \in J$, we also have that the net \newline
$(\langle T_i \la_0(a_j) b , \la_0(c) d \rangle)_{j \in J}$ converges to 0. 
This implies that $\langle T_i \, t(b) , \la_0(c)d \rangle = 0$.

\medskip

The fact that $(T_i)_{i \in I}$ converges strongly to 1, implies that
$\langle t(b) , \la_0(c) d \rangle = 0$. Consequently, $t=0$.
\end{demo}

\begin{notation}
We define $\la$ to be the closure of $\la_0$ for the norm topology on $A$ and 
the strong topology on $L(B,E)$. We also define $N$ to be the domain of $\la$.
\end{notation}

\begin{lemma} \label{lem6.3}
Consider $a_1,a_2 \in N$, $b_1,b_2 \in B$ and $i \in I$.
Then $\langle T_i \la(a_1) b_1 , \la(a_2) b_2 \rangle = b_2^* \, \rho_i(a_2^* 
\, a_1) \, b_1$.
\end{lemma}
\begin{demo}
\begin{enumerate}
\item Choose $c \in N_0$. There exists a net $(d_j)_{j \in J}$ in $N_0$ such 
that $(d_j)_{j \in J}$ converges to $a_2$ and $(\la_0(d_j))_{j \in J}$ 
converges strongly to $\la(a_2)$.

We have clearly that the net $(\langle T_i \la_0(c) b_1 , \la_0(d_j) b_2 
\rangle)_{j \in J}$ converges  to $\langle T_i \la_0(c) b_1 , \la(a_2) b_2 
\rangle$.
We have also that $\langle T_i \la_0(c) b_1 , \la_0(d_j) b_2 \rangle = b_2^* 
\, \rho_i(d_j^* \, c) \, b_1$ for every $j \in J$.
This implies that the net $(\langle T_i \la_0(c) b_1 , \la_0(d_j) b_2 
\rangle)_{j \in J}$ converges to $b_2^* \, \rho_i(a_2^* \, c) \, b_1$.

Combining these two results, we get that
$\langle T_i \la_0(c) b_1 , \la(a_2) b_2 \rangle
= b_2^* \rho_i(a_2^* c) b_1$.

\item There exists a net $(c_k)_{k \in K}$ in $N_0$ such that
$(c_k)_{k \in K}$ converges to $a_1$ and $(\la_0(c_k))_{k \in K}$ converges 
strongly to  $\la(a_1)$.

We have clearly that the net $(\langle T_i \la_0(c_k) b_1 , \la(a_2) b_2 
\rangle)_{k \in K}$ converges  to $\langle T_i \la(a_1) b_1 , \la(a_2) b_2 
\rangle$. By the first part of the proof, we have also that $\langle T_i 
\la_0(c_k) b_1 , \la(a_2) b_2 \rangle = b_2^* \, \rho_i(a_2^* \, c_k) \, b_1$ 
for every $k \in K$.
This implies that the net $(\langle T_i \la_0(c_k) b_1 , \la(a_2) b_2 
\rangle)_{j \in J}$ converges to $b_2^* \, \rho_i(a_2^* \, a_1) \, b_1$.

Combining these two results, we get that
$\langle T_i \la(a_1) b_1 , \la(a_2) b_2 \rangle
= b_2^* \, \rho_i(a_2^* \, a_1) \, b_1$.
\end{enumerate}
\end{demo}

\begin{remark} \rm
So, using proposition \ref{prop4.1} once again, we arrive at the following 
conclusion:

\medskip

The set $N$ is a dense subspace of $A$ and $\la$ is a linear mapping from $N$ 
into $\cL(B,E)$ which is closed for the norm topology on $A$
and the strong topology on $L(B,E)$.
Furthermore, $N_0$ is a core for $\la$ in the norm-strong sense.

\medskip

Consider $a_1,a_2 \in N$. Then $\rho_i(a_2^* \, a_1) = \la(a_2)^* T_i \, 
\la(a_1)$ for every $i \in I$.
This implies that $(\rho_i(a_2^* \,  a_1))_{i \in I}$ converges strictly to 
$\la(a_2)^* \la(a_1)$ for every $a_1,a_2 \in N$.
\end{remark}

\medskip

\begin{lemma}
Consider $x \in N_0$, $a_1,\ldots\!,a_n \in N_0$ and $b_1,\ldots\!,b_2 \in 
B$. \newline Then $\| \sum_{k=1}^n \la_0(x a_k) b_k \| \leq \|x\| \, 
\|\sum_{k=1}^n \la_0(a_k) b_k \|$.
\end{lemma}
\begin{demo}
We have for every $i \in I$ that
\begin{eqnarray*}
& & \langle T_i (\sum_{k=1}^n \la_0(x a_k) b_k) , \sum_{k=1}^n \la_0(x a_k) 
b_k \rangle
=  \sum_{k,l=1}^n \langle T_i \la_0(x a_k) b_k , \la_0(x a_l) b_l \rangle \\
& & \spat =  \sum_{k,l=1}^n b_l^* \, \rho_i((x a_l)^* (x a_k)) \, b_k
\leq \|x\|^2 \, \sum_{k,l=1}^n b_l \, \rho_i(a_l^* \, a_k) \, b_k
\hspace{4cm} \text{(*)} \\
& & \spat = \|x\|^2 \, \sum_{k,l=1}^n \langle T_i \la_0(a_k) b_k , \la_0(a_l) 
b_l \rangle
= \|x\|^2 \,\,  \langle T_i (\sum_{k=1}^n \la_0(a_k) b_k) , \sum_{k=1}^n 
\la_0( a_k) b_k \rangle ,
\end{eqnarray*}
where in inequality (*), we used the complete positivity of $\rho_i$.
Because $(T_i)_{i \in I}$ converges strongly to 1, we get that
$$\langle \sum_{k=1}^n \la_0(x a_k) b_k , \sum_{k=1}^n \la_0(x a_k) b_k 
\rangle \leq \|x\|^2 \, \, \langle \sum_{k=1}^n \la_0( a_k) b_k , 
\sum_{k=1}^n \la_0(a_k) b_k \rangle .$$
\end{demo}

This lemma implies for every $x \in N_0$ the existence of a continuous linear 
operator $L_x$ from $E$ into $E$ such that $L_x (\la_0(a) b) = \la_0(x a) b$ 
for every $a \in N_0$ and $b \in B$.  It is also clear that $\|L_x\| \leq 
\|x\|$ for every $x \in N_0$.

\medskip

It is not difficult to check that the mapping $N_0 \rightarrow {\cal B}(E) : 
x \mapsto L_x$ is a continuous algebra-homomorphism. This justifies the 
following definition.

\begin{notation}
We define $\pi$ to be the continuous algebra-homomorphism from $A$ into 
${\cal B}(E)$ such that $\pi(x) \la_0(a) = \la_0(x a)$ for every $x,a \in 
N_0$.
\end{notation}

\medskip

\begin{proposition} \label{prop6.1}
The set $N$ is a left ideal in $A$ and $\pi(x) \la(a) = \la(x a)$
for every $x \in A$ and $a \in N$.
\end{proposition}
\begin{demo}
\begin{enumerate}
\item Choose $x \in N_0$ and $a \in N$.
Then there exists a net $(a_j)_{j \in J}$ in $N_0$ such that $(a_j)_{j \in 
J}$ converges to $a$ and $(\la_0(a_j))_{j \in J}$ converges strongly to 
$\la(a)$. It is clear that $(x a_j)_{j \in J}$ converges to $x a$.

We have also for every $j \in J$ that $x a_j$ belongs to $N_0$ and
$\la_0(x a_j) = \pi(x) \la(a_j)$. This implies that $(\la_0(x a_j))_{j \in 
J}$ converges strongly to $\pi(x) \la(a)$ in ${\cal B}(B,E)$.

By definition, we find that $x a$ belongs to $N$ and $\la(x a) = \pi(x) 
\la(a)$.

\item Choose $x \in A$ and $a \in N$. Then there exists a sequence  
$(x_k)_{k=1}^\infty$ in $N_0$ such that $(x_k)_{k=1}^\infty$ converges to 
$x$. This implies immediately that $(x_k\,a)_{k=1}^\infty$ converges to $x a$.

By the first part of this proof, we know for every $k \in \N$ that $x_k\,a$ 
belongs to $N$ and $\la(x_k\, a) = \pi(x_k) \la(a)$. So we get that $(\la(x_k 
\, a))_{k=1}^\infty$ converges to $\pi(x) \la(a)$.

The norm-strong closedness of $\la$ implies that $x a$ belongs to $N$
and $\la(x a) = \pi(x) \la(a)$.
\end{enumerate}
\end{demo}

\begin{proposition}
The mapping $\pi$ is a $^*$-homomorphism from $A$ into $\cL(E)$.
\end{proposition}
\begin{demo}
Choose $x \in A$.  Take $a,c \in N$ and $b,d \in B$.

We have for every $i \in I$ that
\begin{eqnarray*}
& & \langle T_i \pi(x) \la(a) b , \la(c) d \rangle
= \langle T_i \la(x a) b , \la(c) d \rangle
=  d^* \, \rho_i(c^* x a) \, b \\
& & \spat = d^* \, \rho_i((x^* c)^* a) \, b
=  \langle T_i \la(a) b , \la(x^* c) d \rangle
=  \langle T_i \la(a) b , \pi(x^*) \la(c) d \rangle
\end{eqnarray*}
From this, we infer that
$$\langle \pi(x) \la(a) b ,  \la(c) d \rangle = \langle \la(a) b , \pi(x^*) 
\la(c) d \rangle . $$

This implies that $\langle \pi(x) v , w \rangle = \langle v , \pi(x^*) w $ 
for every $v,w \in E$. So we arrive at the conclusion that $\pi(x)$ belongs 
to $\cL(E)$ and $\pi(x)^* = \pi(x^*)$.
\end{demo}

\begin{proposition}
The $^*$-homomorphism $\pi$ is non-degenerate.
\end{proposition}
\begin{demo}
Choose an approximate unit $(e_j)_{j \in J}$ of $A$.
Then we have for every $j \in J$ that $0 \leq \pi(e_j) \leq 1$.

Take $i \in I$ , $a,c \in N$ and $b,d \in B$.
We have for every $j \in J$ that
\begin{eqnarray*}
\langle \pi(e_j)  T_i \la(a) b , \la(c) d \rangle & = & \langle T_i \la(a) b 
, \pi(e_j) \la(c) d \rangle \\
& = & \langle T_i \la(a) b , \la(e_j c) d \rangle
= d^* \rho_i(c^* e_j a) \, d  .
\end{eqnarray*}
This implies that $(\langle \pi(e_j) T_i \la(a) b , \la(c) d \rangle)_{j \in 
J}$ converges to $d^* \rho_i(c^* a) b$, which is equal to $\langle T_i \la(a) 
b , \la(c) d \rangle$.

\medskip

Because $(\pi(e_j))_{j \in J}$ is bounded, this last result implies that $( 
\langle \pi(e_j) v , w \rangle )_{j \in J}$ converges to
$\langle v , w \rangle$ for every $v,w \in E$.
Using lemma \ref{lemA1}, we see that $(\pi(e_j))_{j \in J}$ converges 
strongly to 1.
\end{demo}

\begin{proposition} \label{prop6.2}
The set $N$ is a left ideal in $M(A)$ and $\la(x a) = \pi(x) \la(a)$
for every $x \in M(A)$ and $a \in N$.
\end{proposition}
\begin{demo}
Choose $x \in M(A)$ and $a \in N$. Take an approximate unit $(e_k)_{k \in K}$ 
for $A$.

Then $(e_k \, x a)_{k \in K}$ converges to $x a$.
By proposition \ref{prop6.1}, we have for every $k \in K$ that $e_k \, x a$ 
belongs to $N$ and
$$\la(e_k \, x a) = \pi(e_k \, x) \la(a) = \pi(e_k) \pi(x) \la(a) \ .$$
Using the previous proposition, this implies that $(\la(e_k \, x a))_{k \in 
K}$ converges strongly to $\pi(x) \la(a)$.

The norm-strong closedness of $\la$ implies that $x a$ belongs to $N$
and $\la(x a) = \pi(x) \la(a)$.
\end{demo}

\begin{proposition} \label{prop6.3}
The mapping $\la$ is closed for the strict topology on $A$ and the strong 
topology on $L(B,E)$.
\end{proposition}
\begin{demo}
Take $a \in A$, $t \in L(B,E)$ such that there exists a net $(a_j)_{j \in J}$ 
in $N_0$ such that $(a_j)_{j \in J}$ converges strictly to $a$ and
$(\la_0(a_j))_{j \in J}$ converges strongly to $t$.

\begin{list}{}{\setlength{\leftmargin}{.4 cm}}

\item Take $e \in N_0$.

It is clear that $(e a_j)_{j \in J}$ converges to $e a$.
We have for every $j \in J$ that $e a_j$ belongs to $N_0$ and
$\la_0(e a_j) = \pi(e) \la_0(a_j)$. This implies that $(\la_0(e a_j))_{j \in 
J}$ converges strongly to $\pi(e) t$.
By definition, we see that $e a$ belongs to $N$ and $\la(e a) = \pi(e) t$.

\end{list}

We know that there exists a bounded net  $(e_k)_{k \in K}$ in $N_0$ such that 
$(e_k)_{k \in K}$ converges strictly to 1.
This implies immediately that $(e_k \, a)_{k \in K}$ converges to $a$.

By the first part of this proof, we know for every $k \in K$ that $e_k \, a$ 
belongs to $N$ and $\la(e_k a) = \pi(e_k) \la(a)$.
This implies that $(\la(e_k \, a))_{j \in J}$ converges strongly to $t$. The 
norm-strong closedness of $\la$ implies that $a$ belongs to $N$ and $\la(a) = 
t$.

\medskip

From this all, we conclude that $\la$ is  equal to the closure of $\la_0$ for 
the strict topology on $A$ and the norm topology on $L(B,E)$.
\end{demo}

\medskip

\begin{remark} \rm
We are now in a position to use the results
and terminology of section 2 using the ingredients $A,B,E,N,\la,\pi$.
So we use the notations $\cF, \cG, \cH$ and so on.

Consider $i \in I$. In the terminology of section 2, we have that $\rho_i$ 
belongs to ${\cal F}$ and $T_{\rho_i} = T_i$.

Therefore $T_i$ belongs to $\pi(A)'$.

Furthermore, there exists a unique element $v_i \in \cL(B,E)$ such that
$T_i^\frac{1}{2} \la(a) = \pi(a) v_i$ for every $a \in N$. We also
know that $\|v_i\|^2 = \|\rho_i\|$ and that $\rho_i(x) = v_i^* \pi(x) v_i$ 
for every $x \in M(A)$.
\end{remark}

\medskip

\begin{proposition} Consider $a_1,a_2 \in N$.
We have that $(\rho(a_2^* \, a_1))_{\rho \in {\cal G}}$ converges strictly to 
$\la(a_2)^* \la(a_1)$.
\end{proposition}
\begin{demo}
Choose $a \in N$. Take $b \in B$.

We have for every $\rho \in \cG$ that $ \rho(a^* a)  \leq \la(a)^* \la(a)$. 
Hence, because $(b^* \rho_i(a^* a) b)_{i \in I}$ converges to $b^* \la(a)^* 
\la(a)b $,  lemma  \ref{lem2.1} implies that $(b^* \rho(a^* a) b)_{\rho \in 
\cG}$ converges to $b^* \la(a)^* \la(a) b$.

Lemma \ref{lemA3} implies that $(\rho(a^* a))_{\rho \in \cG}$ converges 
strictly to $\la(a)^* \la(a)$. The proposition follows by polarisation.
\end{demo}

Similar as in result \ref{res3.3} one proves the following corollary.

\begin{corollary}
The net $(T_\rho)_{\rho \in \cG}$ converges strongly to 1.
\end{corollary}

\medskip

The last lemma is useful for  the material in the next sections.

\begin{lemma}  \label{lem6.1}
Let $u$ be an element in $M(A)$ such that $u N_0 \subseteq N_0$ and
such that there exists an element $S \in \cL(E)$ such that
$S \la_0(a) = \la_0(a u)$ for every $a \in N_0$.
Then $u N \subseteq N$ and $\la(a u) = S \la(a)$ for every $a \in N$.
\end{lemma}
\begin{demo}
Choose $a \in N$. Then there exists a net $(a_j)_{j \in J}$ such that
$(a_j)_{j \in J}$ converges to $a$ and $(\la_0(a_j))_{j \in J}$ converges 
strongly to $\la(a)$.
It is clear that $(a_j u)_{j \in J}$ converges to $a u$.

We have for every $j \in J$ that $a_j u$ belongs to $N_0$ and
$\la_0(a_j u) = S \la_0(a_j)$. This implies that $(\la_0(a_j u))_{j \in J}$ 
converges strongly to $S \la(a)$.

By definition, we get that $a u$ belongs to $N$ and $\la(a u) = S
\la(a)$.
\end{demo}

\section{A construction procedure for regular \cst-valued weights.}
\label{art7}

In this section, we will propose a construction procedure for
regular \cst-valued weights out of some sort of KSGNS-construction which 
generalizes a similar construction procedure in section 2.2 of \cite{Verd}. 
We  will also prove a result (proposition \ref{prop7.4})  which plays a major 
role in the next section.

\bigskip

For the rest of this section, we fix the following ingredients.

\medskip

Consider two \cst-algebras $A$ and $B$ and a Hilbert-\cst-module $E$ over 
$B$. Let $N$ be a dense left ideal of $A$ and $\la$  a linear map from $N$ 
into $L(B,E)$ such that
\begin{itemize}
\item $\la$ is norm-strongly closed.
\item The set $\langle \la(a) b \mid  a \in N, b \in B \rangle$ is dense in 
$E$.
\end{itemize}

Furthermore, we assume the existence of a net $(u_i)_{i \in I}$ in $M(A)$ 
such that we have for every $i \in I$ that
\begin{enumerate}
\item $N u_i \subseteq N$
\item There exists a unique element $S_i \in \cL(E)$ such that $S_i \la(a) = 
\la(a u_i)$ for every $a \in N$.
\item $\|u_i\| \leq 1$ and $\|S_i\| \leq 1$
\item There exists a unique strict completely positive linear mapping 
$\rho_i$ from $A$ into $M(B)$ such that $b_2^* \rho_i(a_2^* a_1) b_1 = 
\langle S_i \la(a_1) b_1, S_i \la(a_2) b_2 \rangle$ for every $a_1,a_2 \in N$ 
and $b_1,b_2 \in B$.
\end{enumerate}
and such that
\begin{enumerate}
\item $(u_i)_{i \in I}$ converges strictly to 1
\item $(S_i)_{i \in I}$ converges strongly to 1
\end{enumerate}

\medskip

Fix $i \in I$ and define $T_i = S_i^* S_i \in \cL(E)$, it is clear
that $0 \leq T_i \leq 1$.

\smallskip

We also have immediately that $b_2^* \rho_i(a_2^* a_1) b_1 = \langle T_i 
\la(a_1) b_1 , \la(a_2) b_2 \rangle$ for every $a_1,a_2 \in N$ and $b_1,b_2 
\in B$.

\medskip

Because $(S_i)_{i \in I}$ converges strongly to 1, we get that $(\langle T_i 
v , v \rangle)_{i \in I}$ converges to $\langle v , v \rangle$ for every $v 
\in E$. Therefore, lemma \ref{lemA1} implies that $(T_i)_{i \in I}$ converges 
strongly to 1.

\medskip

These properties imply that our ingredients $A,B,E,N,\la$ satisfy the 
conditions of the beginning of the previous section, but we don't have to 
close $\la$ anymore. So, we can use the results of the previous section. We 
will give a summary of them.

\begin{itemize}
\item We have for every $a \in N$ that $\la(a)$ belongs to $\cL(B,E)$.
\item $N$ is a left ideal in $M(A)$
\item There exist a unique non-degenerate $^*$-homomorphism $\pi$ from $A$ 
into $\cL(E)$ such that $\pi(x) \la(a) = \la(x a)$ for every $x \in M(A)$ and 
$a \in N$.
\item The mapping $\la$ is closed for the strict topology on $A$ and the
strong topology on $L(B,E)$.
\item It is possible to introduce the objects ${\cal H},{\cal F},{\cal G}$ 
like in section 2. In this terminology, we have for every $i \in I$ that 
$\rho_i$ belongs to ${\cal F}$ and $T_{\rho_i} = \rho_i$.
\item We have for every $i \in I$ that $T_i$ belongs to $\pi(A)'$. Therefore 
there exist a unique $v_i \in \cL(B,E)$ such that $T_i^\frac{1}{2} \la(a) = 
\pi(a) v_i$. Moreover $\|v_i\|^2 = \|\rho_i\|$. Furthermore, $\rho_i(x) = 
v_i^* \pi(x) v_i$ for every $x \in M(A)$.
\item The net $(T_\rho)_{\rho \in {\cal G}}$ converges strongly to 1.
\item We have for every $a_1,a_2 \in N$ that the net $(\rho_i(a_2^* a_1))_{i 
\in I}$ converges strictly to $\la(a_2)^* \la(a_1)$.
\item We have for every $a_1,a_2 \in N$ that the net $(\rho(a_2^* a_1))_{\rho 
\in {\cal G}}$ converges strictly to $\la(a_2)^* \la(a_1)$.
\end{itemize}

\bigskip

In the rest of this section, we want to construct a regular \cst-valued 
weight with these ingredients. At the same time, we will prove some extra 
properties. Essentially, we will define a \cst-valued weight $\vfi$ on $A$ 
such that $\Nfi = N$ and
$\vfi(b^* a) = \la(b)^* \la(a)$ for every $a,b \in N$. There are three 
problems with this:
\begin{enumerate}
\item Is this mapping well defined? If $a_1,\ldots\!,a_n , b_1,\ldots\!,b_n$ 
are elements in $N$ such that $\sum_{k=1}^n b_k^* a_k = 0$, do we have that 
$\sum_{k=1}^n \la(b_k)^* \la(a_k) = 0$ ?
\item Is the positive part of $N^* N$ a hereditary cone in $A^+$ ?
\item Is the resulting $\vfi$ lower semi-continuous?
\end{enumerate}
The first question can be easily answered. If $\sum_{k=1}^n b_k^* a_k = 0$, 
then
$$\sum_{k=1}^n \la(b_k)^* T_i \, \la(a_k) = \rho_i(\sum_{k=1}^n b_k^* a_k) = 
0$$
for every $i \in I$. Because $(T_i)_{i \in I}$ converges strongly to 1, we 
get that $\sum_{k=1}^n \la(b_k)^* \la(a_k) = 0$.

\medskip

The other problems are more difficult to resolve. But this will be done in 
the rest of this section.

\bigskip\medskip

We start of with the following lemma.

\begin{lemma}
Consider $i \in I$. Then $S_i$ belongs to $\pi(A)'$.
\end{lemma}
\begin{demo}
Choose $x \in A$. Take $a \in N$.

By assumption, we have that $a u_i$ belongs to $N$
and $S_i \la(a) = \la(a u_i)$. This implies that $x(a u_i)$ belongs
to $N$ and $\la(x (a u_i)) = \pi(x) \la(a u_i) = \pi(x) S_i \la(a)$.

We know also that $x a$ belongs to $N$. Again, this implies that $ (x a) u_i$ 
belongs to $N$ and $\la( (x a) u_i ) = S_i \la(x a) = S_i \pi(x) \la(a)$.

Comparing these two results, we get that $\pi(x) S_i  \la(a) = S_i \pi(x) 
\la(a)$.

From this, we infer that $\pi(x) S_i = S_i \pi(x)$.
\end{demo}

This lemma allows us to introduce the following notation (see proposition 
\ref{prop2.1}).

\begin{notation} \label{not7.1}
Consider $i \in I$. Then there exists a unique element $w_i \in \cL(B,E)$ 
such that $S_i \la(a) = \pi(a) w_i$ for every $a \in N$. Moreover, we have 
that $\|w_i\|^2 = \|\rho_i\|$ and $\rho_i(x) = w_i^* \pi(x) w_i$ for all $x 
\in M(A)$.
\end{notation}

Using this result, we can prove the following one.

\begin{result} \label{lem7.3}
Consider $i \in I$. Then $A u_i \subseteq N$ and $\la(a u_i) = \pi(a) w_i$ 
for every $a \in A$.
\end{result}
\begin{demo}
Choose $a \in A$. Then there exists a sequence $(a_j)_{j=1}^\infty$ in $N$ 
such that $(a_j)_{j=1}^\infty$ converges to $a$.

We get immediately that $(a_j u_i)_{j=1}^\infty$ converges to $a u_i$. By the 
previous notations, we have for every $j \in \N$ that $a_j u_i$ belongs to 
$N$ and
$\la(a_j u_i) = \pi(a_j) w_i$. This implies that $(\la(a_j 
u_i))_{j=1}^\infty$ converges to $\pi(a) w_i$.

Because $\la$ is norm-strongly closed, we get that $a u_i$ belongs to $N$ and 
$\la(a u_i) = \pi(a) w_i$.
\end{demo}

\bigskip

\begin{remark} \label{rem7.1} \rm
We want to use some Hilbert space theory to prove some equalities. But 
therefore, we will have to make a transition from Hilbert \cst-modules
to Hilbert spaces. This will be done in the following way.
\begin{itemize}
\item Consider a \cst-algebra $C$ and a Hilbert \cst-module $G$ over $C$.
Let $\om$ be an element in $C_+^*$.

First, we define the positive sesquilinear mapping $( \, , \, )$ from $G 
\times G$ into $\C$ such that $(v,w) = \om(\langle v , w \rangle)$ for every 
$v,w \in G$. Then $G , (\, , \,)$ becomes a semi innerproduct space.

We define $K = \{ v \in G \mid (v,v)=0 \}$. Then $K$ is a subspace of $G$
and $\frac{G}{K}$ is in a natural way a innerproduct space.
Then $\overline{G}$ is defined as the completion of $\frac{G}{K}$,
so $\overline{G}$ is a Hilbert space. The inner product on $\overline{G}$ 
will also be denoted by $\langle \, , \, \rangle$.

For every $v \in G$, we define $\overline{v}$ as the equivalence class of $v$ 
in $\frac{G}{N}$. Then we have the following properties
\begin{enumerate}
\item The mapping $G \rightarrow \overline{G} : v \mapsto \overline{v}$
is linear.
\item The set \ $\{ \overline{v} \mid v \in G \}$ is dense in $\overline{G}$.
\item We have for every $v,w \in G$ that $\langle \overline{v} , \overline{w} 
\rangle = \om( \langle v , w \rangle)$.
\end{enumerate}

(If $G = C$, we get nothing else but the GNS-space of $\om$.)

\item Consider a \cst-algebra $C$ and 2 Hilbert \cst-modules $G_1$,$G_2$ over 
$C$. Let $\om$ be an element in $C_+^*$.
We use $\om$ to create 2 Hilbert spaces in the way described above.

Now let $f$ be a linear mapping from $G_1$ into $G_2$
such that there exists a positive number $M$
satisfying $\langle f(v) , f(v) \rangle \leq M \, \langle v , v \rangle$
for every $v \in G_1$. (This is the case for elements in $\cL(G_1,G_2)$.)

Then it is not difficult to see that there exists a unique continuous linear 
map $\overline{f}$ from $\overline{G_1}$ into $\overline{G_2}$ such that 
$\overline{f}(\overline{v}) = \overline{f(v)}$ for every $v \in G_1$. 
Moreover $\|f\| \leq M$.
\end{itemize}
\end{remark}

We will now give a first application of this transition method.
The following two results are inspired by proposition 2.1.15 of \cite{Verd}.

\begin{lemma} \label{lem7.1}
Let $\om \in B_+^*$ and define the Hilbert space $\overline{E}$ with respect 
to this functional $\om$ as described above.
Consider $a \in A$, $b \in B$ such that
$(b^* \rho_i(a^* a) b)_{i \in I}$ is eventually bounded. Then there exists a 
sequence $(a_n)_{n=1}^\infty$ in $N$  and an element $v \in \overline{E}$ 
such that $(a_n)_{n=1}^\infty$ converges to $a$ and 
$(\,\overline{\la(a_n)b}\,)_{n=1}^\infty$ converges to $v$.
\end{lemma}
\begin{demo}
Using result \ref{lem7.3}, we have for every $i \in I$ that $a u_i$ belongs 
to $N$  and $\la(a u_i) = \pi(a) w_i$, which implies that
\begin{eqnarray*}
\|\,\overline{\la(a u_i) b}\,\|^2 & = &
\langle \, \overline{\la(a u_i) b} , \overline{\la(a u_i) b} \, \rangle = 
\om( \langle \la(a u_i) b , \la(a u_i) b \rangle ) \\
& = & \om( \langle \pi(a) w_i b , \pi(a) w_i b \rangle )
= \om(b^* \, \rho_i(a^* a) \, b) \ \ ,
\end{eqnarray*}
where we used notation \ref{not7.1} in the last equality.
Consequently, the net $(\,\overline{\la(a u_i) b}\,)_{i \in I}$ is eventually 
bounded.

Because we also have that the net $(a u_i)_{i \in I}$ converges to $a$, lemma 
\ref{lemA5} implies the existence of a sequence $(a_n)_{n=1}^\infty$ in the 
convex hull of $\{ a u_i \mid i \in I \} \subseteq N$ and $v \in 
\overline{E}$ such
that $(a_n)_{n=1}^\infty$ converges to $a$ and $(\,\overline{\la(a_n) 
b}\,)_{n=1}^\infty$ converges to $v$.
\end{demo}

\begin{proposition} \label{prop7.3}
Consider $a \in M(A)^+$, $b \in B$ and $x  \in M(B)^+$ such that there exists 
an element $i_0 \in I$ such that $b^* \rho_i(a) b \leq x$ for every $i \in I$ 
with $i \geq i_0$.
Then $b^* \rho(a) b \leq x$ for every $\rho \in \cF$.
\end{proposition}
\begin{demo}
Fix $\eta \in \cF$.

Choose $e \in A^+$ with $\|e\| \leq 1$.  Take $\om \in B_+^*$.

We have for every $i \in I$ with $i \geq i_0$ that $b^* \rho_i( a^\frac{1}{2} 
e a^\frac{1}{2}) b \leq  b^* \rho_i (a) b \leq x$.

Therefore, the previous lemma implies the existence of a sequence 
$(c_n)_{n=1}^\infty$  in N and an element $v \in \overline{E}$ such that 
$(c_n)_{n=1}^\infty$ converges to $e^\frac{1}{2} \, a^\frac{1}{2}$
and $(\overline{\la(c_n)b})_{n=1}^\infty$ converges to $v$. \ \ \ \  \ \ \ (a)

We have for every $\rho \in \cF$  and $n \in \N$ that
$$
\langle \, \overline{T_\rho} \, \overline{\la(c_n)b} , \, 
\overline{\la(c_n)b} \rangle  =\langle \, \overline{T_\rho \la(c_n) b} , 
\overline{\la(c_n)b} \, \rangle
 =   \om( \langle T_\rho \la(c_n) b , \la(c_n) b \rangle)
=   \om( b^* \rho(c_n^* c_n) b )  . $$
From this equation and (a), we infer easily that
$$ \langle \overline{T_\rho} v , v \rangle = \om(b^* \rho(a^\frac{1}{2} e 
a^\frac{1}{2}) b)  \text{\ \ \ \ \ \ \ \ \ (b)}$$
for every $\rho \in \cF$.

In particular, we have for every $i \in I$ with $i \geq i_0$ that
$$ \langle \overline{T_i} v , v \rangle
= \om(b^* \rho_i(a^\frac{1}{2} e a^\frac{1}{2}) b )
\leq \om(x).$$
Because $(\overline{T_i})_{i \in I}$ converges strongly to 1
the previous inequality implies that $\|v\|^2 \leq \om(x)$.

Hence, using (b) with $\rho$ equal to $\eta$, we get that
$$\om(b^* \eta(a^\frac{1}{2} e a^\frac{1}{2}) b) =
\langle \overline{T_\eta} v , v \rangle \leq \|v\|^2 \leq \om(x) .$$

\medskip

Because $\om$ was chosen arbitrarily, we conclude that
$b^* \eta(a^\frac{1}{2} e a^\frac{1}{2}) b \leq x$.

\medskip

Next, we take an approximate unit $(e_k)_{k \in K}$ in $A$.
By the previous part, we know that $b^* \eta(a^\frac{1}{2} e_k a^\frac{1}{2}) 
b \leq x$ for every $k \in K$.
The convergence of $(b^* \eta(a^\frac{1}{2} e_k a^\frac{1}{2})b)_{k \in K}$
to $b^* \eta(a) b$ implies that $b^* \eta(a) b \leq x$.

\end{demo}

\begin{proposition} \label{prop7.2}
Consider $a \in M(A)^+$ and $b \in B$ such that $(b^* \rho_i(a) b)_{i \in I}$ 
converges to an element $x \in B^+$.
Then the net $(b^* \rho(a) b)_{\rho \in \cG}$ converges also to $x$.
\end{proposition}
\begin{demo}
First, we prove that $b^* \rho(a^* a) b \leq x$ for every $\rho \in
\cF$.

Therefore, fix $\eta \in \cF$.

Choose $\vep > 0$. Then there exists an element $i_0 \in I$ such that
$\| b^* \rho_i(a^* a) b - x \| \leq \vep$ for every $i \in I$ with $i \geq 
i_0$.
This implies that $b^* \rho_i(a^* a) b \leq x+\vep \, 1$ for every $i \in I$ 
with $i \geq i_0$.

Therefore, the previous proposition guarantees that
$b^* \eta(a^* a) b \leq x+\vep \, 1$.

Because $\vep$ was chosen arbitrarily, this implies that $b^* \eta(a^* a) b 
\leq x$.

\medskip

Now, because $(b^* \rho_i(a) b)_{i \in I}$ converges to $x$, lemma 
\ref{lem2.1} implies that $(b^* \rho(a) b)_{\rho \in \cG}$ converges to $x$.
\end{demo}

\bigskip\bigskip

In the next part we will prove a major result (proposition \ref{prop7.4}) of 
this paper which allows us to resolve both problems mentioned in the 
beginning of this section. At the same time, this result  will be very useful 
in the following section. We will split the proof of proposition 
\ref{prop7.4} up in several parts.

\bigskip

Fix a dense right ideal $D$ of $B$. We define the set $${\cal P} = \{ a \in 
M(A)^+ \mid \text{We have for every } d \in D \text{ that the net } (d^* 
\rho(a) d)_{\rho \in {\cal G}} \text{ is convergent in } B \} .$$  
\label{pag1}

By lemma \ref{lem5.1}, we know that ${\cal P}$ is a hereditary cone in 
$M(A)^+$.

We will denote ${\cal N} = \{ a \in M(A) \mid a^* a \in {\cal P} \}$ and
${\cal M} = {\cal N}^* {\cal N} = \text{span } P$.
As usual, ${\cal N}$ is a left ideal in $M(A)$, ${\cal M}$ is a sub-$^*$-
algebra of $M(A)$ and ${\cal M}^+ = {\cal P}$.

\medskip

The following lemma follows from polarisation.

\begin{lemma} Consider $a_1,a_2 \in {\cal N}$ and $b_1,b_2 \in D$.
Then we have that the net $(b_2^* \rho(a_2^* a_1) b_1)_{\rho \in {\cal G}}$ 
is convergent in $B$.
\end{lemma}

We are now going to construct a Hilbert-\cst-module over $B$ in a similar way 
as we did for the KSGNS-construction for \cst-valued weights.

\medskip

We define the complex vector space $F = {\cal N} \od D$. Moreover
\begin{itemize}
\item We turn $F$ into a right $B$-module such that $(a \ot b) c = a \ot (b 
c)$ for every $a \in {\cal N}$, $b \in D$ and $c \in B$.
\item We turn $F$ into a semi-innerproduct module over $B$
such that $(b_2^* \rho(a_2^* a_1) b_1)_{\rho \in {\cal G}}$ converges to
\newline $\langle a_1 \ot b_1 , a_2 \ot b_2 \rangle$ for every $a_1,a_2 \in 
{\cal N}$ and $b_1,b_2 \in D$.
\end{itemize}
As before, we define $L = \{ x \in F \mid \langle x , x \rangle =  0 \}$,
then $L$ is a submodule of $F$. Then $\frac{F}{L}$ is turned into a 
innerproduct module over $B$. We define $M(A) \otp D$ to be the completion of 
$\frac{F}{L}$, so $M(A) \otp D$ is a Hilbert \cst-module over $B$.

For every $a \in {\cal N}$ and every $b \in D$, we define $a \otp b$ to be 
the equivalence class of $a \ot b$ in $\frac{F}{L}$. Then we have the 
following properties:
\begin{enumerate}
\item The function ${\cal N} \times D \rightarrow M(A) \otp D : (a,b) \mapsto 
a \otp b$ is bilinear.
\item The set $\langle a \otp b \mid a \in {\cal N}, b \in D \rangle$ is 
dense in $M(A) \otp D$.
\item For every $a \in {\cal N}$, $b \in D$ and $c \in B$, we have that
$(a \otp b) c = a \otp (b c)$.
\item For every $a_1,a_2 \in {\cal N}$ and $b_1,b_2 \in D$, we have that
the net $(b_2^* \rho(a_2^* a_1) b_1)_{\rho \in {\cal G}}$ converges to 
\newline $\langle a_1 \otp b_1 , a_2 \otp b_2 \rangle$.
\end{enumerate}

From the 4th property, we have immediately that
\begin{equation}
\sum_{i,j=1}^n b_j^* \, \rho(a_j^* \, a_i) \, b_i
\leq \langle \sum_{i=1}^n a_i \otp b_i , \sum_{i=1}^n a_i \otp b_i \rangle 
\label{eq7.1}
\end{equation}
for every $\rho \in \cG$, $a_1,\dots\!,a_n \in \cN$ and $b_1,\ldots\!,b_n \in 
B$. Of course, this inequality remains true for $\rho \in \cF$.

Also, the 4th property implies that $\langle (x a_1) \otp b_1 , a_2 \otp b_2 
\rangle = \langle a_1 \otp b_1 , (x^* a_2) \otp b_2 \rangle$ for every 
$a_1,a_2 \in {\cal N}$ and $b_1,b_2 \in D$.
We will use this fact in the following lemma.

\begin{lemma}
Consider $x \in M(A)$, $a_1,\ldots\!,a_n \in {\cal N}$ and $b_1,\ldots\!,b_n 
\in D$.
Then
$$ \| \sum_{i=1}^n (x a_i) \otp b_i \| \leq \|x\| \, \| \sum_{i=1}^n a_i \otp 
b_i \| .$$
\end{lemma}
\begin{demo}
By the remark before the lemma, we get immediately that
$$\langle \sum_{i=1}^n (z a_i) \otp b_i , \sum_{i=1}^n (z a_i) \otp b_i 
\rangle = \langle \sum_{i=1}^n (z^* z a_i) \otp b_i , \sum_{i=1}^n a_i \otp 
b_i \rangle $$
for every $z \in M(A)$

We know that there exist an element $y \in M(A)$ such that $\|x\|^2 1 - x^* x 
= y^* y$.
So we see that
\begin{eqnarray*}
& & \|x\|^2 \langle \sum_{i=1}^n  a_i \otp b_i , \sum_{i=1}^n  a_i \otp b_i 
\rangle - \langle \sum_{i=1}^n (x a_i) \otp b_i , \sum_{i=1}^n (x a_i) \otp 
b_i \rangle  \\
& & \spat = \|x\|^2 \langle \sum_{i=1}^n  a_i \otp b_i , \sum_{i=1}^n  a_i 
\otp b_i \rangle - \langle \sum_{i=1}^n (x^* x a_i) \otp b_i , \sum_{i=1}^n 
a_i \otp b_i \rangle  \\
& & \spat =  \langle \sum_{i=1}^n (y^* y a_i) \otp b_i , \sum_{i=1}^n a_i 
\otp b_i \rangle  \\
& & \spat =  \langle \sum_{i=1}^n (y a_i) \otp b_i , \sum_{i=1}^n  (y a_i) 
\otp b_i \rangle \geq 0.
\end{eqnarray*}
This implies that
$$  \|x\|^2 \langle \sum_{i=1}^n  a_i \otp b_i , \sum_{i=1}^n  a_i \otp b_i 
\rangle \geq \langle \sum_{i=1}^n (x a_i) \otp b_i , \sum_{i=1}^n (x a_i) 
\otp b_i \rangle .$$
\end{demo}

Choose $x \in M(A)$. Then there exists a unique continuous linear mapping
$L_x$ from $M(A) \otp D$ into $M(A) \otp D$ such that $L_x(a \otp b) = (x a) 
\otp b$ for every $a \in M(A)$ and $b \in D$.

\medskip

Using the remark before the previous lemma, it follows for an element $x \in 
M(A)$ that $\langle L_x v , w \rangle = \langle v , L_{x^*} w \rangle $ for 
every $v,w \in M(A) \otp D$. This implies that $L_x$ belongs to
$\cL(M(A) \otp D)$ and $L_x^* = L_{x^*}$. Therefore, the following notation 
is justified.

\begin{notation}
We define the mapping $\th$ from $A$ into $\cL(M(A) \otp D)$ such that
$\th(x) (a \otp b) = (x a) \otp b$ for $x \in A$, $a \in {\cal N}$ and
$b \in D$.
Then $\th$ is a $^*$-homomorphism.
\end{notation}

\begin{lemma}
The $^*$-homomorphism $\th$ is non-degenerate and $\th(x) (a \otp b) = (x a) 
\otp b$ for every $a \in {\cal N}$, $b \in D$ and $x \in M(A)$.
\end{lemma}
\begin{demo}
Choose an approximate unit $(e_k)_{k \in K}$ for $A$.

Take $a_1,\ldots\!,a_n \in {\cal N}$, $b_1,\ldots\!,b_n \in D$.
Choose $\vep > 0$.

By definition of the inner-product $\langle \, , \, \rangle$, there exists an 
element $\rho \in {\cal G}$ such that
$$\| \sum_{i,j=1}^n b_j^* \rho(a_j^* a_i) b_i - \langle \sum_{i=1}^n a_i \otp 
b_i , \sum_{i=1}^n a_i \otp b_i \rangle \| \leq \frac{\vep}{2} .$$

Because $\rho$ is strictly continuous on bounded sets, we also have the 
existence of an element $k_0 \in K$ such that
$$\| \sum_{i,j=1}^n b_j^* \rho(a_j^* e_k a_i) b_i -
 \sum_{i,j=1}^n b_j^* \rho(a_j^* a_i) b_i \| \leq \frac{\vep}{2} $$
for every $k \in K$ with $k \geq k_0$.

Choose $l \in K$ with $l \geq k_0$. Combining the previous two results,
we get that
$$\| \sum_{i,j=1}^n b_j^* \rho(a_j^* e_l a_i) b_i - \langle \sum_{i=1}^n a_i 
\otp b_i , \sum_{i=1}^n a_i \otp b_i \rangle \| \leq \vep $$

Inequality \ref{eq7.1} implies that
$$\sum_{i,j=1}^n b_j^* \rho(a_j^* e_l a_i) b_i \leq \langle \sum_{i=1}^n 
(e_l^\frac{1}{2} a_i) \otp b_i ,
\sum_{i=1}^n (e_l^\frac{1}{2} a_i) \otp b_i  \rangle .$$
Therefore, we get that
\begin{eqnarray*}
\sum_{i,j=1}^n b_j^* \rho(a_j^* e_l a_i) b_i
& \leq & \langle \th(e_l^\frac{1}{2}) (\sum_{i=1}^n a_i \otp b_i) ,
\th(e_l^\frac{1}{2}) (\sum_{i=1}^n a_i \otp b_i) \rangle \\
&  = & \langle \th(e_l) (\sum_{i=1}^n a_i \otp b_i) ,
\sum_{i=1}^n a_i \otp b_i \rangle  \\
& \leq &\langle \sum_{i=1}^n a_i \otp b_i ,
\sum_{i=1}^n a_i \otp b_i \rangle .
\end{eqnarray*}
From this all, we infer that
\begin{eqnarray*}
& &  \| \langle \th(e_l) (\sum_{i=1}^n a_i \otp b_i) ,
\sum_{i=1}^n a_i \otp b_i \rangle  -
\langle \sum_{i=1}^n a_i \otp b_i ,
\sum_{i=1}^n a_i \otp b_i \rangle  \| \\
& & \spat \leq \| \sum_{i,j=1}^n b_j^* \rho(a_j^* e_l a_i) b_i -
 \langle \sum_{i=1}^n a_i \otp b_i ,
\sum_{i=1}^n a_i \otp b_i \rangle  \| \leq \vep .
\end{eqnarray*}
Hence, we see that the net $(\langle \th(e_k) (\sum_{i=1}^n a_i \otp b_i) , 
\sum_{i=1}^n a_i \otp b_i \rangle )_{k \in K}$ converges to
$ \langle \sum_{i=1}^n a_i \otp b_i , \sum_{i=1}^n a_i \otp b_i \rangle$.

Because $(\th(e_k))_{k \in K}$ is bounded, this implies that
$(\langle \th(e_k) v , v \rangle)_{k \in K}$ converges to $\langle v , v 
\rangle$ for every $v \in M(A) \otp D$. Using lemma \ref{lemA1} once more, we 
get
that $(\th(e_k))_{k \in K}$ converges strongly to 1. So $\th$ is non-
degenerate.

\medskip

It is easy to check that the mapping $M(A) \rightarrow \cL(M(A) \otp D) : x 
\mapsto L_x$ is a $^*$-homomorphism extending $\th$.
By unicity of the extension, we must have that $\th(x) = L_x$ for every $x 
\in M(A)$
\end{demo}

From the beginning of this section. We know that $(\rho(a_2^* a_1))_{\rho \in 
{\cal G}}$ converges strictly to $\la(a_2)^* \la(a_1)$ for every $a_1,a_2 \in 
N$.

This implies immediately that $N$ is a subset of ${\cal N}$ and
$\langle \la(a_1) b_1 , \la(a_2) b_2 \rangle = \langle a_1 \otp b_1 , a_2 
\otp b_2 \rangle$ for every $a_1,a_2 \in N$ and $b_1,b_2 \in D$.
From this, we get the existence of a unique isometry $U$ from $E$ into $M(A) 
\otp D$ such that $U (\la(a) b) = a \otp b$ for every
$a \in N$ and $d \in D$. \label{pag2}

Then $U$ is $B$-linear and $\langle U v , Uw \rangle = \langle v , w \rangle 
$ for every $v,w \in E$. Ultimately, we will prove that $U$ is also 
surjective and this will be the main result of this section.

\bigskip

We want to prove the equality in lemma \ref{lem7.5} (shortly after we have 
proven it, it will become clear why we need it). In order to do so, we make a 
transition from Hilbert \cst-module theory to Hilbert space theory. For the 
next part, we will fix a positive linear functional $\om$ on $B$ and 
introduce Hilbert spaces with respect to $\om$ as described in remark 
\ref{rem7.1}. The same notations will be used.

\medskip

\begin{lemma}
Consider $a_1,\ldots\!,a_m \, , \, c_1,\ldots\!,c_n \in {\cal N}$ and
$b_1,\ldots\!,b_m \, , \, d_1,\ldots\!,d_n \in D$. Let $\rho$ be an element  
\newline in  ${\cal F}$.
Then
$$|\sum_{j=1}^m \sum_{k=1}^n \om(d_k^* \rho(c_k^* a_j) b_j)| \leq 
\|\sum_{j=1}^m \overline{a_j \otp b_j} \| \, \| \sum_{k=1}^n \overline{c_k 
\otp d_k} \| . $$
\end{lemma}
\begin{demo}
It is possible to define a sesquilinear mapping $s$ from $M(A) \od B$ into 
$\C$ such that $s(a \ot b , c \ot d) = \om(d^* \rho(c^* a) b)$ for every $a,c 
\in M(A)$ and $b,d \in B$.

It is easy to check that $s$ is positive. Therefore, the Cauchy-Schwarz 
inequality implies that
\begin{eqnarray*}
|\sum_{j=1}^m \sum_{k=1}^n \om(d_k^* \rho(c_k^* a_j) b_j)|^2
& = & | s(\sum_{j=1}^m a_j \ot b_j , \sum_{k=1}^n c_k \ot d_k )|^2 \\
& \leq &  s(\sum_{j=1}^m a_j \ot b_j,\sum_{j=1}^m a_j \ot b_j) \, \,
s(\sum_{k=1}^n c_k \ot d_k,\sum_{k=1}^n c_k \ot d_k) \\
& =  &  \sum_{j,k=1}^m \om(b_k^* \rho(a_k^* a_j) b_j) \,
\sum_{j,k=1}^n \om(d_k^* \rho(c_k^* c_j) d_j) \\
& \leq & \om(\langle \sum_{j=1}^m a_j \otp b_j ,
\sum_{j=1}^m a_j \otp b_j \rangle) \, \,
\om(\langle \sum_{k=1}^n c_k \otp d_k ,
\sum_{k=1}^n c_k \otp d_k \rangle) \hspace{1cm} (*) \\
& = & \| \sum_{j=1}^m \overline{a_j \otp b_j} \|^2 \, \,
\| \sum_{k=1}^n \overline{c_k \otp d_k} \|^2
\end{eqnarray*}
\end{demo}
where we used inequality \ref{eq7.1} in inequality (*).
Consider $\rho \in {\cal F}$. By the previous lemma, there exist a unique
continuous positive sesquilinear form $s_\rho$ on $\overline{M(A) \otp D}$ 
such that $s_\rho(\overline{a \otp b} , \overline{c \otp d}) = \om(d^* 
\rho(c^* a) b )$ for every $a,b \in {\cal N}$ and $b,d \in D$.
The previous lemma also implies that $\|s_\rho\| \leq 1$.
These remarks justify the following notation.

\begin{notation}
Consider $\rho \in {\cal F}$. Then there exists a unique element $D_\rho$ in 
${\cal B}(\overline{M(A) \otp D})$ such that
$\langle D_\rho \overline{a \otp b} , \overline{c \otp d} \rangle
= \om(d^* \rho(c^* a) b) $ for every $a,c \in {\cal N}$ and $b,d \in D$.
We also have that $0 \leq D_\rho \leq 1$.
\end{notation}

We have for every $a,c \in {\cal N}$ and $b,d \in D$ that the net
$(\langle D_\rho \overline{a \otp b} , \overline{c \otp d}\rangle)_{\rho \in 
{\cal G}}$ converges to $\om(\langle a \otp d, c \otp  d \rangle)$
and this last expression equals $\langle \overline{a \otp b} , \overline{c 
\otp d} \rangle$.

Because $(D_\rho)_{\rho \in {\cal G}}$ is bounded, this implies that
$(\langle D_\rho v ,v \rangle)_{\rho \in {\cal G}}$ converges to $\langle v , 
v \rangle$ for every $v \in \overline{A \otp B}$.

Referring to lemma \ref{lemA1}, we get that $(D_\rho)_{\rho \in {\cal G}}$ 
converges strongly to 1.

\begin{lemma}
Consider a bounded net $(a_k)_{k \in K}$ in ${\cN}$, $a \in \cN$, $b \in D$ 
and $v \in \overline{M(A) \otp D}$ such that
$(a_k)_{k \in K}$ converges strictly to $a$ and $(\overline{a_k \otp b})_{k 
\in K}$ converges to $v$. Then $\overline{a \otp b} = v$.
\end{lemma}
\begin{demo}
Choose $\rho \in \cF$. We have for every $k \in K$
that
\begin{eqnarray*}
\| D_\rho^\frac{1}{2}\,\overline{a_k \otp b} - D_\rho^\frac{1}{2} 
\,\overline{a \otp b} \|^2 & = &
\| D_\rho^\frac{1}{2}\, \overline{(a_k - a) \otp b} \|^2
= \langle D_\rho\, \overline{(a_k - a) \otp b} , \overline{(a_k - a) \otp b} 
\rangle  \\
& = & \om(b^* \, \rho((a_k-a)^* (a_k-a)) \, b) .
\end{eqnarray*}
This implies that $(D_\rho^\frac{1}{2} \,\overline{a_k \otp b} )_{k \in K}$ 
converges to $D_\rho^\frac{1}{2} \, \overline{a \otp b}$.

It is also clear that $(D_\rho^\frac{1}{2} \, \overline{a_k \otp b})_{k \in 
K}$ converges to $D_\rho^\frac{1}{2} v$. Therefore, $D_\rho^\frac{1}{2} \, 
\overline{a \otp b} = D_\rho^\frac{1}{2} v$.

Because $(D_\rho^\frac{1}{2})_{\rho \in  \cG}$ converges strongly to 1, we 
get that $\overline{a \otp b} = v$.
\end{demo}

\begin{proposition}
We have that $\overline{U}$ is a unitary transformation from $\overline{E}$ 
to $\overline{M(A) \otp D}$.
\end{proposition}
\begin{demo}
Because $\langle U x , U x \rangle = \langle x , x \rangle$ for every $x \in 
E$, we get that $\overline{U}$ is isometric. We turn to the surjectivity of 
$\overline{U}$.

\begin{enumerate}
\item Choose $a \in \cN \cap A$ and $b \in D$.

Inequality \ref{eq7.1} implies for every $i \in I$ that $b^* \rho_i(a^* a) b 
\leq \langle a \otp b , a \otp b \rangle$.

Therefore, lemma \ref{lem7.1} implies the existence of a sequence 
$(a_n)_{n=1}^\infty$ in $N$ and $v \in \overline{E}$ such
that $(a_n)_{n=1}^\infty$ converges to $a$ and $(\overline{\la(a_n) b})_{n 
\in N}$ converges to $v$.

We have for every $n \in \N$ that $$ \overline{U} \, \overline{\la(a_n) b} = 
\overline{ U \la(a_n) b } = \overline{a_n \otp b} , $$
which implies that $(\overline{a_n \otp b})_{n=1}^\infty$ converges to 
$\overline{U} v$.

Using the previous lemma, we see that $\overline{a \otp b} = \overline{U} v $.
\item Next, we choose an element $a \in \cN$.
Take an approximate unit $(e_k)_{k \in K}$ in $A$. We have for every $k \in 
K$ that $e_k \, a$ belongs to $\cN \cap A$, so $\overline{e_k \, a \otp b}$ 
belongs to $\text{Ran }\overline{U}$ by the first part of the proof.

We have for every $k \in K$ that $e_k \, a \otp b =  \th(e_k) (a \otp b)$. 
Therefore, the non-degeneracy of $\th$ implies that $(e_k \, a \otp b)_{k \in 
K}$ converges to $a \otp b$. So we get that
$(\overline{e_k \, a \otp b})_{k \in K}$ converges to $\overline{a \otp b}$.

Because $\text{Ran }\overline{U}$ is closed in $\overline{M(A) \otp D}$, we 
conclude that $\overline{a \otp b}$ belongs to $\text{Ran }\overline{U}$.
\end{enumerate}

Because $\text{Ran }\overline{U}$ is closed in $\overline{M(A) \otp D}$ and 
the set $\langle \overline{a \otp b} \mid a \in \cN , b \in D \rangle$ is 
dense in $\overline{M(A) \otp D}$, we get that $\text{Ran }\overline{U} = 
\overline{M(A) \otp D}$.
\end{demo}

\begin{lemma}
Consider $a \in \cN$, $b \in D$, $c \in M(A)$,$d \in B$ and $\rho \in \cF$. 
\newline
Then $\om(d^* \rho(c^* a)\, b) = \om (\langle a \otp b , U T_\rho^\frac{1}{2} 
\pi(c) v_\rho d \rangle )$.
\end{lemma}
\begin{demo}
It is enough to prove the equality for $c \in N$ and $d \in D$.
We have for every $a_1,a_2 \in N$, $b_1,b_2 \in B$ that
\begin{eqnarray*}
\langle \, \overline{T_\rho} \, \overline{\la(a_1) b_1} , \overline{\la(a_2) 
b_2} \, \rangle &  = &
\langle \, \overline{T_\rho \la(a_1) b_1} , \overline{\la(a_2) b_2} \, \rangle
= \om( \langle T_\rho \la(a_1) b_1 , \la(a_2) b_2 \rangle )  \\
& = & \om( b_2^* \rho(a_2^* \, a_1) \, b_1)
= \langle D_\rho \, \overline{a_1 \otp b_1} , \overline{a_2 \otp b_2} \, 
\rangle \\
& = & \langle D_\rho \, \overline{U \la(a_1) b_1} , \overline{U \la(a_2) b_2} 
\, \rangle
= \langle D_\rho \, \overline{U} \, \overline{\la(a_1) b_1} , \overline{U} \, 
\overline{\la(a_2) b_2} \, \rangle \\
& = & \langle \overline{U}^* D_\rho \overline{U} \, \overline{\la(a_1) b_1} , 
\overline{\la(a_2) b_2} \, \rangle
\end{eqnarray*}
This implies that $\overline{T_\rho} = \overline{U}^* D_\rho \overline{U}$.

Therefore, we have that
\begin{eqnarray*}
\om(d^* \rho(c^* a) b) & = &
\langle D_\rho \, \overline{a \otp b} , \overline{c \otp d} \, \rangle
= \langle \, \overline{a \otp b} , D_\rho \, \overline{c \otp d} \, \rangle \\
& = & \langle \, \overline{a \otp b} , \overline{U} \, \overline{T_\rho} \, 
\overline{U}^* \,
\overline{c \otp d} \, \rangle
= \langle \, \overline{a \otp b} , \overline{U} \, \overline{T_\rho} \, 
\overline{\la(c) d} \, \rangle \\
& = & \langle \, \overline{a \otp b} , \overline{U T_\rho \la(c) d} \, \rangle
= \om(\langle a \otp b , U T_\rho \la(c) d \rangle) \\
& = & \om(\langle a \otp b , U T_\rho^\frac{1}{2} \pi(c) v_\rho d \rangle)
\end{eqnarray*}
where we used notation \ref{not2.1} in the last equality.
\end{demo}

Remembering that this lemma is true for every $\om \in B_+^*$, we get the 
lemma which we wanted to prove.

\begin{lemma} \label{lem7.5}
Consider $a \in \cN$, $b \in D$, $c \in M(A)$, $d \in B$ and $\rho \in \cF$. 
\newline Then $d^* \rho(c^* a)\, b = \langle a \otp b , U T_\rho^\frac{1}{2} 
\pi(c) v_\rho d \rangle $.
\end{lemma}

Before using this lemma, we need to introduce some extra noations.

\begin{lemma}
Consider $\rho \in {\cal F}$, $a_1,\ldots\!,a_n \in {\cal N}$ and
$b_1,\ldots\!,b_n \in D$. Then
$$\langle \sum_{j=1}^n \pi(a_j) v_\rho b_j , \sum_{j=1}^n \pi(a_j) v_\rho b_j 
\rangle \leq  \langle \sum_{j=1}^n a_j \otp b_j , \sum_{j=1}^n a_j \otp b_j 
\rangle .$$
\end{lemma}
\begin{demo}
We have that
$$ \langle \sum_{j=1}^n \pi(a_j) v_\rho b_j , \sum_{j=1}^n \pi(a_j) v_\rho 
b_j \rangle   = \sum_{j,k=1}^n b_k^* \rho(a_k^* \, a_j) b_l
\leq \langle \sum_{j=1}^n a_j \otp b_j , \sum_{j=1}^n a_j \otp b_j \rangle . 
$$
\end{demo}

This lemma justifies the following notation.

\begin{notation}
Consider $\rho \in {\cal F}$. We define the continuous linear mapping  
$F_\rho$ from $M(A) \otp D$ into $E$ such that $F_\rho(a \otp b) = \pi(a) 
v_\rho b$ for every $a \in {\cal N}$ and $b \in D$.
We have that $F_\rho$ is $B$-linear and $\langle F_\rho v , F_\rho v \rangle 
\leq \langle v , v \rangle$ for every $v \in M(A) \otp D$.
\end{notation}

\medskip

\begin{definition}
Consider $\rho \in \cF$. We define $R_\rho = U T_\rho^\frac{1}{2} F_\rho$. 
Therefore, $R_\rho$ is a continuous $B$-linear mapping from $M(A) \otp D$ 
into $M(A) \otp D$ such that $\|R_\rho\| \leq 1$.
\end{definition}

Now we will use lemma \ref{lem7.5} to prove that $R_\rho$ belongs to 
$\cL(M(A) \otp D)$.

\begin{lemma} \label{lem7.2}
Consider $\rho \in \cF$. We have for every $a_1,a_2 \in \cN$ and $b_1,b_2 \in 
D$ that $\langle R_\rho a_1 \otp b_1 , a_2 \otp b_2 \rangle = b_2^* \, 
\rho(a_2^* \, a_1) \, b_1$.
This implies that $R_\rho$ belongs to $\cL(M(A) \otp D)$ and $0 \leq R_\rho 
\leq 1$.
\end{lemma}
\begin{demo}
By lemma \ref{lem7.5}, we have for every $a_1,a_2 \in \cN$ and $b_1,b_2 \in 
D$ that
\begin{eqnarray*}
b_1^* \, \rho(a_1^* \, a_2) \, b_2 & = & \langle a_2 \otp b_2 , U 
T_\rho^\frac{1}{2} \pi(a_1) v_\rho b_1 \rangle \\
& = & \langle a_2 \otp b_2 , U T_\rho^\frac{1}{2} F_\rho \, a_1 \otp b_1 
\rangle \\
& = & \langle a_2 \otp b_2 , R_\rho \, a_1 \otp b_1 \rangle .
\end{eqnarray*}
Taking the adjoint of this equation, gives the equality
$$\langle R_\rho \, a_1 \otp b_1 , a_2 \otp b_2 \rangle
= b_2^* \, \rho(a_2^* \, a_1) \, b_1 .$$
Using the complete positivity of $\rho$, this  equality implies that
$\langle R_\rho v , v \rangle \geq 0$ for every $v \in M(A) \otp D$.

This equality implies that $R_\rho$ belongs to $\cL(M(A) \otp D)$ and
$R_\rho \geq 0$.
\end{demo}

\begin{lemma}
We have that $(R_\rho)_{\rho \in \cG}$ converges strongly to 1.
\end{lemma}
\begin{demo}
Choose $a,c \in \cN$ and $b,d \in D$. The definition of $\langle \, , \, 
\rangle$ and the previous lemma imply that
$(\langle R_\rho a \otp b , c \otp d \rangle)_{\rho \in \cG}$ converges
to $\langle a \otp b , c \otp d \rangle$.

Because $(R_\rho)_{\rho \in \cG}$ is bounded, this implies that $(\langle 
R_\rho v , w \rangle)_{\rho \in \cG}$ converges to $\langle v , w \rangle$ 
for every $v,w \in M(A) \otp D$

Lemma \ref{lemA1} guarantees that $(R_\rho)_{\rho \in \cG}$ converges 
strongly to 1.
\end{demo}

At last, we can prove that $U$ is a unitary transformation from $E$ to $M(A) 
\otp D$.

\begin{proposition} \label{prop7.4}
We have that $U$ is a unitary tansformation from $E$ to $M(A) \otp D$.
\end{proposition}
\begin{demo}
Take $v \in M(A) \otp D$.
For every $\rho \in \cG$, the definition of $R_\rho$ guarantees that  $R_\rho 
v$ belongs to $\text{Ran }U$.
The previous lemma implies that $(R_\rho v)_{\rho \in \cG}$ converges to $v$. 
Because $\text{Ran }U$ is closed in $M(A) \otp D$, we get that $v$ belongs to 
$\text{Ran }U$.
\end{demo}

\begin{lemma}
\begin{itemize}
\item Consider $a \in M(A)$. Then $U^* \th(a) U = \pi(a)$.
\item Consider $\rho \in \cG$. Then $U^* R_\rho U = T_\rho$.
\end{itemize}
\end{lemma}
\begin{demo}
\begin{itemize}
\item We have for every $b \in N$ and $c \in D$ that
\begin{eqnarray*}
U^* \th(a) U \la(b) c & = & U^* \th(a) (b \otp c) = U^* (a b \otp c) \\
& = & \la(a b) c = \pi(a) \la(b) c .
\end{eqnarray*}
This implies that $U^* \th(a) U = \pi(a)$.
\item We have for every $a_1,a_2 \in N$ and $b_1,b_2 \in D$ that
\begin{eqnarray*}
\langle U^* R_\rho U \la(a_1) b_1 , \la(a_2) b_2 \rangle & = &
\langle R_\rho U \la(a_1) b_1 , U \la(a_2) b_2 \rangle \\
& = &\langle R_\rho (a_1 \otp b_1) , a_2 \otp b_2 \rangle
=  b_2^* \, \rho(a_2^* \, a_1) \, b_1  \text{\ \ \ \ \ \ \ \ \ (*)} \\
& = & \langle T_\rho \la(a_1) b_1 , \la(a_2) b_2 \rangle ,
\end{eqnarray*}
where in equality (*), we used lemma \ref{lem7.2}.

This implies that $U^* R_\rho U = T_\rho$.
\end{itemize}
\end{demo}

\begin{lemma} \label{lem7.4}
Consider $\rho \in \cF$ and $S \in \pi(A)' \cap \cL(E)$ such that $S^* S = 
T_\rho$. Let $v$ be the unique element in $\cL(B,E)$ such that
$S \la(c) = \pi(c) v$ for every $c \in N$.
Then $S U^* (a \otp b) = \pi(a) v b$ for every $a \in \cN$ and $b \in D$.
\end{lemma}
\begin{demo}
Choose $a \in \cN$ and $b \in D$.

\medskip

Using the previous lemma and lemma \ref{lem7.2}, we find  for every $c \in 
\cN \cap A$ that
\begin{eqnarray*}
\| S U^* (c \otp b)\|^2
& = & \| \langle S U^* (c \otp b) , S U^* (c \otp b) \rangle \|
= \|\langle T_\rho U^*(c \otp b) , U^* (c \otp b) \rangle \| \\
& = & \| \langle U T_\rho U^* (c \otp b) , c \otp b \rangle \|
= \| \langle R_\rho (c \otp b) , c \otp b \rangle \| \\
& = & \| b^* \rho(c^* c) b \| \leq \|b\|^2 \, \|c\|^2 \, \|\rho\| .
\end{eqnarray*}
So we get that the mapping $\cN \cap A \rightarrow E : c \mapsto S U^* (c 
\otp b)$ is continuous.
We also have for every $c \in N$ that $S U^* (c \otp b) = S \la(c) b = \pi(c) 
v b$. Because $N$ is dense in $A$, we get that
$S U^* (c \otp b) = \pi(c) v b$ for every $c \in \cN \cap A$.

\medskip

Take an approximate unit $(e_k)_{k \in K}$ for $A$. We have for every $k \in 
K$ that $e_k \, a$ belongs to $\cN \cap A$, which implies
that $S U^* (e_k \, a \otp b) = \pi(e_k \, a) v b$ by the first part of the 
proof.

So we get for every $k \in K$ that
\begin{eqnarray*}
\pi(e_k) (S U^* (a \otp b)) & = & S \pi(e_k) U^* (a \otp b)
= S U^* \th(e_k) (a \otp b) \\
& = & S U^* (e_k \, a \otp b) = \pi(e_k \, a) v b
= \pi(e_k) ( \pi(a) v b ) .
\end{eqnarray*}
The non-degeneracy of $\pi$ implies that $S U^* (a \otp b) = \pi(a) v b$.
\end{demo}

As a special case of the previous lemma, we get the following one.

\begin{lemma} \label{lem7.6}
Consider $\rho \in \cF$. Then $T_\rho^\frac{1}{2} U^* (a \otp b) = \pi(a) 
v_\rho b$ for every $a \in \cN$, $b \in D$.
\end{lemma}

\bigskip\bigskip

Now we will apply this theory for the first time.

\begin{proposition} \label{prop7.1}
Consider $a \in A$. Then $a$ belongs to $N$ $\Leftrightarrow$
The net $(b^* \rho(a^* a) b)_{\rho \in \cG}$ is convergent for every $b \in 
B$.
\end{proposition}
\begin{demo}
One implication is trivial. We prove the other one.

Therefore, suppose that $(b^* \rho(a^* a) b)_{\rho \in \cG}$ is convergent 
for every $b \in B$. We will apply the previous theory  (which starts at page 
\pageref{pag1}) with $D = B$.
By the definition of $\cN$, it is clear that $a$ belongs to $\cN$.

Choose $b \in B$.
By result \ref{lem7.3} and lemma \ref{lem7.4}, we have for every $i \in I$ 
that $a u_i$ belongs to $N$ and $$\la(a u_i) b = \pi(a) w_i b = S_i U^* (a 
\otp b) $$
(The mapping $U$ was introduced at page \pageref{pag2}).

This implies that $(\la(a u_i) b)_{\rho \in \cG}$ converges to
$U^* (a \otp b)$.

It is also clear that $(a u_i)_{i \in i}$ converges to $a$.
Therefore, the norm-strong closedness of
$\la$ implies that $a$ belongs to $N$ (and $\la(a) b = U^*(a \otp b)$ for 
every $b \in B$).
\end{demo}

The remarks at the beginning of this section guarantee that the following 
mapping $\vfi$ is well defined:

\begin{definition} \label{def7.1}
We define the linear mapping $\vfi$ from $N^* N$ into $M(B)$
such that $\vfi(\sum_{j=1}^n b_j^* \, a_j) = \sum_{j=1}^n \la(b_j)^* 
\la(a_j)$ for every $n \in \N$ and all $a_1,\ldots\!,a_n \in N$.
\end{definition}

\begin{proposition}
We have that $\vfi$ is a densely defined \cst-valued weight from $A$ to $B$ 
such that $\Nfi = N$ and $\vfi(b^* a) = \la(b)^* \la(a)$ for every $a,b \in 
N$.
\end{proposition}
\begin{demo}
Define
$$P = \{ x  \in A^+ \mid \text{We have for every } b \in B
\text{ that the net } (b^* \rho(a^* a) b)_{\rho \in \cG} \text{ is convergent 
in } B \} .$$
Then $P$ is an hereditary cone in $A^+$  (this follows from lemma 
\ref{lem5.1}).
From proposition \ref{prop7.1}, we know that $N = \{ a \in A^+ \mid a^* a \in 
P \}$. Define $M = \text{span } P = N^* N$. Then $\vfi$ is a linear mapping 
from $M$ into $M(B)$.

\medskip

We have for every $n \in \N$ and all $a_1,\ldots\!,a_n \in N$ and
$b_1,\ldots\!,b_n \in B$ that
\begin{eqnarray*}
& &\sum_{j,k=1}^n b_k^* \, \vfi(a_k^* \, a_j) \, b_j
= \sum_{j,k=1}^n b_k^* \la(a_k)^* \la(a_j) b_j \\
& & \spat = \sum_{j,k=1}^n  \langle \la(a_j) b_j , \la(a_k) b_k \rangle = 
\langle \sum_{j=1}^n \la(a_j) b_j ,
\sum_{j=1}^n \la(a_j) b_j \rangle \geq 0 .
\end{eqnarray*}

So, we see that $\vfi$ is a \cst-valued weight from $A$ to $B$ such that 
$\Nfi = N$. We have also immediately that $\vfi$ is densely defined.
\end{demo}

\begin{corollary}  \label{corol7.2}
The triplet $(E,\la,\pi)$ is a KSGNS-construction for $\vfi$.
This implies that $\cH_\vfi = \cH$, $\cF_\vfi = \cF$ and $\cG_\vfi=\cG$.
\end{corollary}

The last statement of this corollary follows immediately from notation
\ref{not3.1}.

\medskip

From the beginning of this section, we know for every $a \in N$ that 
$(\rho(a^* a))_{\rho \in \cG}$ converges strictly to $\la(a)^*
\la(a)$ which is equal to $\vfi(a^* a)$ by definition.
Referring to proposition \ref{prop7.1}, the definition of $\vfi$ and 
definition \ref{def3.1}, we have also:

\begin{corollary}
The \cst-valued weight $\vfi$ is lower semi-continuous.
\end{corollary}

Of course, we have also:

\begin{corollary} \label{corol7.1}
The \cst-valued weight $\vfi$ is regular and has $(u_i)_{i \in I}$ as 
truncating net.
\end{corollary}

In the next section, we will use the theory, introduced in this section, even 
more.

\section{Regular \cst-valued weights.} \label{art8}

In this section, we will prove some results about regular \cst-valued 
weights. In fact, we have already proven the most difficult steps in the 
previous section.

\bigskip

We will work with two \cst-algebras $A$ and $B$
and a regular \cst-valued weight $\vfi$ from $A$ into $M(B)$.
Take a KSGNS-construction $(E,\la,\pi)$ for $\vfi$.

\medskip

We will also fix a truncating net $(u_i)_{i \in I}$ for $\vfi$. First,
we recapitulate some notations and properties.

Choose $i \in I$.
 \begin{itemize}
\item We define $S_i$ to be the unique element in $\cL(E) \cap \pi(A)'$
such that $S_i \la(a) = \la(a u_i)$ for every $a \in \Nfi$.
Furthermore, we put $T_i = S_i^* S_i  \in \cL(E)^+ \cap \pi(A)'$.
\item We have that $\|u_i\| \leq 1$, $\|S_i\| \leq 1$ and $0 \leq T_i \leq 1$.
\item Define the strict completely mapping $\rho_i$ from $A$ into $M(B)$ such 
that $\rho_i(a_2^* a_1) = \la(a_2)^* T_i \la(a_1)$ for every $a_1,a_2 \in A$.
Then, $\rho_i$ belongs to $\cF_\vfi$ and $T_{\rho_i} = T_i$.
\item We define $v_i,w_i \in \cL(B,E)$ such that $\pi(a) v_i = 
T_i^\frac{1}{2} \la(a)$ and $\pi(a) w_i = S_i \la(a)$ for every $a \in \Nfi$. 
\newline
Then the equality $\|v_i\|^2 = \|w_i\|^2 = \|\rho_i\|$ holds. We have also 
that $\rho_i(x) = v_i^* \pi(x) v_i = w_i^* \pi(x) w_i$ for every $x \in 
M(A)$. This implies that $v_i^* v_i = w_i^* w_i$.
\item We have also that $(u_i)_{i \in I}$ converges strictly to 1
and that $(S_i)_{i \in I}$, $(T_i)_{i \in I}$ converge strongly to 1.
\end{itemize}

\medskip

Lemma \ref{lem7.3} of the previous section implies that

\begin{result}
Consider $i \in I$. We have that $A u_i \subseteq \Nfi$ and $\la(a u_i) = 
\pi(a) w_i$ for every $a \in A$.
\end{result}

This result can be generalized to

\begin{result}
Consider $i \in I$. We have that $M(A) u_i \subseteq \overline{\cN}_\vfi$
and $\la(a u_i) = \pi(a) w_i$ for every $a \in M(A)$.
\end{result}
\begin{demo}
Choose $a \in M(A)$. Then there exists a bounded net $(a_k)_{k \in K}$
in $A$ such that $(a_k)_{k \in K}$ converges strictly to $a$.
Then $(a_k u_i)_{k \in K}$ converges strictly to $a u_i$.

By the previous result, we have for every $k \in K$ that
$a_k u_i$ belongs to $\Nfi$ and $\la(a_k u_i) = \pi(a_k) w_i$.
This implies that $(\la(a_k u_i))_{k \in K}$ converges strongly to
$\pi(a) w_i$.

The strictly-strongly closedness of $\overline{\la}$ implies that
$a u_i$ belongs to $\overline{\cN}_\vfi$ and $\la(a u_i) = \pi(a) w_i$ 
(definition \ref{def5.1}).
\end{demo}

\begin{corollary}
Consider $i \in I$. Then $u_i$ belongs to $\overline{\cN}_\vfi$
and $\la(u_i) = w_i$. Therefore, $\vfi(u_i^* u_i) = w_i^* w_i = v_i^* v_i$.
\end{corollary}

Proposition \ref{prop5.1} implies the following one.

\begin{proposition}
Consider $i \in I$. We have for every $a \in \overline{\cN}_\vfi$ that
$T_i^\frac{1}{2} \la(a) = \pi(a) v_i$ and
$S_i \la(a) = \pi(a) w_i = \la(a u_i)$
We have for every $a_1,a_2 \in \overline{\cN}_\vfi$ that
$\rho_i(a_2^* \, a_1) = \la(a_2)^* T_i \la(a_1)$.
\end{proposition}

The next proposition is a generalization of proposition 2.1.15 of \cite{Verd} 
for \cst-valued weights.

\begin{proposition} \label{prop8.2}
\begin{enumerate}
\item Consider $a \in M(A)^+$. Then $a$ belongs to $\overline{\cM}_\vfi^+$
$\Leftrightarrow$ The net $(b^* \rho_i(a) b)_{i \in I}$ is convergent in $B$ 
for every $b \in B$.
\item Consider $a \in \overline{\cM}_\vfi$. Then $(\rho_i(a))_{i \in I}$ 
converges strictly to $\vfi(a)$.
\end{enumerate}
\end{proposition}
\begin{demo}
\begin{itemize}
\item Choose $a \in M(A)^+$ such that $(b^* \rho_i(a) b)_{i \in I}$ is 
convergent in $B$ for every $b \in B$. Proposition \ref{prop7.2} of the 
previous section implies that $(b^* \rho(a) b)_{\rho \in \cG_\vfi}$ is 
convergent for every $b \in B$. By definition \ref{def5.2}, we get that
$a$ belongs to $\overline{\cM}_\vfi^+$.
\item Choose $a_1,a_2 \in \overline{\cN}_\vfi$.
By the previous result, we have for every $i \in I$ that
$\rho_i(a_2^* \, a_1) = \la(a_2)^* T_i \la(a_1)$
This implies that $(\rho_i(a_2^* \, a_1))_{i \in I}$ converges strictly to 
$\la(a_2)^* \la(a_1) = \vfi(a_2^* \, a_1)$.
\end{itemize}
These two results imply the proposition.
\end{demo}

\begin{corollary}
Consider $a \in A^+$. Then $a$ belongs to $\Mfi^+$
$\Leftrightarrow$ The net $(b^* \rho_i(a) b)_{i \in I}$ is convergent in $B$ 
for every $b \in B$.
\end{corollary}

\medskip

The proof of the following proposition is similar to the proof of proposition 
\ref{prop5.2}.

\begin{proposition}
Let $a$ be an element in $\overline{\cN}_\vfi$ such that $\vfi(a^*a)$ belongs 
to $B$. Then we have for every $i \in I$ that $\rho_i(a^* a)$ belongs to $B$ 
and the net $(\rho_i(a^* a))_{i \in I}$ converges to $\vfi(a^* a)$.
\end{proposition}

\begin{corollary}
Let $a$ be an element in $\overline{\cN}_\vfi$. Then $\vfi(a^* a)$ belongs to 
$B$ $\Leftrightarrow$ We have for every $i \in I$ that $\rho_i(a^* a)$ 
belongs to $B$ and $(\rho_i(a^* a))_{i \in I}$ is convergent in $B$.
\end{corollary}

\bigskip

In many cases, we do not really work with $\Nfi$, but rather with a core of 
$\la$. The following propositions give certain possible cores for $\la$. We 
show that they even behave better than just being a core.

\begin{proposition} \label{prop8.3}
Consider a dense subspace  $K$ of $A$ and define the set $L = \langle a u_i 
\mid a \in K, i \in I \rangle$. \newline Let $a \in \Nfi$ such that $\la(a) 
\neq 0$. Then there exists a  net $(a_j)_{j \in J}$ in $L$ such that 
\begin{enumerate}
\item $\|a_j\| \leq \|a\|$ and $\|\la(a_j)\| \leq \|\la(a)\|$ for every $j 
\in J$.
\item $(a_j)_{j \in J}$ converges to $a$ and $(\la(a_j))_{j \in J}$ converges 
strongly$^*$ to $\la(a)$.
\end{enumerate}
\end{proposition}
\begin{demo}
From the beginning of this section, we know for every $i \in I$ that $\pi(a) 
w_i = \la(a u_i) = S_i \la(a)$,
which implies that $(\pi(a) w_i)_{i \in I}$ converges strongly to $\la(a)$. 
Because $\la(a) \neq 0$, this implies the existence of $i_0 \in I$ such that 
$\pi(a) w_i \neq 0$ for every $i \in I$ with $i \geq i_0$.
Define $I_0 = \{ i \in I \mid i \geq i_0 \}$.

\medskip

For the moment, fix $i \in I_0$. There exists a sequence $(c_n)_{n=1}^\infty$ 
in $K$ such that $(c_n)_{n=1}^\infty$ converges to $a$.
Because $a \neq 0$ and $\pi(a) w_i \neq 0$, there exist $n_0 \in \N$ such 
that we have for all $n \in \N$ with $n \geq n_0$ that $c_n \neq 0$ and
$\pi(c_n) w_i \neq 0$.

\medskip

For every $n \in \N$ with $n \geq n_0$, we define $\lambda_n = 
\frac{\|a\|}{\|c_n\|}$ and $\mu_n = \frac{\|\pi(a) w_i\|}{\|\pi(c_n) w_i\|}$.
Then we have that the sequences $(\lambda_n)_{n=n_0}^\infty$ and 
$(\mu_n)_{n=n_0}^\infty$ converge to 1.
For every $n \in \N$ with $n \geq n_0$, we define $d_n = 
\min(\lambda_n,\mu_n) \,\,  c_n$, so $d_n$ belongs to $K$.

We get immediately that $(d_n)_{n=n_0}^\infty$ converges to $a$.
Moreover, we have for every $n \in \N$  with $n \geq n_0$ that
$$\|d_n\| \leq \lambda_n \, \|c_n\| = \|a\| \hspace{1.25cm} \text{ and } 
\hspace{1.25cm} \|\pi(d_n)w_i\| \leq \mu_n \, \|\pi(c_n) w_i\| = \|\pi(a) 
w_i\| = \|S_i \la(a)\| \leq \|\la(a)\| $$

This implies for every $i \in I_0$ and $m \in \N$ the existence of an element 
$b_{m,i} \in K$
such that $\|b_{m,i}\| \leq \|a\|$, $\|\pi(b_{m,i}) w_i\| \leq \|\la(a)\|$ 
and $\|b_{m,i} - a\| \leq
\frac{1}{m} \frac{1}{\|w_i\|+1}$.

\medskip

We define $J=\N \times I_0$ and put on $J$ the product ordering.
For every $j=(m,i) \in J$, we put $a_j = b_{m,i} \, u_i \in L$.
Then we get for every $j=(m,i) \in J$ that
$$\|a_j\| \leq \|b_{m,i}\| \, \|u_i\| \leq \|a\| \hspace{1.5cm} \text{
 and  } \hspace{1.5cm} \|\la(a_j)\| = \|\la(b_{m,i} u_i)\| = \| \pi(b_{m,i}) 
w_i\| \leq \| \la(a) \| \ . $$

Now we are going to prove the convergence properties.

\begin{enumerate}
\item Choose $\vep > 0$. Then
\begin{itemize}
\item There exists $i_1 \in I_0$ such that $\|a u_i - a\| \leq 
\frac{\vep}{2}$ for every $i \in I_0$ with $i \geq i_1$.
\item There exists an element $m_1 \in \N$ such that $\frac{1}{m_1} \leq 
\frac{\vep}{2}$.
\end{itemize}
Therefore, we have for every $j=(m,i) \in J$ with $j \geq (m_1,i_1)$ that
\begin{eqnarray*}
\|a_j - a \| & = & \|b_{m,i} \,  u_i - a\| \leq  \|b_{m,i} \, u_i - a u_i\| + 
\|a u_i - a \| \\
& \leq & \|b_{m,i} - a \| \, \|u_i\| + \frac{\vep}{2}
 \leq  \frac{1}{m} + \frac{\vep}{2} \leq \frac{1}{m_1} + \frac{\vep}{2} \leq 
\vep .
\end{eqnarray*}
Hence, we see that $(a_j)_{j \in J}$ converges to $a$.
\item Choose $c \in B$. Take $\vep > 0$. Then
\begin{itemize}
\item There exists $i_1 \in I_0$ such that $\|S_i \la(a) c  - \la(a) c\| \leq 
\frac{\vep}{2}$ for every $i \in I_0$ with $i \geq i_1$.
\item There exists an element $m_1 \in \N$ such that $\frac{1}{m_1} \leq 
\frac{\vep}{2} \frac{1}{\|c\|+1}$.
\end{itemize}
Therefore, we have for every $j=(m,i) \in  J$ with $j \geq (m_1,i_1)$ that
\begin{eqnarray*}
\|\la(a_j) c - \la(a) c\| & = & \| \la(b_{m,i}  \, u_i) c - \la(a) c \|  \leq 
\| \la(b_{m,i} \, u_i)c - \la(a u_i) c \| + \|\la(a u_i) c - \la(a) c\| \\
& = & \| \pi(b_{m,i}) w_i c - \pi(a) w_i c \| + \| S_i \la(a) c - \la(a) c\|
\leq  \| b_{m,i} - a \| \, \|w_i\| \, \|c\| + \frac{\vep}{2} \\
& \leq & \frac{1}{m} \frac{1}{\|w_i\|+1} \|w_i\| \, \|c\| + \frac{\vep}{2} 
\leq \frac{1}{m_1} \|c\| + \frac{\vep}{2}  \leq \vep .
\end{eqnarray*}
This implies that $(\la(a_j))_{j \in J}$ converges strongly to $\la(a)$
\item Because $(\la(a_j))_{j \in J}$ is bounded and $(a_j)_{j \in J}$ 
converges to $a$, lemma \ref{lem5.2} implies that $(\la(a_j))_{j \in J}$ 
converges strongly$^*$ to $\la(a)$.
\end{enumerate}
\end{demo}

If $\la(a) = 0$, we can conclude something which is a little bit weaker:

\begin{proposition}
Consider a dense subspace of $K$ of $A$ and define the set $L = \langle a u_i 
\mid a \in K, i \in I \rangle$. \newline Let $a \in \Nfi$ and $M$ a positive 
number such that $\|\la(a)\| < M$. Then there exists a  net $(a_j)_{j \in J}$ 
in $L$ such that \begin{enumerate}
\item $\|a_j\| \leq \|a\|$ and $\|\la(a_j)\| \leq M$ for every $j \in J$.
\item $(a_j)_{j \in J}$ converges to $a$ and $(\la(a_j))_{j \in J}$ converges 
strongly$^*$ to $\la(a)$.
\end{enumerate}
\end{proposition}
\begin{demo}
Because of the previous proposition, we only have to deal with the case that 
that $\la(a)=0$. We can also assume that $a \neq 0$ (If $a=0$, the 
proposition is trivially true).

We have for every $i \in I$ that $\pi(a) w_i = \la(a u_i) = S_i \la(a) =0$.

\medskip

For the moment, fix $i \in I$. There exists a sequence $(c_n)_{n=1}^\infty$ 
in $K$ such that $(c_n)_{n=1}^\infty$ converges to $a$.
Because $a \neq 0$ and $\pi(a) w_i=0$, there exist $n_0 \in \N$ such that we 
have for all $n \in \N$ with $n \geq n_0$ that $c_n \neq 0$ and
$\|\pi(c_n) w_i\| \leq M$.

\medskip

For every $n \in \N$ with $n \geq n_0$, we define $\lambda_n = 
\frac{\|a\|}{\|c_n\|}$.
Then we have that the sequence $(\lambda_n)_{n=n_0}^\infty$ converges to 1.

For every $n \in \N$ with $n \geq n_0$, we define $d_n = \min(\lambda_n,1) 
\,\, c_n$, so $d_n$ belongs to $K$.

We get immediately that $(d_n)_{n=n_0}^\infty$ converges to $a$.
Moreover, we have for every $n \in \N$  with $n \geq n_0$ that
$$\|d_n\| \leq \lambda_n \, \|c_n\| = \|a\| \hspace{1.5cm} \text{ and  } 
\hspace{1.5cm} \|\pi(d_n)w_i\| \leq \|\pi(c_n) w_i\| \leq M  \ .$$

This implies for every $i \in I_0$ and $m \in \N$ the existence of an element 
$b_{m,i} \in K$
such that $\|b_{m,i}\| \leq \|a\|$, $\|\pi(b_{m,i}) w_i\| \leq M$ and 
$\|b_{m,i} - a\| \leq
\frac{1}{m} \frac{1}{\|w_i\|+1}$.

\medskip

We define $J=\N \times I$ and put on $J$ the product ordering.
For every $j=(m,i) \in J$, we put $a_j = b_{m,i} \, u_i \in L$.
We have immediately for every $j=(m,i) \in J$ that
$$\|a_j\| \leq \|b_{m,i}\| \, \|u_i\| \leq \|a\| \hspace{1.5cm} \text{ and } 
\hspace{1.5cm} \|\la(a_j)\| = \|\la(b_{m,i} u_i)\| = \| \pi(b_{m,i}) w_i\| 
\leq M \ . $$

Similarly as in the previous proposition, one proves the convergence 
properties.
\end{demo}

We would like to prove similar properties for elements in 
$\overline{\cN}_\vfi$.

\begin{proposition}
Consider a dense subspace of $K$ of $A$ and define the set $L = \langle a u_i 
\mid a \in K, i \in I \rangle$. \newline Let $a \in \overline{\cN}_\vfi$ such 
that $\la(a) \neq 0$. Then there exists a  net $(a_j)_{j \in J}$ in $L$ such 
that \begin{enumerate}
\item $\|a_j\| \leq \|a\|$ and $\|\la(a_j)\| \leq \|\la(a)\|$ for every $j 
\in J$.
\item $(a_j)_{j \in J}$ converges strictly to $a$ and $(\la(a_j))_{j \in J}$ 
converges strongly$^*$ to $\la(a)$.
\end{enumerate}
\end{proposition}
\begin{demo}
Take an approximate unit $(e_k)_{k \in K}$ for $A$. We have for every $k \in 
K$ that $e_k \, a$ belongs to $\Nfi$ and $\la(e_k \, a) = \pi(e_k) \la(a)$.

This implies that $(\la(e_k \, a))_{k \in K}$ converges strongly to $\la(a)$. 
Because $\la(a) \neq 0$, we get the existence of an element $k_0$ in $K$ such 
that $\la(e_k \, a) \neq 0$ for every $k \in K$ with $k \geq k_0$. Put $K_0 = 
\{ k \in K \mid k \geq k_0 \}$.

\medskip

We define $J = \{ (k,n,F) \mid k \in K_0 , \, n \in \N , \, F \text{ a finite 
subset of } B \}$. We put on $J$ an order $\leq$ such that  we have for every 
$j_1 = (k_1,n_1,F_1) \in K$ and $j_2 = (k_2,n_2,F_2) \in K$ that
$$j_1 \leq j_2 \hspace{1.5cm} \Leftrightarrow  \hspace{1.5cm} k_1 \leq k_2\ 
,\  n_1 \leq n_2 \text{\ and \ }
F_1 \subseteq F_2 .$$
In this way, $J$ becomes a directed set.

\medskip

Choose $j=(k,n,F) \in J$.
By proposition \ref{prop8.3}, there exists an element $a_j \in L$ such that
\begin{itemize}
\item $\|a_j\| \leq \|e_k \, a\|$ and $\|\la(a_j)\| \leq \|\la(e_k \, a)\|$.
\item We have that $\|a_j - a e_k \| \leq \frac{1}{n}$.
\item For every $b \in F$, we have that $\|\la(a_j) b - \la(a \, e_k) b \| 
\leq \frac{1}{n}$.
\end{itemize}
Then we have immediately that
$$ \|a_j\| \leq \|a\| \hspace{1.5cm} \text{ and } \hspace{1.5cm}
\|\la(a_j)\| \leq  \| \pi(e_k) \la(a) \| \leq \|\la(a)\| \ . $$

Now we turn to the convergence properties.

\begin{enumerate}
\item Choose $c \in A$. Take $\vep > 0$. Then
\begin{itemize}
\item There exists $k_1 \in K_0$ such that $\| e_k \, a c - a c \| \leq 
\frac{\vep}{2}$  for every $k \in K_0$ with $k \geq k_1$.
\item There exists $n_1 \in \N$ such that $\frac{1}{n_1} \leq 
\frac{\vep}{2(\|c\|+1)}$.
\end{itemize}
Put $j_1 = (k_1,n_1,\emptyset) \in J$. Then we have for every $j=(k,n,F) \in 
J$ with $j \geq j_1$ that $k \geq k_1$ and $n \geq n_1$, therefore
\begin{eqnarray*}
\| a_j \, c - a c \| & \leq & \| a_j \, c - e_k \, a c \| + \| e_k \, a c - a 
c \| \leq \| a_j - e_k \, a \| \, \|c\| + \frac{\vep}{2} \\
& \leq & \frac{1}{n} \, \|c\| + \frac{\vep}{2} \leq \frac{1}{n_1} \, \|c\| + 
\frac{\vep}{2} \leq \vep
\end{eqnarray*}
Hence we see that $(a_j \, c)_{j \in J}$ converges to $a c$.
Similarly, one proves that $(c a_j)_{j \in J}$ converges to $c a$.
\item Choose $b \in B$. Take $\vep > 0$. Then
\begin{itemize}
\item There exists $k_1 \in K_0$ such that $\| \pi(e_k) \la(a) b - \la(a) b 
\| \leq \frac{\vep}{2}$ for every $k \in K_0$ with $k \geq k_1$.
\item There exists $n_1 \in \N$ such that $\frac{1}{n_1} \leq \frac{\vep}{2}$.
\end{itemize}
Put $j_1 = (k_1,n_1,\{b\}) \in J$. Then we have for every $j=(k,n,F) \in J$ 
with $j \geq j_1$ that $k \geq k_1$, $n_1 \geq n$ and $b \in F$, therefore
\begin{eqnarray*}
\| \la(a_j) b - \la(a) b \|  & \leq & \| \la(a_j) b - \la(e_k \, a) b \|+ \| 
\la(e_k \, a) b - \la(a) b \|  \\ & \leq & \frac{1}{n} + \| \pi(e_k) \la(a) b 
- \la(a) b \| \leq  \frac{1}{n_1}  + \frac{\vep}{2}  \leq \frac{\vep}{2} + 
\frac{\vep}{2} = \vep .
\end{eqnarray*}
From this, we get that $(\la(a_j))_{j \in J}$ converges strongly to $\la(a)$.
\item Because $(\la(a_j))_{j \in J}$ is bounded and $(a_j)_{j \in J}$ 
converges strictly to $a$, lemma \ref{lem5.2} implies that \newline 
$(\la(a_j))_{j \in J}$ converges strongly$^*$ to $\la(a)$.
\end{enumerate}
\end{demo}

The following proposition can be proven in a similar way.

\begin{proposition}
Consider a dense subspace of $K$ of $A$ and define the set $L = \langle a u_i 
\mid a \in K, i \in I \rangle$. \newline Let $a \in \overline{\cN}_\vfi$ and 
$M$ a positive number such that $\|\la(a)\| < M$. Then there exists a  net 
$(a_j)_{j \in J}$ in $L$ such that \begin{enumerate}
\item $\|a_j\| \leq \|a\|$ and $\|\la(a_j)\| \leq M$ for every $j \in J$.
\item $(a_j)_{j \in J}$ converges strictly to $a$ and $(\la(a_j))_{j \in J}$ 
converges strongly$^*$ to $\la(a)$.
\end{enumerate}
\end{proposition}

\bigskip\bigskip

So far, we have extended our weight to $\overline{\vfi}$ in such a way that 
$\overline{\vfi}$ takes values in $M(B)$. However, in some cases it is 
interesting to extend $\vfi$ even further  and let it take values in the set 
of elements affiliated to $B$. Now we will take a first step in this 
direction.

\begin{lemma}
Consider $a \in M(A)^+$ and define the set
$$D = \{ b \in B \mid (b^* \rho(a) b)_{\rho \in \cG_\vfi} \text{ is 
convergent in } B \} .$$
Then $D$ is a right ideal in $M(B)$.
\end{lemma}
\begin{demo}
It is immediately clear that $0$ belongs to $D$, that $D$ is closed under
scalar multiplication and that $D$ is closed under multiplication from the 
right with elements of $M(B)$. The addition requires a little bit of 
explanation

Choose $b,c \in D$.
As usual, we get for any $\rho,\eta \in \cG_\vfi$ with $\rho \leq \eta$ that
$$ (b+c)^* (\eta - \rho)(a) (b+c)
\leq 2 \, b^* (\eta - \rho)(a) b + 2 \, c^* (\eta-\rho)(a) c . $$
This implies that $(\,(b+c)^* \rho(a) (b+c)\,)_{\rho \in \cG_\vfi}$ is Cauchy 
and hence convergent in $B$. Therefore, $b+c$ belongs to $D$.
\end{demo}

\begin{definition} \label{def8.2}
We define $\tNfi$ to be the set of elements $a$ in $M(A)$
such that the set
$$\{ b \in B \mid (b^* \rho(a^* a) b)_{\rho \in \cG_\vfi} \text{ is 
convergent in } B \} $$
is dense in $B$
\end{definition}

\begin{remark} \rm
It is easy to see that $\overline{\cN}_\vfi$ is a subset of $\tNfi$.

\medskip

Choose $a \in \overline{\cN}_\vfi$. Take $b \in B$.
Then we have for every $\rho \in \cG_\vfi$ that $\pi(a) v_\rho b = 
T_\rho^\frac{1}{2} \la(a) b$ (corollary \ref{cor5.1}). This implies that 
$(\pi(a) v_\rho b)_{\rho \in \cG_\vfi}$ converges to $\la(a) b$.

\medskip

Next we want to define for every $a \in \tNfi$ an operator $\la(a)$ from $B$ 
into $E$ (which is not necessarily everywhere defined). The previous 
discussion will justify the following definition.
\end{remark}

\medskip

\begin{definition} \label{def8.1}
Consider $a \in \tNfi$. Then we define the mapping $\la(a)$ from within $B$ 
into $E$ such that:
\begin{enumerate}
\item The domain of $\la(a)$ is equal to
$\{ b \in B \mid \text{The net } (\pi(a) v_\rho b)_{\rho \in \cG_\vfi}
\text{ is convergent in } E \} $.
\item We have for  every $b \in D(\la(a))$ that
$(\pi(a) v_\rho b)_{\rho \in \cG_\vfi}$ converges to $\la(a)(b)$.
\end{enumerate}
It is not difficult to check that $\la(a)$ is a $B$-linear map from within 
$B$ into $E$
\end{definition}

The following result is rather nice and depends on the results of the 
previous section.

\begin{proposition} \label{prop8.1}
Consider $a \in \tNfi$. Then
\begin{enumerate}
\item The domain of $\la(a)$ is equal to
$\{ b \in B \mid (b^* \rho(a^* a) b)_{\rho \in \cG_\vfi} \text{ is convergent 
in } B \} \ $ .
\item Let $\rho \in \cF_\vfi$ and $S \in \pi(A)' \cap \cL(E)$ such that
$S^* S = T_\rho$. Define $v$ to be the unique element in $\cL(B,E)$ such
that $S \la(c) = \pi(c) v$ for every $c \in \Nfi$.
Then $S \la(a) b = \pi(a) v b$ for every $b \in D(\la(a))$.
\end{enumerate}
\end{proposition}
\begin{demo}
We define the set
$D = \{ b \in B \mid (b^* \rho(a^* a) b)_{\rho \in \cG_\vfi} \text{ is 
convergent in } B \} \ $.
Because $a$ belongs to $\tNfi$, we have that $D$ is a dense
right ideal in $B$.

Therefore, we can use the construction of $M(A) \otp D$ from the previous 
section (and which begins at page \pageref{pag1}). We will use the notations 
of  that section. The mapping $U$ is introduced at page \pageref{pag2}.

By definition, $a$ belongs to $\cN$.

\begin{itemize}
\item Choose $b \in D(\la(a))$ . By notation \ref{not2.1}, we have for every 
$\rho \in \cG_\vfi$ that $b^* \rho(a^* a) b = \langle \pi(a) v_\rho b , 
\pi(a) v_\rho b \rangle$.
This implies immediately that $(b^* \rho(a^* a) b)_{\rho \in \cG_\vfi}$
is convergent in $B$. Therefore, $b$ belongs to $D$.
\item Choose $b \in D$. Then we can look at $a \otp b \in M(A) \otp D$. By 
lemma \ref{lem7.6}, we know hat
$\pi(a) v_\rho b = T_\rho^\frac{1}{2} U^* (a \otp b) $
for every $\rho \in \cG_\vfi$.
This implies that $(\pi(a) v_\rho b)_{\rho \in \cG_\vfi}$ converges
to $U^* (a \otp b)$.
By definition, we get that $b$ belongs to $D(\la(a))$ and $\la(a) b
= U^* (a \otp b)$.
\end{itemize}

So we get that $D(\la(a)) = D$ and $\la(a) b = U^*(a \otp b)$
for every $b \in D$. Therefore, statement 2) of the proposition follows from 
lemma \ref{lem7.4}.
\end{demo}

\begin{remark} \rm
By the second statement of the proposition, we have for every $a \in \tNfi$ 
that the mapping $D(\la(a)) \rightarrow E : b \mapsto S \la(a) b$
is continuous.
\end{remark}

\medskip

\begin{corollary}  \label{corol8.1}
Consider $a \in \tNfi$. Then we have for every $\rho \in \cG_\vfi$ and $b \in 
D(\la(a))$  \newline that
$T_\rho^\frac{1}{2} \la(a) b = \pi(a) v_\rho b$.
\end{corollary}

\medskip

We have also the following result.

\begin{corollary}
Consider $a \in M(A)$. Then $a$ belongs to $\tNfi$
\newline $\Leftrightarrow$ the set \
$\{ b \in B \mid \text{The net } (\pi(a) v_\rho b)_{\rho \in \cG_\vfi}
\text{ is convergent in } E \} $ \ is dense in $B$.
\end{corollary}
\begin{demo}
The implication from the left to the right follows from proposition
\ref{prop8.1}.

On the other hand, suppose that the set above is dense.

Let $b$ be an element in $B$ such that $(\pi(a) v_\rho b)_{\rho \in 
\cG_\vfi}$ is convergent in $E$. Because $b^* \rho(a^* a) b =$ \newline 
$\langle \pi(a) v_\rho b , \pi(a) v_\rho b \rangle$ for every $\rho \in 
\cG_\vfi$ (see notation \ref{not2.1}), this implies that the net $(b^* 
\rho(a^* a) b)_{\rho \in \cG_\vfi}$ is convergent in $B$.

Using this result, the density of the set above implies that $a$ is an 
element of $\tNfi$.
\end{demo}

\medskip

\begin{proposition}
Consider $a \in \tNfi$. Then $\la(a)$ is a closed densely defined $B$-linear 
map from within $B$ in $E$.
\end{proposition}
\begin{demo}
It follows from the previous result that $\la(a)$ is densely defined. We turn 
to the closedness.

Take a sequence $(b_n)_{n=1}^\infty$ in $D(\la(a))$  and elements
$b \in B$, $x \in E$ such that $(b_n)_{n=1}^\infty$ converges to $b$
and $(\la(a) b_n)_{n=1}^\infty$ converges to $x$.

\begin{list}{}{\setlength{\leftmargin}{.4 cm}}

\item Choose $\eta \in \cG_\vfi$. By notation \ref{not2.1},
we have for every $n \in \N$ that $$ \langle T_\eta^\frac{1}{2} \la(a) b_n , 
T_\eta^\frac{1}{2} \la(a) b_n \rangle = \langle \pi(a) v_\eta b_n , \pi(a) 
v_\eta b_n \rangle .$$
From this, we get easily that
$$  \langle T_\eta^\frac{1}{2} x , T_\eta^\frac{1}{2} x \rangle
= \langle \pi(a) v_\eta b , \pi(a) v_\eta b \rangle
= b^* \eta(a^* a) b .$$

\end{list}

Therefore, the net $(b^* \rho(a^* a) b)_{\rho \in \cG_\vfi}$ converges
to $\langle x , x \rangle$. By proposition \ref{prop8.1}, we see that $b$ 
belongs
to $D(\la(a))$.

\begin{list}{}{\setlength{\leftmargin}{.4 cm}}

\item Choose $\eta \in \cG_\vfi$.
It is clear that $(T_\eta^\frac{1}{2} \la(a) b_n)_{n=1}^\infty$ converges
to $T_\eta^\frac{1}{2} x$.

Because the mapping $D(\la(a)) \rightarrow E : c \mapsto T_\eta^\frac{1}{2} 
\la(a) c$ is continuous, we have also that
$(T_\eta^\frac{1}{2} \la(a) b_n)_{n=1}^\infty$ converges to 
$T_\eta^\frac{1}{2} \la(a) b$. Hence, $T_\eta^\frac{1}{2} x = 
T_\eta^\frac{1}{2} \la(a) b$.

\end{list}

The fact that $(T_\rho)_{\rho \in \cG_\vfi}$ converges strongly to 1,
gives us that $\la(a) b = x$.

\medskip

Consequently, we have proven that $\la(a)$ is closed.
\end{demo}

Another consequence of proposition \ref{prop8.1} is the following result.

\begin{result}
Consider $a \in \tNfi$, $\rho \in \cF_\vfi$ and let $S$ be an element in 
$\cL(E) \cap \pi(A)'$ such that $S^* S = T_\rho$. Define $v$ to be the 
element in $\cL(B,E)$ such that $S \la(c) = \pi(c) v $ for every $c \in 
\Nfi$. We have for every $x \in E$ that $S^* x$ belongs to $D(\la(a)^*)$ and 
$\la(a)^* ( S^* x) = v^* \pi(a^*) x$.
\end{result}
\begin{demo}
By proposition \ref{prop8.1}, we have for every $b \in D(\la(a))$ that
$$ \langle \la(a) b , S^* x \rangle = \langle S \la(a) b , x \rangle
= \langle \pi(a) v b , x \rangle = \langle b , v^* \pi(a^*)  x \rangle .$$
The result follows immediately.
\end{demo}

A special case hereof is the following result.

\begin{result}
Consider $a \in \tNfi$, $\rho \in \cF_\vfi$ and $x \in E$. Then 
$T_\rho^\frac{1}{2} x$ belongs to $D(\la(a)^*)$ and 
$\la(a)^*(T_\rho^\frac{1}{2} x) = v_\rho^* \pi(a^*) x$.
\end{result}

\medskip

So we get the following density result.

\begin{corollary}
Consider $a \in \tNfi$. Then $\la(a)^*$ is densely defined.
\end{corollary}

\medskip\medskip

Using proposition \ref{prop8.1} once more, we find:

\begin{result} \label{res8.1}
Consider $a_1,a_2 \in \tNfi$ and $b_1 \in D(\la(a_1))$, $b_2 \in 
D(\la(a_2))$. Then $\langle T_\rho \la(a_1) b_1 , \la(a_2) b_2 \rangle
= b_2^* \, \rho(a_2^* \, a_1) b_1$ for every $\rho \in \cF_\vfi$.
\end{result}
\begin{demo}
By corollary \ref{corol8.1}, we have that
\begin{eqnarray*}
\langle T_\rho \la(a_1) b_1 , \la(a_2) b_2 \rangle
& = & \langle T_\rho^\frac{1}{2} \la(a_1) b_1 , T_\rho^\frac{1}{2} \la(a_2) 
b_2  \rangle
= \langle \pi(a_1) v_\rho \, b_1 , \pi(a_2) v_\rho \, b_2 \rangle \\
& = & b_2^* \, v_\rho^* \, \pi(a_2^* \, a_1) v_\rho b_1
= b_2^* \, \rho(a_2^* \, a_1) b_1
\end{eqnarray*}
where we used notation \ref{not2.1} in the last equality.
\end{demo}

\medskip

This implies immediately the following convergence result.

\begin{result} \label{res8.2}
Consider $a_1,a_2 \in \tNfi$ and $b_1 \in D(\la(a_1))$, $b_2 \in 
D(\la(a_2))$. Then the net $(b_2^* \, \rho(a_2^* \, a_1) b_1)_{\rho \in 
\cG_\vfi}$ converges to $\langle \la(a_1) b_1 , \la(a_2) b_2 \rangle$.
\end{result}

\medskip

\begin{remark} \rm
By using result \ref{res8.1}, it is easy to see that the following holds.

Consider $a \in \tNfi$. We have for every $c \in \tNfi$, $b \in D(\la(c))$ 
and $\rho \in \cF_\vfi$ that $T_\rho \la(c) b$ belongs to
$D(\la(a)^*)$ and $\la(a)^*(T_\rho \la(c) b) =  \rho(a^* c) b$.
\end{remark}

\bigskip

\begin{remark} \rm
Consider $a \in \tNfi$. We see that $\la(a)$ is a densely defined
closed $B$-linear operator from $B$ into $E$ such that $\la(a)^*$ is densely 
defined. So, it behaves rather well.

The question remains open whether it behaves perfect, i.e. whether
$\la(a)$ is a regular operator in the sense of \cite{Lan} (which is true if 
$B$ is commutative). Therefore, we will give the following definition.
\end{remark}

\begin{definition}
We define the set $\hat{\cN}_\vfi = \{ a \in \tNfi \mid \la(a) \text{ is 
regular } \}$.
\end{definition}

\medskip

\begin{remark} \rm
By the results proven so far, we have that $$\hat{\cN}_\vfi
= \{ a \in \tNfi \mid  1+\la(a)^* \la(a) \text{ has dense range in } B \}.$$ 
It seems to be an interesting question to look for simpler conditions on 
elements of $\tNfi$ in order for them to belong to $\hat{\cN}_\vfi$.

\smallskip\smallskip

Of course, the following is true. Consider $a \in \tNfi$ such that there 
exist a positive element $\sde$ afilliated with $B$ such that $\sde \subseteq 
\la(a)^* \la(a)$. Then $a$ belongs to $\hat{N}_\vfi$ and
$\la(a)^* \la(a) = \sde$.
\end{remark}

\medskip

\begin{proposition}
We have the following properties:
\begin{enumerate}
\item Consider $a,b \in M(A)$ such that $a^* a = b^* b$.
Then $a$ belongs to $\tilde{\cN}_\vfi$ $\Leftrightarrow$ $b$ belongs to 
$\tilde{\cN}_\vfi$
\item Consider $a,b \in \tilde{\cN}_\vfi$ such that $a^* a = b^* b$.
Then $D(\la(a)) = D(\la(b))$ and $\langle \la(a)x , \la(a) y \rangle
= \langle \la(b)x , \la(b)y \rangle$ for every $x,y \in D(\la(a))$.
\item Consider $a,b \in \tilde{\cN}_\vfi$ such that $a^*a = b^* b$,
then $\la(a)^* \la(a) = \la(b)^* \la(b)$.
\item Consider $a,b \in M(A)$ such that $a^* a = b^* b$.
Then $a$ belongs to $\hat{\cN}_\vfi$ $\Leftrightarrow$ $b$ belongs to
$\hat{\cN}_\vfi$.
\end{enumerate}
\end{proposition}
\begin{demo}
\begin{enumerate}
\item This follows immediately from the definition of  $\tilde{\cN}_\vfi$.
\item Proposition \ref{prop8.1}.1 implies immediately that $D(\la(a)) = 
D(\la(b))$.

Choose $x,y \in D(\la(a))$.
Result \ref{res8.1} implies that
$$\langle T_\rho \la(a) x , \la(a) y \rangle = y^* \rho(a^* a) x = y^*
\rho(b^* b) x = \langle T_\rho \la(b) x , \la(b) x \rangle $$
for every $\rho \in \cG_\vfi$.

Because $(T_\rho)_{\rho \in \cG_\vfi}$ converges strongly to 1,
this implies that $\langle \la(a) x , \la(a) y \rangle = \langle \la(b) x , 
\la(b) y \rangle $.
\item Choose $x \in D(\la(a)^* \la(a))$. It is clear that $x$ belongs to
$D(\la(a))$. So by 2, we have that $x$ belongs to $D(\la(b))$.

By 2 we have also for every  $y \in D(\la(b))$ that $y$ belongs to
$D(\la(a))$ and
$$ \langle \la(b) x , \la(b) y \rangle = \langle \la(a) x , \la(a) y \rangle 
= \langle \la(a)^*(\la(a) x) , y \rangle . $$
This implies by definition that $\la(b)x$ belongs to $D(\la(b)^*)$ and
$\la(b)^*( \la(b) x ) = \la(a)^* (\la(a) x)$.

From this discussion, we infer that $\la(a)^* \la(a) \subseteq \la(b)^* 
\la(b)$. Similarly, we have that $\la(b)^* \la(b) \subseteq \la(a)^* \la(a)$.

\item This follows from 1 and the fact that we have for every $c \in 
\tilde{\cN}_\vfi$ that $\la(c)$ is regular if and only if $1+\la(c)^* \la(c)$ 
has dense range in $B$.
\end{enumerate}
\end{demo}

\medskip

\begin{definition}
We define the set $\hat{\cM}_\vfi^+ = \{ a \in M(A)^+ \mid a^\frac{1}{2}
\text{ belongs to } \hat{\cN}_\vfi \}$.
For every $a \in \hat{\cM}_\vfi^+$, we put $\vfi(a) = \la(a^\frac{1}{2})^* 
\la(a^\frac{1}{2})$, so $\vfi(a)$ is a positive element affiliated with $B$.
\end{definition}

It is clear that $\overline{\cM}_\vfi^+ \subseteq \hat{\cM}_\vfi^+$ and that 
for every $a \in \overline{\cM}_\vfi^+$ the notation $\vfi(a)$ is consistent 
with the notation introduced before.

\medskip

The previous proposition implies the following one.

\begin{proposition}
Consider $a \in M(A)$. Then
\begin{enumerate}
\item $a$ belongs to $\hat{\cN}_\vfi$ $\Leftrightarrow$
$a^* a$ belongs to $\hat{\cM}_\vfi^+$.
\item If $a$ belongs to $\hat{\cN}_\vfi$, then $\vfi(a^* a) = \la(a)^* 
\la(a)$.
\end{enumerate}
\end{proposition}

\medskip

We have also the following generalizations of proposition \ref{prop8.2}.

\begin{proposition}
Consider $a \in \tNfi$.
\begin{enumerate}
\item Let $b$ be an element in $B$. Then $b$ belongs to the domain of
$\la(a)$  \newline
$\Leftrightarrow$ The net $(b^* \rho_i(a^* a) b)_{i \in I}$ is convergent in 
$B$. \newline
$\Leftrightarrow$ The net $(\pi(a) v_i b)_{i \in I}$ is convergent in $B$
\newline
$\Leftrightarrow$ The net $(\pi(a) w_i b)_{i \in I}$ is convergent in $B$
\item Let $b$ be an element in $D(\la(a))$. Then
$(\pi(a) v_i b)_{i \in I}$ and $(\pi(a) w_i b)_{i \in I}$ converge both to 
$\la(a) b$.
\end{enumerate}
\end{proposition}
\begin{demo}
We prove everything with respect to the elements $v_i \, (i\!\in\!I)$. The 
case with the elements $w_i \, (i\!\in\!I)$ is completely similar.

\begin{itemize}
\item Suppose that $b$ belongs to $D(\la(a))$.
We know by corollary \ref{corol8.1} that $\pi(a) v_i b= T_i^\frac{1}{2} 
\la(a) b$ for every $i \in I$. Therefore, the net $(\pi(a) v_i b)_{i \in I}$ 
converges to
$\la(a) b$.

\item Suppose that $(\pi(a) v_i b)_{i \in I}$ converges.
We have for every $i \in I$ that $b^* \rho_i(a^* a) b = \langle \pi(a) v_i b 
, \pi(a) v_i b \rangle$.

This implies that $(b^* \rho_i(a^* a) b)_{i \in I}$ converges.

\item Suppose that $(b^* \rho_i(a^* a) b)_{i \in I}$ converges.
By proposition \ref{prop7.2}, we know that $(b^* \rho(a^* a) b)_{\rho \in 
\cG_\vfi}$ converges.
Therefore, $b$ belongs to $D(\la(a))$ (proposition \ref{prop8.1}).
\end{itemize}

Combining these results, the proposition follows.
\end{demo}

Similarly, we have that

\begin{proposition}
Consider $a \in M(A)$. Then $a$ belongs to $\tNfi$
\newline $\Leftrightarrow$
the set
$\{ b \in B \mid (b^* \rho_i(a^* a) b)_{i \in I} \text{ is convergent in } B 
\} $
is dense in $B$.
\newline $\Leftrightarrow$
the set
$\{ b \in B \mid (\pi(a) v_i b)_{i \in I} \text{ is convergent in } E \} $ is 
dense in $B$
\newline $\Leftrightarrow$
the set
$\{ b \in B \mid (\pi_i(a) w_i b)_{i \in I} \text{ is convergent in } E \} $ 
is dense in $B$
\end{proposition}

\section{The tensor product of regular \cst-valued weights.}
\label{art9}

In this section, we are going to define and prove some properties
about the tensor product of two regular \cst-valued weights. (cfr. 
\cite{Verd} for weights). This procedure can be easily adapted to the case of 
the tensor product of more regular \cst-valued weights.

\medskip

For the rest of this section, we will fix \cst-algebras $A_1$,$A_2$,
$B_1$,$B_2$ and regular \cst-valued weights $\vfi_1$ from $A_1$ into $M(B_1)$ 
and $\vfi_2$ from $A_2$ into $M(B_2)$.

\medskip

Let us also take KSGNS-constructions $(E_1,\la_1,\pi_1)$ for $\vfi_1$ and 
$(E_2,\la_2,\pi_2)$ for $\vfi_2$.

\medskip

We will also fix a truncating net $(u_i^1)_{i \in I_1}$ for $\vfi_1$
and a truncating net $(u_i^2)_{i \in I_2}$ for $\vfi_2$

\begin{list}{}{\setlength{\leftmargin}{.4 cm}}

\item Take $k \in \{1,2\}$ and $i \in I_k$.

We define $S_{i}^k$ to be the unique element in $\cL(E_k)$ such
that $S_{i}^k \la_k(a) = \la_k(a \, u_i^k)$ for every $a \in \cN_{\vfi_k}$. 
Furthermore, we put $T_i^k = (S_i^k)^* (S_i^k)$.

Moreover, $\rho_i^k$ will denote the strict completely positive mapping from 
$A_k$ into $M(B_k)$ such that
$\rho_i^k(a_2^* \, a_1) = \la_k(a_2)^* T_i^k \, \la_k(a_1)$ for every 
$a_1,a_2 \in \cN_{\vfi_k}$.

\end{list}

\medskip

We have introduced the necessary notations, so we can start with the 
construction of the tensor product.

\begin{itemize}

\item The set $\cN_{\vfi_1} \od \cN_{\vfi_2}$ is a dense subalgebra of $A_1 
\ot A_2$. Furthermore, $\la_1 \od \la_2$
is a linear mapping from $\cN_{\vfi_1} \od \cN_{\vfi_2}$ into $\cL(B_1 \ot 
B_2, E_1 \ot E_2)$ such that
the set $\langle \, (\la_1 \od \la_2)(c)d \mid c \in \cN_{\vfi_1} \od 
\cN_{\vfi_2} , d \in B_1 \ot B_2 \rangle$ is dense in $E_1 \ot E_2$.

\item Choose $i \in I_1 \times I_2$. Define $T_i = T_{i_1}^1 \ot T_{i_2}^2 
\in \cL(E_1 \ot E_2)$. Then  $0 \leq T_i \leq 1$.

Furthermore, we put $\rho_i = \rho_{i_1}^1 \ot \rho_{i_2}^2$ which is a 
strict completely positive mapping from $A_1 \ot A_2$ into $M(B_1 \ot B_2)$.

It is easy to check that $\langle T_ i (\la_1 \od \la_2)(c_1) d_1 , (\la_1 
\od \la_2)(c_2) d_2 \rangle = d_2^* \, \rho_i(c_2^* \, c_1) \, d_1$ for every 
$c_1,c_2 \in \cN_{\vfi_1} \od \cN_{\vfi_2}$ and $d_1,d_2 \in B_1 \od B_2$.
This implies that $\langle T_ i (\la_1 \od \la_2)(c_1) d_1 , (\la_1 \od 
\la_2)(c_2) d_2 \rangle = d_2^* \, \rho_i(c_2^* \, c_1) \, d_1$ for every 
$c_1,c_2 \in \cN_{\vfi_1} \od \cN_{\vfi_2}$ and $d_1,d_2 \in B_1 \ot B_2$.

\item It is also clear that $(T_i)_{i \in I_1 \times I_2}$ converges strongly 
to 1.

\end{itemize}

The previous discussion implies that the elements $A_1 \ot A_2$, $B_1 \ot 
B_2$, $E_1 \ot E_2$, $\cN_{\vfi_1} \od \cN_{\vfi_2}$, $\la_1 \od \la_2$ 
satisfy the conditions of the beginning of section \ref{art6}.
Therefore, $\la_1 \od \la_2$ is closable for the norm topology on $A_1 \ot 
A_2$ and the strong topology on $L(B_1 \ot B_2, E_1 \ot E_2)$  (lemma 
\ref{lem6.2}).
We define $\la_1 \ot \la_2$ to be this closure of $\la_1 \od \la_2$.
We denote the domain of $\la_1 \ot \la_2$ by $N$.

\medskip

Then $N$ is a dense subspace of $A_1 \ot A_2$ which is a left ideal of $M(A_1 
\ot A_2)$ (proposition \ref{prop6.2}).

The remark after lemma \ref{lem6.3} and proposition \ref{prop6.3} imply that 
$\la_1 \ot \la_2$ is a linear mapping
$N$ into \newline $\cL(B_1 \ot B_2, E_1 \ot E_2)$ which
is closed for the strict topology
on $A_1 \ot A_2$ and the strong topology on \newline
$L(B_1 \ot B_2, E_1 \ot E_2)$.

It is easy to check that $(\pi_1 \ot \pi_2)(x) (\la_1 \od \la_2)(c) =
(\la_1 \od \la_2)(x c)$ for every $x,c \in \cN_{\vfi_1} \od \cN_{\vfi_2}$. 
Using proposition \ref{prop6.2}, we see that
$(\pi_1 \ot \pi_2)(x) (\la_1 \ot \la_2)(c) =
(\la_1 \ot \la_2)(x c)$ for every $x \in M(A_1 \ot A_2)$ and $c \in N$.

\medskip

Now we want to go a little bit further in our construction procedure.

\begin{itemize}
\item Choose $i \in I_1 \times I_2$. Then we define $u_i = u_{i_1}^1 \ot 
u_{i_2}^2 \in M(A_1 \ot A_2)$, so we have clearly that $\|u_i\| \leq 1$.
Moreover, define $S_i = S_{i_1}^1 \ot S_{i_2}^2$, then $S_i \in \cL(E_1 \ot 
E_2)$ and $\|S_i\| \leq 1$. We have also that $T_i = S_i^* S_i$.

It is straightforward to check for every $c \in \cN_{\vfi_1} \od 
\cN_{\vfi_1}$ that $c u_i$ belongs to $\cN_{\vfi_1} \od \cN_{\vfi_2}$ and 
\newline $S_i \, (\la_1 \od \la_2)(c) = (\la_1 \od \la_2)(c u_i)$.
Hence, lemma \ref{lem6.1} implies that $N u_i \subseteq N$ and that
\newline $S_i \, (\la_1 \ot \la_2)(c) = (\la_1 \ot \la_2)(c u_i)$ for every 
$c \in N$

\item It is clear that $(S_i)_{i \in I_1 \times I_2}$ converges strongly to 1 
and that $(u_i)_{i \in I_1 \times I_2}$ converges strictly to 1.
\end{itemize}

We are now in a position to define the \cst-valued weight $\vfi_1 \ot \vfi_2$.

\begin{definition} \label{def9.1}
The elements $A_1 \ot A_2$, $B_1 \ot B_2$, $E_1 \ot E_2$ , $N$
and $\la_1 \ot \la_2$ satisfy the conditions of the beginning of section 
\ref{art7}.
Therefore, we can use definition \ref{def7.1} to define the tensor product 
$\vfi_1 \ot \vfi_2$, using these ingredients.
\end{definition}

By corollary \ref{corol7.1}, we have the following proposition.

\begin{proposition}
We have that $\vfi_1 \ot \vfi_2$ is a regular \cst-valued weight from
$A_1 \ot A_2$ into $M(B_1 \ot B_2)$.
\end{proposition}

Corollary \ref{corol7.2} implies the following proposition.

\begin{proposition}
We have that $(E_1 \ot E_2, \la_1 \ot \la_2, \pi_1 \ot \pi_2)$ is a KSGNS-
construction for $\vfi_1 \ot \vfi_2$.
\end{proposition}

Using this proposition, it is not difficult to check that our definition of 
$\vfi_1 \ot \vfi_2$ is independent of the used construction procedure.

This last proposition determines $\vfi_1 \ot \vfi_2$ completely and 
definition \ref{def9.1} becomes in a certain sense irrelevant. In fact, we 
want to forget the discussion before the definition.  To keep things clear, 
we gather the useful results of this discussion.

\begin{proposition}
The mapping $\la_1 \ot \la_2$ is a linear mapping from $\cN_{\vfi_1 \ot 
\vfi_2}$ into $\cL(B_1 \ot B_2, E_1 \ot E_2)$ which is closed for the strict 
topology on $A_1 \ot A_2$ and the strong topology on $L(B_1  \ot B_2, E_1 \ot 
E_2)$.

Moreover, $\cN_{\vfi_1} \od \cN_{\vfi_2}$ is a norm-strong core for $\la_1 
\ot \la_2$ and $(\la_1 \ot \la_2)(a_1 \ot a_2) = \la_1(a_1) \ot \la_2(a_2)$ 
for every $a_1 \in \cN_{\vfi_1}$ and $a_2 \in \cN_{\vfi_2}$.
\end{proposition}

For every $i \in I_1 \times I_2$, we have defined the element $u_i = 
u_{i_1}^1 \ot u_{i_2}^2 \in M(A_1 \ot A_2)$.
Furthermore, we have defined the elements $S_i = S_{i_1}^1 \ot S_{i_2}^2$
and $T_i = T_{i_1}^1 \ot T_{i_2}^2$ in $\cL(E_1 \ot E_2)$.
So $T_i = S_i^* S_i$.

We also defined the element $\rho_i = \rho_{i_1}^1 \ot \rho_{i_2}^2$, which 
is a strict completely positive mapping from $A_1 \ot A_2$ into $M(B_1 \ot 
B_2)$.

Concerning this elements, we have the following properties.

\begin{proposition}
The net $(u_i)_{i \in I_1 \times I_2}$ is  truncating for $\vfi_1 \ot \vfi_2$.
\end{proposition}

From this, we get immediately that $\vfi_1 \ot \vfi_2$ is strongly regular if 
$\vfi_1$ and $\vfi_2$ are strongly regular.

\begin{proposition} \label{prop9.1}
Consider $i \in I_1 \times I_2$. Then we have the following properties:
\begin{itemize}
\item We have for every $c \in \cN_{\vfi_1 \ot \vfi_2}$ that
$S_i \, (\la_1 \ot \la_2)(c) = (\la_1 \ot \la_2)(c u_i)$.
\item For every $c_1,c_2 \in \cN_{\vfi_1 \ot \vfi_2}$, we have the equality 
$\rho_i(c_2^* \, c_1) = (\la_1 \ot \la_2)(c_2)^* T_i (\la_1 \ot \la_2)(c_1)$.
\item We have the inclusion $(\cN_{\vfi_1} \od \cN_{\vfi_2}) \, u_i
\subseteq \cN_{\vfi_1} \od \cN_{\vfi_2}$
\end{itemize}
\end{proposition}

Take $k \in \{1,2\}$ and $i \in I_k$. Then $v_i^k$ and $w_i^k$ will denote 
the unique elements in $\cL(B_k,E_k)$ such that
$(T_i^k)^\frac{1}{2} \la_k(a) = \pi_k(a) v_i^k$ and
$S_i^k \la_k(a) = \pi_k(a) w_i^k$ for every $a \in \cN_{\vfi_k}$.

\medskip

For every $i \in I_1 \times I_2$, we define the elements $v_i$ and $w_i$ in 
$\cL(B_1 \ot B_2, E_1 \ot E_2)$ such that $v_i = v_{i_1}^1 \ot v_{i_2}^2$ and 
$w_i =  w_{i_1}^1 \ot w_{i_2}^2$.
It will not be a great surprise that these elements will be the corresponding 
objects for $\vfi_1 \ot \vfi_2$. More precisely:

\begin{result}
Consider  $i \in I_1 \times I_2$. Then we have that
$T_i^\frac{1}{2} \la(c) = (\pi_1 \ot \pi_2)(c) v_i$
and $S_i \la(c) = (\pi_1 \ot \pi_2)(c) w_i$ for every
$c \in \cN_{\vfi_1 \ot \vfi_2}$.
\end{result}

\begin{result} \label{res9.1}
Consider $\om_1 \in \cF_{\vfi_1}$ and $\om_2 \in \cF_{\vfi_2}$. Then $\om_1 
\ot \om_2$ belongs to $\cF_{\vfi_1 \ot \vfi_2}$ and
$T_{\om_1 \ot \om_2} = T_{\om_1} \ot T_{\om_2}$.
\end{result}

Both results follow by checking equalities for elements in $\cN_{\vfi_1} \od 
\cN_{\vfi_2}$ and then using the fact that $\cN_{\vfi_1} \od \cN_{\vfi_2}$ is 
a norm-strong core for $\la_1 \ot \la_2$.
(like in the proof of lemma \ref{lem6.1}).

\begin{corollary}
Consider $\om_1 \in \cG_{\vfi_1}$ and $\om_2 \in \cG_{\vfi_2}$. Then $\om_1 
\ot \om_2$ belongs to $\cG_{\vfi_1 \ot \vfi_2}$
\end{corollary}

\bigskip

The following lemma will be very useful to us.

\begin{lemma} \label{lem9.1}
Consider $c \in M(A_1 \ot A_2)^+$, $d \in B_1 \ot B_2$ and $x \in B_1 \ot 
B_2$ such that the net \newline
$(\,d^*(\om_1 \ot \om_2)(c)d\,)_{ \om \in \cG_{\vfi_1} \times \cG_{\vfi_2}}$ 
converges to $x$.
Then $(\,d^* \om(c) d\,)_{\om \in \cG_{\vfi_1 \ot \vfi_2}}$ converges to $x$.
\end{lemma}
\begin{demo}
Because $(\,d^*(\om_1 \ot \om_2)(c)d\,)_{\om \in \cG_{\vfi_1} \times 
\cG_{\vfi_2}}$ is an increasing net in $B$ which converges to $x$, we get that
$d^* (\om_1 \ot \om_2)(c) d \leq x$ for every $\om_1 \in \cG_{\vfi_1}$ and 
$\om_2 \in \cG_{\vfi_2}$.
From this, it is easy to conclude that $d^* (\om_1 \ot \om_2)(c) d \leq x$ 
for every $\om_1 \in \cF_{\vfi_1}$ and $\om_2 \in \cF_{\vfi_2}$.

In particular, we see that $d^* \rho_i(c) d \leq x$
for every $i \in I_1 \times I_2$.
From proposition \ref{prop7.3}, we infer that $d^* \rho(c) d \leq x$
for every $\rho \in \cF_{\vfi_1 \ot \vfi_2}$.

Because we also have that
$(\,d^*(\om_1 \ot \om_2)(c)d\,)_{ \om \in \cG_{\vfi_1} \times \cG_{\vfi_2}}$ 
converges to $x$, lemma \ref{lem2.1} implies that
\newline $(d^* \om(c) d)_{\rho \in \cG_{\vfi_1 \ot \vfi_2}}$ converges to $x$.
\end{demo}

We will find a first application of this lemma in the next theorem, which 
gives a nice characterization of $\vfi_1 \ot \vfi_2$.

\begin{theorem}  \label{theo9.1}
We have the following properties.
\begin{enumerate}
\item Consider $c \in M(A_1 \ot A_2)^+$. Then $c$ belongs to 
$\overline{\cM}_{\vfi_1 \ot \vfi_2}^+$ $\Leftrightarrow$ The net 
$(\,d^*(\om_1 \ot \om_2)(c)d\,)_{ \om \in \cG_{\vfi_1} \times \cG_{\vfi_2}}$ 
is convergent in $B_1 \ot B_2$ for every $d \in B_1 \ot B_2$.
\item Let $c$ be an element in $\overline{\cM}_{\vfi_1 \ot \vfi_2}$.
Then the net $(\,(\om_1 \ot \om_2)(c)\,)_{\om \in \cG_{\vfi_1} \times 
\cG_{\vfi_2}}$ converges strictly to  \newline $(\vfi_1 \ot \vfi_2)(c)$.
\end{enumerate}
\end{theorem}
\begin{demo}
\begin{itemize}
\item Suppose that the net $(\,d^*(\om_1 \ot \om_2)(c)d\,)_{ \om \in 
\cG_{\vfi_1} \times \cG_{\vfi_2}}$ is convergent in $B_1 \ot B_2$ for every 
$d \in B_1 \ot B_2$.

By the previous lemma, we see that $(d^* \om(c) d)_{\om \in \cG_{\vfi_1 \ot 
\vfi_2}}$ is convergent for every $d \in B_1 \ot B_2$.
By definition, this implies that $c$ belongs to $\overline{\cM}_{\vfi_1 \ot 
\vfi_2}^+$.

\item Choose $c_1,c_2 \in \cN_{\vfi_1 \ot \vfi_2}$.
We know by result \ref{res9.1} that
$$ (\om_1 \ot \om_2)(c_2^* \, c_1)= (\la_1 \ot \la_2)(c_2)^* (T_{\om_1} \ot 
T_{\om_2}) (\la_1 \ot \la_2)(c_1) $$
for every $\om \in \cG_{\vfi_1} \times \cG_{\vfi_2}$.

Because the net $(T_{\om_1} \ot T_{\om_2})_{\om \in \cG_{\vfi_1} \times 
\cG_{\vfi_2}}$ converges strongly to $1$, this implies that the net \newline
$(\,(\om_1 \ot \om_2)(c_2^* \, c_1)\,)_{\om \in \cG_{\vfi_1} \times 
\cG_{\vfi_2}}$ converges strictly to $(\la_1 \ot \la_2)(c_2)^*
(\la_1 \ot \la_2)(c_1)$ which is equal to
$(\vfi_1 \ot \vfi_2)(c_2^* \, c_1)$.
\end{itemize}
The proposition follows easily from these two items.
\end{demo}

\begin{corollary}
Consider $c \in (A_1 \ot A_2)^+$. Then $c$ belongs to $\cM_{\vfi_1 \ot 
\vfi_2}^+$ $\Leftrightarrow$ The net \newline $(\,d^*(\om_1 \ot 
\om_2)(c)d\,)_{ \om \in \cG_{\vfi_1} \times \cG_{\vfi_2}}$ is convergent in 
$B_1 \ot B_2$ for every $d \in B_1 \ot B_2$.
\end{corollary}

\medskip

\begin{corollary} \label{corol9.1}
Consider $a_1 \in \overline{\cM}_{\vfi_1}$ and $a_2 \in 
\overline{\cM}_{\vfi_2}$.
Then $a_1 \ot a_2$ belongs to $\overline{\cM}_{\vfi_1 \ot \vfi_2}$
and \newline $(\vfi_1 \ot \vfi_2)(a_1 \ot a_2) = \vfi(a_1) \ot \vfi(a_2)$.
\end{corollary}

In fact, this corollary follows easily from theorem \ref{theo9.1} for $a_1 
\in \overline{\cM}_{\vfi_1}^+$ and $a_2 \in \overline{\cM}_{\vfi_2}^+$.
Because any element of $\overline{\cM}_{\vfi_1}$ can be written as a linear 
combination of elements in $\overline{\cM}_{\vfi_1}^+$ (and similarly for 
$\overline{\cM}_{\vfi_2}$), the corollary follows.

\begin{corollary}
Consider $a_1 \in \cM_{\vfi_1}$ and $a_2 \in \cM_{\vfi_2}$. Then
$a_1 \ot a_2$ belongs to $\cM_{\vfi_1 \ot \vfi_2}$.
\end{corollary}

\bigskip

A result which is very much related to this one, is the following one.

\begin{result}
Consider $a_1 \in \overline{\cN}_{\vfi_1}$ and $a_2 \in 
\overline{\cN}_{\vfi_2}$. Then
$a_1 \ot a_2$ belongs to $\overline{\cN}_{\vfi_1 \ot \vfi_2}$ and
$(\la_1 \ot \la_2)(a_1 \ot a_2) = \la(a_1) \ot \la(a_2)$.
\end{result}
\begin{demo}
Because $a_k^* a_k$ belongs to $\overline{\cM}_{\vfi_k}^+$ ($k=1,2$),
corollary \ref{corol9.1} implies that $a_1^* a_1 \ot a_2^* a_2$ belongs to 
$\overline{\cM}_{\vfi_1 \ot \vfi_2}^+$. Therefore, $a_1 \ot a_2$ belongs to 
$\overline{\cN}_{\vfi_1 \ot \vfi_2}$.

Choose $e_1 \in A_1$ and $e_2 \in A_2$. Then $e_k \, a_k$ belongs to 
$\cN_{\vfi_k}$ and $\la_k(e_k \, a_k) = \pi(e_k) \la_k(a_k)$ ($k=1,2$).
By definition, we get that $e_1 a_1 \ot e_2 a_2$ belongs to $\cN_{\vfi_1 \ot 
\vfi_2}$ and
$$(\la_1 \ot \la_2)(e_1 a_1 \ot e_2 a_2)
= \la_1(e_1 a_1) \ot \la(e_2 a_2)
= (\pi_1(e_1) \ot \pi_2(e_2))\, (\la_1(a_1) \ot \la(a_2)) .$$
This last equality implies that
$$(\pi_1(e_1) \ot \pi_2(e_2))(\la_1 \ot \la_2)(a_1 \ot a_2)
= (\pi_1(e_1) \ot \pi_2(e_2))\, (\la_1(a_1) \ot \la(a_2)) .$$
Using the non-degeneracy of $\pi_1$ and $\pi_2$, we get that
$(\la_1 \ot \la_2)(a_1 \ot a_2) = \la_1(a_1) \ot \la_2(a_2)$.
\end{demo}

\bigskip

Using proposition \ref{prop8.3} (and the subsequent results) and proposition 
\ref{prop9.1}.3 , we get the following results.

\begin{proposition}
Consider $x \in \cN_{\vfi_1 \ot \vfi_2}$ such that $(\la_1 \ot \la_2)(x) \neq 
0$. Then there exists a  net $(x_j)_{j \in J}$ in $\cN_{\vfi_1} \od 
\cN_{\vfi_2}$ such that \begin{enumerate}
\item $\|x_j\| \leq \|x\|$ and $\|(\la_1 \ot \la_2)(x_j)\| \leq \|(\la_1 \ot 
\la_2)(x)\|$ for every $j \in J$.
\item $(x_j)_{j \in J}$ converges to $x$ and $(\,(\la_1 \ot \la_2)(x_j)\,)_{j 
\in J}$ converges strongly$^*$ to $(\la_1 \ot \la_2)(x)$.
\end{enumerate}
\end{proposition}

\begin{proposition}
Consider $x \in \cN_{\vfi_1 \ot \vfi_2}$ and let $M$ be a positive number 
such that $\|(\la_1 \ot \la_2)(x)\| < M $. Then there exists a  net $(x_j)_{j 
\in J}$ in $\cN_{\vfi_1} \od \cN_{\vfi_2}$ such that \begin{enumerate}
\item $\|x_j\| \leq \|x\|$ and $\|(\la_1 \ot \la_2)(x_j)\| \leq M$ for every 
$j \in J$.
\item $(x_j)_{j \in J}$ converges to $x$ and $(\,(\la_1 \ot \la_2)(x_j)\,)_{j 
\in J}$ converges strongly$^*$ to $(\la_1 \ot \la_2)(x)$.
\end{enumerate}
\end{proposition}

\begin{proposition}
Consider $x \in \overline{\cN}_{\vfi_1 \ot \vfi_2}$ such that $(\la_1 \ot 
\la_2)(x) \neq 0$. Then there exists a  net $(x_j)_{j \in J}$ in 
$\cN_{\vfi_1} \od \cN_{\vfi_2}$ such that \begin{enumerate}
\item $\|x_j\| \leq \|x\|$ and $\|(\la_1 \ot \la_2)(x_j)\| \leq \|(\la_1 \ot 
\la_2)(x)\|$ for every $j \in J$.
\item $(x_j)_{j \in J}$ converges strictly to $x$ and $(\,(\la_1 \ot 
\la_2)(x_j)\,)_{j \in J}$ converges strongly$^*$ to $(\la_1 \ot \la_2)(x)$.
\end{enumerate}
\end{proposition}

\begin{proposition}
Consider $x \in \overline{\cN}_{\vfi_1 \ot \vfi_2}$ and let $M$ be a positive 
number such that $\|(\la_1 \ot \la_2)(x)\| < M $. Then there exists a  net 
$(x_j)_{j \in J}$ in $\cN_{\vfi_1} \od \cN_{\vfi_2}$ such that 
\begin{enumerate}
\item $\|x_j\| \leq \|x\|$ and $\|(\la_1 \ot \la_2)(x_j)\| \leq M$ for every 
$j \in J$.
\item $(x_j)_{j \in J}$ converges strictly to $x$ and $(\,(\la_1 \ot 
\la_2)(x_j)\,)_{j \in J}$ converges strongly$^*$ to $(\la_1 \ot \la_2)(x)$.
\end{enumerate}
\end{proposition}

\bigskip\bigskip

In the last part of this section, we want to prove some results
about elements in $\tilde{\cN}_{\vfi_1 \ot \vfi_2}$.

\begin{proposition} \label{prop9.2}
We have the following properties:
\begin{enumerate}
\item Consider $c \in \tilde{\cN}_{\vfi_1 \ot \vfi_2}$ and $d \in B_1 \ot 
B_2$ Then $d$ belongs to $D((\la_1 \ot \la_2)(c))$ $\Leftrightarrow$
The net \newline  $(\,d^* (\om_1 \ot \om_2)(c^* c) d \,)_{\om \in 
\cG_{\vfi_1} \times \cG_{\vfi_2}}$ is convergent in $B_1 \ot B_2$.
\item Consider $c_1,c_2 \in \tilde{\cN}_{\vfi_1 \ot \vfi_2}$ and $d_1
\in D((\la_1 \ot \la_2)(c_1))$, $d_2 \in D((\la_1 \ot \la_2)(c_2))$. Then we 
have that the net
$(\,d_2^* \, (\om_1 \ot \om_2)(c_2^* \, c_1) \, d_1 \,)_{\om \in \cG_{\vfi_1} 
\times \cG_{\vfi_2}}$ converges to
$\langle (\la_1 \ot \la_2)(c_1) d_1 , (\la_1 \ot \la_2)(c_2) d_2 \rangle$.
\end{enumerate}
\end{proposition}
\begin{demo}
\begin{itemize}
\item Choose $c \in \tilde{\cN}_{\vfi_1 \ot \vfi_2}$. Suppose that $d$ is an 
element in $B_1 \ot B_2$ such that the net \newline $(\,d^* (\om_1 \ot 
\om_2)(c^* c) d \,)_{\om \in \cG_{\vfi_1} \times \cG_{\vfi_2}}$ is convergent 
in $B_1 \ot B_2$.
By lemma \ref{lem9.1}, we know that the net \newline $(\, d^* \om(c^* c) d 
\,)_{\om \in \cG_{\vfi_1 \ot \vfi_2}}$ is convergent in $B_1 \ot B_2$.
Using proposition \ref{prop8.1}, we see that $d$ belongs to $D((\la_1 \ot 
\la_2)(c))$.
\item Choose $c_1,c_2 \in \tilde{\cN}_{\vfi_1 \ot \vfi_2}$ and $d_1
\in D((\la_1 \ot \la_2)(c_1))$, $d_2 \in D((\la_1 \ot \la_2)(c_2))$.
By  result \ref{res8.1} and result \ref{res9.1}, we have for every $\om \in 
\cG_{\vfi_1} \times \cG_{\vfi_2}$  that
$$d_2^* \, (\om_1 \ot \om_2)(c_2^* \, c_1) \, d_1 = \langle (T_{\om_1} \ot 
T_{\om_2}) (\la_1 \ot \la_2)(c_1) d_1 , (\la_1 \ot \la_2)(c_2) d_2 \rangle .$$
This implies immediately that the net $(\,d_2^* \, (\om_1 \ot \om_2)(c_2^* \, 
c_1) \, d_1 \,)_{\om \in \cG_{\vfi_1} \times \cG_{\vfi_2}}$ converges to 
\newline
$\langle (\la_1 \ot \la_2)(c_1) d_1 , (\la_1 \ot \la_2)(c_2) d_2 \rangle$.
\end{itemize}
Combining these two results, the proposition follows.
\end{demo}

\begin{proposition}
Consider $c \in M(A_1 \ot A_2)$. Then $c$ belongs to $\tilde{\cN}_{\vfi_1 \ot 
\vfi_2}$ $\Leftrightarrow$ The set
$$\{ d \in B_1 \ot B_2 \mid \text{ the net } (\, d^* (\om_1 \ot \om_2)(c^* c) 
d \,)_{\om \in \cG_{\vfi_1} \times \cG_{\vfi_2}} \text{ is convergent in } 
B_1 \ot B_2 \} $$
is dense in $B_1 \ot B_2$.
\end{proposition}

One implication of this proposition follows from the previous proposition. 
The other one follows from lemma \ref{lem9.1}.

\begin{proposition}
Consider $a_1 \in \tilde{\cN}_{\vfi_1}$ and $a_2 \in \tilde{\cN}_{\vfi_2}$. 
Then $a_1 \ot a_2$ belongs to $\tilde{\cN}_{\vfi_1 \ot \vfi_2}$, \ \ 
$\la_1(a_1) \od \la_2(a_2)$ is closable and its closure is a restriction of 
$(\la_1 \ot \la_2)(a_1 \ot a_2)$.
\end{proposition}
\begin{demo}
\begin{itemize}
\item Choose $d_1 \in D(\la_1(a_1))$, $d_2 \in D(\la_2(a_2))$.

We have for every $\om_1 \in \cG_{\vfi_1}$ and $\om_2 \in \cG_{\vfi_2}$
that
$$ (d_1 \ot d_2)^* (\om_1 \ot \om_2)((a_1 \ot a_2)^* (a_1 \ot a_2)) (d_1 \ot 
d_2)
= d_1^* \, \om_1(a_1^* \, a_1) \, d_1 \ot d_2^* \, \om_2(a_2^* \, a_2) \, d_2 
.$$
Using result \ref{res8.2}, this implies that the net
$$\bigl( \, (d_1 \ot d_2)^* (\om_1 \ot \om_2)((a_1 \ot a_2)^* (a_1 \ot a_2)) 
(d_1 \ot d_2) \, \bigr)_{\om \in \cG_{\vfi_1} \ot \cG_{\vfi_2}}$$
converges to $\langle \la_1(a_1) d_1 , \la_1(a_1) d_1 \rangle \ot
\langle \la_2(a_2) d_2 , \la_2(a_2) d_2 \rangle$.

\medskip

Therefore, we conclude from the two previous propositions that $a_1 \ot a_2$ 
belongs to $\tilde{\cN}_{\vfi_1 \ot \vfi_2}$ and that $D(\la_1(a_1)) \od 
D(\la_2(a_2))$ is a subset of  $D((\la_1 \ot \la_2)(a_1 \ot a_2))$.

\item Choose $d_1 \in D(\la_1(a_1))$ and $d_2 \in D(\la_2(a_2))$. From the 
first part we already know that $d_1 \ot d_2$ belongs
to $D((\la_1 \ot \la_2)(a_1 \ot a_2))$.

Choose $b_1 \in \cN_{\vfi_1}$, $b_2 \in \cN_{\vfi_2}$, $c_1 \in B_1$ and $c_2 
\in B_2$. Then we have for every $\om_1 \in \cG_{\vfi_1}$ and $\om_2 \in 
\cG_{\vfi_2}$ that
$$ (c_1 \ot c_2)^* (\om_1 \ot \om_2)((b_1 \ot b_2)^* (a_1 \ot a_2)) (d_1 \ot 
d_2)
= c_1^* \, \om_1(b_1^* \, a_1) \, d_1 \ot c_2^* \, \om_2(b_2^* \, a_2) \, d_2 
.$$

Using result \ref{res8.2}, this implies that the net
$$\bigl( \, (c_1 \ot c_2)^* (\om_1 \ot \om_2)((b_1 \ot b_2)^* (a_1 \ot a_2)) 
(d_1 \ot d_2) \, \bigr)_{\om \in \cG_{\vfi_1} \ot \cG_{\vfi_2}}$$
converges to $\langle \la_1(a_1) d_1 , \la_1(b_1) c_1 \rangle \ot
\langle \la_2(a_2) d_2 , \la_2(b_2) c_2 \rangle$ which is equal to
\newline $\langle \la_1(a_1) d_1 \ot \la_2(a_2) d_2 , \la_1(b_1) c_1 \ot 
\la_2(b_2) c_2 \rangle$. \ \ \ \ \ \ (a)

On the other hand, proposition  \ref{prop9.2} guarantees that the net 
$$\bigl( \, (c_1 \ot c_2)^* (\om_1 \ot \om_2)((b_1 \ot b_2)^* (a_1 \ot a_2)) 
(d_1 \ot d_2) \, \bigr)_{\om \in \cG_{\vfi_1} \times \cG_{\vfi_2}}$$
converges to $\langle (\la_1 \ot \la_2)(a_1 \ot a_2)(d_1 \ot d_2) ,
(\la_1 \ot \la_2)(b_1 \ot b_2)(c_1 \ot c_2) \rangle$ which by definition 
equals $\langle (\la_1 \ot \la_2)(a_1 \ot a_2)(d_1 \ot d_2) , \la_1(b_1)c_1 
\ot  \la_2(b_2)c_2 \rangle$. \ \ \ \ \ \  (b)

Combining (a) and (b), we get that
$$ \langle \la_1(a_1) d_1 \ot \la_2(a_2) d_2 , \la_1(b_1) c_1 \ot \la_2(b_2) 
c_2 \rangle = \langle (\la_1 \ot \la_2)(a_1 \ot a_2)(d_1 \ot d_2) , 
\la_1(b_1)c_1  \ot  \la_2(b_2)c_2 \rangle .$$

From this, we infer that
$(\la_1 \ot \la_2)(a_1 \ot a_2)(d_1 \ot d_2) = \la_1(a_1) d_1 \ot \la_2(a_2) 
d_2$.

\medskip

Hence, we have proven that $\la_1(a_1) \od \la_2(a_2) \subseteq (\la_1 \ot 
\la_2)(a_1 \ot a_2)$. Because $(\la_1 \ot \la_2)(a_1 \ot a_2)$ is closed, the 
lemma follows.
\end{itemize}
\end{demo}

\medskip

\begin{remark} \rm Consider  Hilbert \cst-modules $F_1, F_2$ over a \cst-
algebra $C$. Let $t$ be a regular operator from within $F_1$ into $F_2$. A 
truncating sequence for $t$ is by definition a sequence $(e_n)_{n=1}^\infty$
in $\cL(F_1)$ such that
\begin{enumerate}
\item We have for every $n \in \N$ that $\|e_n\| \leq 1$.
\item The net $(e_n)_{n=1}^\infty$ converges strictly to 1.
\item We have for every $n \in \N$ that $e_n \, |t| \subseteq |t| \, e_n$
and $|t| \, e_n$ belongs to $\cL(F_1)$.
\end{enumerate}
We should mention that, by using the functional calculus for $|t|$, we get 
the existence of a truncating sequence for $t$.

\medskip

Let us take a truncating sequence $(e_n)_{n=1}^\infty$ for $t$.

\medskip

We know that $D(t^* t)$ is a core for $t$ and $|t|$. It is easy to check that 
$\|t(x)\| = \|\,|t|(x)\|$ for every $x \in D(t^* t)$. Using these two facts 
and the closedness of $t$ and $|t|$, it is not difficult to check that $D(t) 
= D(|t|)$ and that $\|t(x)\| = \| \, |t|(x) \|$ for every $x \in D(t)$. We 
will use this result a few times during this remark.

\medskip

For the moment, fix $n \in \N$. Because $D(t) = D(|t|)$, we have for every $x 
\in A$ that $e_n x$ belongs to $D(t)$.

Moreover, we have for every $x \in A$ that $\| t(e_n x) \|  =
\|\,|t|(e_n x) \| $.

This equality implies that the mapping $t e_n$ is bounded and $\|t e_n\| = \| 
\, |t| e_n \|$. It is not difficult to check that $e_n^* t^* \subseteq (t 
e_n)^*$, which implies that $(t e_n)^*$ is densely defined.

These two facts imply that $t e_n$ belongs to $\cL(F_1,F_2)$. We also get 
that $e_n^* t^*$ is bounded and $(t e_n)^* = \overline{e_n^* t^*}$.

\medskip

Choose $x \in A$. Then
\begin{itemize}
\item If $x$ belongs to $D(t)$, then $(t(e_n x))_{n=1}^\infty$ converges to 
$t(x)$.

\begin{list}{}{\setlength{\leftmargin}{.4 cm}}

\item Because $x$ belongs to $D(t)$, $x$ also belongs to $D(|t|)$.

Moreover, we have for every $n \in \N$ that
$$\| t(e_n x) - t(x) \| = \| \, |t|(e_n x - x)  \|
= \| e_n \, |t|(x) - |t|(x) \| \ ,$$
which implies that $(t(e_n x))_{n=1}^\infty$ converges to $t(x)$.

\end{list}

\item We have that $x$ belongs to $D(t)$ $\Leftrightarrow$ The net
$(t(e_n x))_{n=1}^\infty$ is convergent.

\begin{list}{}{\setlength{\leftmargin}{.4 cm}}

\item One implication follows from the first part, the other implication 
follows from the closedness of $t$ and the fact that $(e_n x)_{n=1}^\infty$ 
converges to $x$.

\end{list}

\end{itemize}

\end{remark}

The above results will enable us to prove the following theorem.

\begin{theorem}
Consider $a_1 \in \hat{\cN}_{\vfi_1}$ and $a_2 \in  \hat{\cN}_{\vfi_2}$.
Then $a_1 \ot a_2$ belongs to $\hat{\cN}_{\vfi_1 \ot \vfi_2}$ and
\newline $(\la_1 \ot \la_2)(a_1 \ot a_2) = \la_1(a_1) \ot \la_2(a_2)$.
\end{theorem}
\begin{demo}
By the previous lemma, we know that $a_1 \ot a_2$ belongs to 
$\tilde{\cN}_{\vfi_1 \ot \vfi_2}$ and
$\la_1(a_1) \ot \la_2(a_2) \subseteq (\la_1 \ot \la_2)(a_1 \ot a_2)$.

Take truncating sequences $(e_n)_{n=1}^\infty, (f_n)_{n=1}^\infty$ for 
$\la_1(a_1)$, $\la_2(a_2)$ respectively. It is not too difficult to check 
that $(e_n \ot f_n)_{n=1}^\infty$ is a truncating sequence for $\la_1(a_1) 
\ot \la_2(a_2)$.

Choose $m \in \N$. We will first prove that $(\la_1 \ot \la_2)(a_1 \ot 
a_2)(e_m \ot f_m) = (\la_1(a_1) \ot \la_2(a_2))(e_m \ot f_m)$.

\begin{list}{}{\setlength{\leftmargin}{.4 cm}}

\item Choose $y \in B_1 \ot B_2$. Then there exists a sequence 
$(y_j)_{j=1}^\infty$ in $B_1 \od B_2$ such that $(y_j)_{j=1}^\infty$
converges to $y$.
We clearly have that $(\,(e_m \ot f_m) \, y_j)_{j=1}^\infty$ converges to
$(e_m \ot f_m)\, y$.

By the previous lemma, we have for every $j \in \N$ that
$(e_m \ot f_m)\, y_j$ belongs to $D((\la_1 \ot \la_2)(a_1 \ot a_2))$ and
$$ (\la_1 \ot \la_2)(a_1 \ot a_2)\bigl(\,(e_m \ot f_m) \, y_j\bigr)
= [(\la_1(a_1) \ot \la_2(a_2))(e_m \ot f_m)] \, y_j .$$
Because $[(\la_1(a_1) \ot \la_2(a_2))(e_m \ot f_m)]$ is bounded, this implies 
that  the net
$$\bigl(\, (\la_1 \ot \la_2)(a_1 \ot a_2)\bigl(\,(e_m \ot f_m)\, y_j\bigr) \,
\bigr)_{j=1}^\infty$$ converges to
$[(\la_1(a_1) \ot \la_2(a_2))(e_m \ot f_m)] \, y$.

Therefore, the closedness of $(\la_1 \ot \la_2)(a_1 \ot a_2)$ implies that 
$(e_m \ot f_m)\,y$ belongs to $D((\la_1 \ot \la_2)(a_1 \ot a_2))$
and $(\la_1 \ot \la_2)(a_1 \ot a_2)\bigl( \,(e_m \ot f_m)\, y \bigr)
= (\la_1(a_1) \ot \la_2(a_2))\bigl( \,(e_m \ot f_m)\,  y \bigr)$.

\end{list}

Choose $x \in D(\la_1(a_1) \ot D(\la_2(a_2))$.
By the previous discussion, we have for every $n \in \N$ that $(e_n \ot 
f_n)\, x$ belongs to $D((\la_1 \ot \la_2)(a_1 \ot a_2))$ and
$$ (\la_1 \ot \la_2)(a_1 \ot a_2)\bigl(\,(e_n \ot f_n)\,x \bigr)
= (\la_1(a_1) \ot \la_2(a_2))\bigl( \, (e_n \ot f_n) \, x \bigr) .$$
Because $(e_n \ot f_n)_{n=1}^\infty$ is truncating for $\la_1(a_1) \ot 
\la_2(a_2)$, the remarks preceding this proposition imply that the sequence
$(\,(\la_1 \ot \la_2)(a_1 \ot a_2)\bigl( \,(e_n \ot f_n)\, x 
\bigr)\,)_{n=1}^\infty$ converges to $(\la_1(a_1) \ot \la_2(a_2))(x)$.

Evidently, we have also that $(\, (e_n \ot f_n)\,x )_{n=1}^\infty$ converges 
to $x$. Therefore, the closedness of $(\la_1 \ot \la_2)(a_1 \ot a_2)$ implies 
that $x$ belongs to $D((\la_1 \ot \la_2)(a_1 \ot a_2))$
and $(\la_1 \ot \la_2)(a_1 \ot a_2) \, x = (\la_1(a_1) \ot \la_2(a_2)) \, x$.
\end{demo}

\begin{corollary}
Consider $a_1 \in \hat{\cM}_{\vfi_1}^+$ and $a_2 \in \hat{\cM}_{\vfi_2}^+$. 
Then $a_1 \ot a_2$ belongs to $\hat{\cM}_{\vfi_1 \ot \vfi_2}^+$ and \newline
$(\vfi_1 \ot \vfi_2)(a_1 \ot a_2) = \vfi_1(a_1) \ot \vfi_2(a_2)$.
\end{corollary}

\section{Appendix 1 : Miscellaneous results.}

In this appendix, we prove some general results, which will be used
several times in this paper.

\begin{lemma} \label{lemA7} Consider a Hilbert \cst-module $E$ over a \cst-
algebra $A$. Let $T$ be a positive element in $\cL(E)$. Then we have that 
$\|T v\|^2 \leq \|T\| \, \|\langle T v,v \rangle \|$.
\end{lemma}

The proof of this lemma is very simple because
$$\langle T v , T v \rangle = \langle  T^2 v , v \rangle \leq \|T\| \, 
\langle T v , v \rangle , $$
where we used the fact that $T^2 \leq \|T\| \, T$.

We will apply this little result in two situations.

\begin{lemma}  \label{lemA1} Consider a Hilbert \cst-module $E$ over a \cst-
algebra $A$.
Let $(T_i)_{i \in I}$ be a net in $\cL(E)^+$ and $T$ an element in $\cL(E)^+$ 
such that $T_i \leq T$ for every $i \in I$.
Then $(T_i)_{i \in I}$ converges strongly to $T$ if and only if
$(\langle T_i v , v \rangle)_{i \in I}$ converges to $\langle T v, v \rangle$ 
for every $v \in E$.
\end{lemma}
\begin{demo}
We have for every $i \in I$ that $\|T_i\| \leq \|T\|$.

Using the previous lemma, we get for every $i \in I$ and $v \in E$ that
$$\|T v - T_i v\|^2 \leq \|T - T_i\| \, \|\langle T v - T_i v , v  \rangle \| 
\leq 2 \|T\| \, \|\langle T v , v \rangle - \langle T_i v , v \rangle \| .$$

Now, the lemma follows easily.
\end{demo}

\begin{lemma} \label{lemA2} Consider a Hilbert \cst-module $E$ over a \cst-
algebra $A$.
Let $(T_i)_{i \in I}$ be an increasing net in $\cL(E)^+$.
Then $(T_i)_{i \in I}$ is strongly convergent in $\cL(E)^+$ if and only if
the net $(\langle T_i v , v \rangle)_{i \in I}$ is convergent for every $v 
\in E$.
\end{lemma}
\begin{demo}
One implication is trivial, we prove the other one.

Therefore, suppose that $(\langle T_i v , v \rangle)_{i \in I}$
is convergent for every $v \in E$.

First, we prove that $(T_i)_{i \in I}$ is bounded.

\begin{list}{}{\setlength{\leftmargin}{.4 cm}}

\item Choose $u \in E$. Because $(\langle T_i u ,u \rangle)_{i \in I}$ is 
convergent, there exist a positive number $N_u$ and an element $i_0$ such 
that $\| \langle T_i u , u \rangle \| \leq N_u$ for ever $i \in I$ with $i 
\geq i_0$.

For every $j \in I$ there exist an element $i \in I$ with $i \geq i_0$
and $i \geq j$, implying that
$$ \| \langle T_j u , u \rangle \| \leq \| \langle T_i u , u \rangle \| \leq 
N_u .$$
So we get that the net $(\langle T_i u , u \rangle)_{i \in I}$ is bounded.

\medskip

Let us fix $w \in E$. By polarisation and the previous result, we have for 
every $v \in E$ that the net $(\langle T_i w , v \rangle)_{i \in I}$ is 
bounded.
Using the uniform boundedness principle, we get that the net
$(T_i w)_{i \in I}$ is bounded.

\medskip

Applying the uniform boundedness principle once again, we see that the net 
$(T_i)_{i \in I}$ is bounded. Hence there exist a strictly positive
number $M$ such that $\|T_i\| \leq M$ for every $i \in I$.

\end{list}

Choose $u \in E$. Take $\vep > 0$. Then there exist an element $i_1 \in I$ 
such that $\|\langle T_i u , u \rangle - \langle T_{i_1} u , u \rangle \|  
\leq \frac{\vep}{2M}$ for every $i \in I$ with $i \geq i_1$. Using lemma 
\ref{lemA7}, we have for every $i \in I$ with $i \geq i_1$ that
$$\| T_i u - T_{i_1} u \|  \leq  \|T_i - T_{i_1}\| \, \| \langle (T_i - 
T_{i_1}) u , u \rangle \|
 \leq  2 M \, \frac{\vep}{2M} = \vep .$$

So we see that $(T_i u)_{i \in I}$ is Cauchy and hence convergent in $E$.

\medskip

From this all, we infer the existence of a mapping $T$ from $E$ into $E$ such 
that $(T_i(w))_{i \in I}$ converges to $T(w)$ for every $w \in E$. It follows 
immediately that $\langle T v , w \rangle = \langle v , T w \rangle$ for 
every $v,w \in E$. This implies that $T$ belongs to $\cL(E)$
and $T^* = T$. Moreover, it  is also clear that $\langle T v , v \rangle \geq 
0$ for every $v \in E$, which implies that $T \geq 0$.
\end{demo}

We want to restate this results in the \cst-algebra case.

\begin{lemma} \label{lemA3}
Consider a \cst-algebra $A$. Let $(a_i)_{i \in I}$ be a net in $M(A)$ and $a$ 
an element in $M(A)$ such that $a_i \leq a$ for every $i \in I$.
Then $(a_i)_{i \in I}$ converges strictly to $a$ if and only if
$(b^* a_i b)_{i \in I}$ converges to $b^* a b$ for every $b \in A$.
\end{lemma}

\begin{lemma} \label{lemA4}
Consider a \cst-algebra $A$. Let $(a_i)_{i \in I}$ be an increasing net in 
$M(A)^+$. Then $(a_i)_{i \in I}$ is strictly convergent in $M(A)^+$
if and only if the net $(b^* a_i b)_{i \in I}$ is convergent for every $b  
\in A$.
\end{lemma}

The following lemma is due to Jan verding (see lemma A.1.2 of \cite{Verd}).

\begin{lemma} \label{lemA5}
Let $E$ be a normed space, $H$ a Hilbert space and $\la$ linear mapping
from within $E$ into $H$. Let $(x_i)_{i \in I}$ be a net in $D(\la)$ and $x$ 
an element in $E$ such that $(x_i)_{i\in I}$ converges to $x$ and
$(\la(x_i))_{i\in I}$ is bounded. Then there exists a sequence 
$(y_n)_{n=1}^\infty$ in the convex hull of $\{ x_i \mid i \in I \}$
and an element $v \in H$ such that $(y_n)_{n=1}^\infty$ converges to $y$
and $(\la(y_n))_{n=1}^\infty$ converges to $v$.
\end{lemma}
\begin{demo}
By the Banach-Alaoglu theorem, there exists a subnet $(x_{i_j})_{j\in J}$ of 
$(x_i)_{i\in I}$ and $v \in H$ such that
$(\la(x_{i_j}))_{j \in J}$ converges to $v$ in the weak topology on $H$. (For 
this, we need $H$ to be a Hilbert space.)

\medskip

Fix $n \in \N$. Then there exists $j_n \in J$ such that $\|x_{i_j} - x\| \leq 
\frac{1}{n}$ for all $j \in J$ with $j \geq j_n$.

Now $v$ belongs to the weak-closed convex hull of the set $\{\la(x_{i_j}) 
\mid j \in J \text{ such that } j \geq j_n \}$, which is the same as  the 
norm-closed convex hull.

Therefore, there exist $\lambda_1,\ldots\!,\lambda_m \in \R^+$
with $\sum_{k=1}^m \lambda_k = 1$ and elements $\alpha_1,\ldots\!,\alpha_m 
\in J$  with $\alpha_1,\ldots\!,\alpha_m \geq j_n$  such that
$$ \| v - \sum_{k=1}^m \lambda_k \la(x_{i_{\alpha_k}}) \| \leq \frac{1}{n} .$$

Put $y_n =  \sum_{k=1}^{l} \lambda_k x_{i_{\alpha_k}}$. Then
$y_n \in D(\la)$, and $\la(y_n) = \sum_{k=1}^m \lambda_k 
\la(x_{i_{\alpha_k}})$.

Therefore, we have immediately that $\|v - \la(y_n)\| \leq \frac{1}{n}$.

Furthermore,
$$\|x - y_n\| = \left\|\sum_{k=1}^{m}\lambda_k(x - x_{i_{\alpha_k}})\right\| 
\leq
\sum_{k=1}^{m} \lambda_k  \frac{1}{n} = \frac{1}{n}$$

\medskip

Therefore, we find that $(y_n)_{n=1}^\infty$ converges to $y$ and
that $(\la(y_n))_{n=1}^\infty$ converges to $v$.
\end{demo}

\begin{lemma} \label{lemA6} Consider a \cst-algebra $A$ and a dense left 
ideal $N$ in $A$. Let $s$ be a positive sesquilinear form on $N$ such that 
$s(a b_1,b_2) = s(b_1,a^* b_2)$ for all $a\in A$ and all $b_1,b_2\in N$.
Moreover, suppose that there exists a positive linear functional $\th$ on $A$ 
such $s(b,b)\leq \th(b^*b)$ for every $b \in N$.

Then there exists a unique positive linear functional $\omega$ on $A$  with 
$\om \leq \th$  such that
$\omega(b_2^* \, b_1) = s(b_1,b_2)$ for every $b_1,b_2 \in N$.
\end{lemma}

\begin{demo}
Let $(\pi,H,v)$ be GNS-object for $\th$ ($v$ is a cyclic vector).

Because $s$ is a positive sesquilinear form on N, we can use the Cauchy-
Schwarz inequality for $s$. So we have for every $b_1,b_2\in N$ that

$$ |s(b_1,b_2)|^2 \leq  s(b_1,b_1)\,\,s(b_2,b_2)
\leq  \th(b_1^*\,b_1)\,\,\th(b_2^*\,b_2)
=  \|\pi(b_1) v\|^2\,\, \|\pi(b_2) v\|^2  .$$
Therefore we can define a continuous positive sesquilinear form $t$
on $H$ such that
$t(\pi(b_1) v,\pi(b_2)v) = s(b_1,b_2)$ for every $b_1,b_2 \in N$.
It is clear that $t$ is positive and $\|t\| \leq 1$.

So there exist an element $T \in B(H)$ with $0\leq T \leq 1$ such that
$t(x,y) = \langle T x ,y \rangle$ for every $x,y \in H$.

This implies that $\langle T \pi(b_1) v , \pi(b_2) v \rangle = s(b_1,b_2)$ 
for every $b_1,b_2 \in N$.

\medskip

Next we show that $T$ belongs to $\pi(A)'$. Therefore, choose $a \in N$.

We have for every $b_1,b_2 \in N$ that
\begin{eqnarray*}
\langle T \pi(a)\, \pi(b_1)v ,\pi(b_2)v \rangle & = &
\langle T \pi(a b_1)v ,\pi(b_2)v \rangle
= t(\pi(a b_1)v ,\pi(b_2) v) \\
& = & s(a b_1,b_2) =  s(b_1,a^* b_2) \\
& = & t(\pi(b_1)v,\pi(a^* b_2) v) =  \langle T \pi(b_1) v,\pi(a^* b_2)v 
\rangle \\
& = & \langle T \pi(b_1)v , \pi(a^*)\pi(b_2)v \rangle
=  \langle \pi(a) T \pi(b_1) v , \pi(b_2)v \rangle .
\end{eqnarray*}
This implies that $T \pi(a) = \pi(a) T$.

\medskip

Now we define the continuous linear functional $\omega$ on $A$ such
$\omega(x) = \langle T \pi(x)v , v \rangle $ for every $x \in A$.

Using the fact that $T$ belongs to $\pi(A)'$, we have for every $b_1,b_2 \in 
N$ that $$\om(b_2^*\,b_1) = \langle T \pi(b_2^* \, b_1) v , v \rangle
= \langle T \pi(b_1) v , \pi(b_2) v \rangle = s(b_1,b_2). $$

Consequently, we have for every  $b \in N$ that $\om(b^* \, b) = s(b,b)$,
implying that $0 \leq \om(b^* b) \leq \th(b^* b)$.
This implies easily that $0 \leq \om \leq \th$.

\end{demo}

We have even proven a stronger result where $N$ is not assumed to be dense 
(of course the unicity is not longer valid in this case).
This proof can be found in lemma A.1.3 of \cite{Verd}.

\section{Appendix 2 : A small technical result.}
\label{art4}

In this appendix, we will prove a technical result which will be used
in several sections.

Consider a Hilbert \cst-module $E$ over a \cst-algebra $B$
and let $D$ be subset of $E$ such that its linear span is dense in $E$.
Let $(T_i)_{i \in I}$ be a net in $\cL(E)$ such that $(T_i)_{i \in I}$ 
converges strongly$^*$ to 1.

Suppose that $t$ is a mapping from $B$ into $E$ such that
for every $i \in I$ and every $v \in D$ there exists an element
$x(v,i) \in B$ such that $\langle T_i \, t(b) , v \rangle  = x(v,i) \, b$ for 
every $b \in B$.

We want to prove that $t$ belongs to $\cL(B,E)$.

\begin{lemma}
We have that $t$ is a continuous $B$-linear map from $B$ into $E$.
\end{lemma}
\begin{demo}
\begin{itemize}
\item Choose $b_1,b_2 \in B$.

Fix $j \in I$. We have for every $v \in D$ that
\begin{eqnarray*}
\langle T_j \, t(b_1 b_2) , v \rangle & = & x(v,j) \, (b_1 b_2)
= (x(v,j) b_1) b_2 = \langle T_j \, t(b_1) , v \rangle \, b_2 \\
& = & \langle T_j \, ( t(b_1) b_2 ) , v \rangle .
\end{eqnarray*}
Because the linear span of $D$ is dense in $E$, this implies
that $T_j t(b_1 b_2) = T_j (t(b_1) b_2)$.

Because $(T_i)_{i \in I}$ converges strongly to 1, this implies that
$t(b_1 b_2) = t(b_1) b_2$.

\item The linearity of $t$ is proven in a completely similar way.

\item Choose a sequence $(b_n)_{n=1}^\infty$ in $B$, $b \in B$ and $w \in E$ 
such that $(b_n)_{n=1}^\infty \rightarrow b$ and
$(t(b_n))_{n=1}^\infty \rightarrow w$.

Fix $j \in I$. Take $v \in D$.

We have for every $n \in \N$ that $\langle T_j \, t(b_n) , v \rangle
= x(v,j)\, b_n$. This implies that $(\langle T_j \, t(b_n) , v 
\rangle)_{n=1}^\infty$ converges to $x(v,j)\, b$, which is equal to
$\langle T_j \, t(b) , v \rangle$.
It is also clear that $(\langle T_j \, t(b_n) , v \rangle)_{n=1}^\infty$
converges to $\langle T_j \, w , v \rangle$.

These two results imply that $\langle T_j \, t(b) , v \rangle = \langle T_j 
w, v \rangle$.

Of course, this implies that $t(b) = w$.

Therefore, we have proven that $t$ is closed. The closed graph theorem 
implies that $t$ is continuous.

\end{itemize}
\end{demo}

\begin{lemma}
Consider an approximate unit $(e_k)_{k \in K}$ for $B$.
For every $k \in K$, we define the element $S_k \in \cL(B,E)$ such that
$S_k(b) = t(e_k) \, b$ for every $b \in B$.

Then we have for every $v \in E$ that $(S_k^*(v))_{k \in K}$ is convergent in 
$B$.
\end{lemma}
\begin{demo}
Remember from the previous lemma that $t$ is a continuous B-linear map
from $B$ into $E$. In particular, we have for every $k$ in $K$ that $\|S_k\| 
\leq \|t\|$.

Choose $w \in D$. Take $\vep > 0$.

Then there exists an element $j \in I$ such that $\|T_j^* w - w\| \leq 
\frac{\vep}{3} \frac{1}{1+\|t\|}$.

We have for every $k \in K$ that
$$S_k^* (T_j^* w) = \langle T_j^* w , t(e_k) \rangle = \langle w , T_j \, 
t(e_k) \rangle = \langle T_j \, t(e_k) , w \rangle ^*
= (x(j,w) e_k)^* = e_k \, x(j,w)^* ,$$
which implies that $(S_k^*(T_j^* w))_{k \in K}$ converges to $ x(j,w)^*$.
Consequently, there exist an element $k_0 \in K$ such that
$\| S_{k_1}^*(T_j^* w) - S_{k_2}^*(T_j^* w) \| \leq \frac{\vep}{3}$
for every $k_1,k_2 \in K$ with $k_1,k_2 \geq k_0$.

Therefore, we have for every $k_1,k_2 \in K$ with $k_1,k_2 \geq k_0$ that
\begin{eqnarray*}
\| S_{k_1}^*(w) - S_{k_2}^*(w) \|
& \leq & \| S_{k_1}^*(w) - S_{k_1}^*(T_j^* w) \|
+ \| S_{k_1}^*(T_j^* w) - S_{k_2}^*(T_j^* w) \| \\
& & + \| S_{k_2}^*(T_j^* w) - S_{k_2}^*(w) \| \\
& \leq & \|S_{k_1}^*\| \, \|w - T_j^* w \|
+ \frac{\vep}{3} + \|S_{k_2}^*\| \, \|w - T_j^* w \| \\
& \leq & \|t\| \, \frac{\vep}{3} \frac{1}{1+ \|t\|} + \frac{\vep}{3} +
\|t\| \, \frac{\vep}{3} \frac{1}{1+ \|t\|} \leq \vep .
\end{eqnarray*}

Hence, we see that $(S_k^*(w))_{k \in K}$ is Cauchy and hence convergent in 
$B$.

From this, we get immediately that $(S_k^*(v))_{k \in K}$ is convergent for 
every $v$ in the linear span of $D$.

Because this linear span is dense in $E$ and $(S_k^*)_{k \in K}$ is bounded, 
we get that $(S_k^*(v))_{k \in K}$ is convergent for every
$v \in E$.
\end{demo}

\medskip

Now we can formulate the final result.

\begin{proposition} \label{prop4.1} We have that $t$ belongs to $\cL(B,E)$.
\end{proposition}
\begin{demo}
Take an approximate unit $(e_k)_{k \in K}$ for $B$.

For every $k \in K$, we define the element $S_k \in \cL(B,E)$ such
that $S_k(b) = t(e_k) \, b $ for every $b \in B$.

Choose $c \in B$. Because $t$ is $B$-linear, we have for ever $k \in K$
that $S_k(c) = t(e_k) c = t(e_k c)$. The continuity of $t$ implies that
$(S_k(c))_{k \in K}$ converges to $t(c)$.

By the previous lemma, we know that there exist a linear mapping $r$ from
$E$ into $B$ such that $(S_k^*(v))_{k \in K}$ converges to $r(v)$ for every 
$v \in E$.

Combining these two results, we get that $\langle t(b) , v \rangle
= \langle b , r(v) \rangle$ for every $b \in B$ and $v \in E$.
This implies that $t$ belongs to $\cL(B,E)$.
\end{demo}

\end{document}